\documentclass[fullpage]{article}

\usepackage{times}
\usepackage{xspace}
\usepackage[breaklinks,colorlinks,
	    linkcolor=blue,citecolor=blue,urlcolor=black]{hyperref}
\usepackage{xurl}
\usepackage{amsmath,amssymb}
\usepackage{stmaryrd}
\usepackage[noabbrev,capitalize]{cleveref}
\usepackage{xstring}
\usepackage{graphicx}
\usepackage{subcaption}
\usepackage{nicematrix}		% baseline=b option for inference rule premises

\hypersetup{breaklinks=true}

\newcommand{\com}[1]{}

\newcommand{\later}[1]{}
\long\def\note#1{}	% exclude
\long\def\baf#1{\note{{\bf BAF: } [{\color{blue} \em #1}]}}	% Bryan's notes
\newcommand{\ie}{i.e.,\xspace}
\newcommand{\eg}{e.g.,\xspace}

\newcommand{\gd}{GD\xspace}		% grounded deduction
\newcommand{\ga}{GA\xspace}		% grounded arithmetic
\newcommand{\bga}{BGA\xspace}		% basic grounded arithmetic
\newcommand{\cga}{CGA\xspace}		% constructive grounded arithmetic
\newcommand{\rga}{RGA\xspace}		% reflective grounded arithmetic
\newcommand{\pcf}{PCF\xspace}		% programming computable functions
\newcommand{\ppf}{PPF\xspace}		% programming parallel functions
\newcommand{\pra}{PRA\xspace}		% primitive recursive arithmetic
\newcommand{\pa}{PA\xspace}		% Peano arithmetic

\newcommand{\gdl}{G\"odel\xspace}

\newcommand{\limp}{\rightarrow}		% logical implication
\newcommand{\liff}{\leftrightarrow}	% logical biconditional (iff)
\newcommand{\leqv}{\equiv}		% logical equivalence (not definition)
\newcommand{\ldef}{\equiv}		% logical equivalence by definition
\newcommand{\subs}[3]{[{#1}/{#2}]{#3}}
\newcommand{\ttmore}{\dots}		% notation for more unnamed variables
\newcommand{\ttc}[2]{{#1}\langle{#2}\rangle}	% term #1 with only vars #2
\newcommand{\tto}[2]{\ttc{#1}{#2,\ttmore}} 	% #1 with at least vars #2

\newcommand{\kwstyle}[1]{\textbf{#1}\xspace}
\newcommand{\kcase}{\kwstyle{case}}
\newcommand{\kof}{\kwstyle{of}}
\newcommand{\kif}{\kwstyle{if}}
\newcommand{\kifz}{\kwstyle{ifz}}
\newcommand{\klet}{\kwstyle{let}}
\newcommand{\kletrec}{\kwstyle{letrec}}
\newcommand{\kin}{\kwstyle{in}}
\newcommand{\ky}{\kwstyle{Y}}		% the Y combinator for higher-order recursion
\newcommand{\kisS}{\kwstyle{isS}}
\newcommand{\kisP}{\kwstyle{isP}}
\newcommand{\kisplus}{\kwstyle{is$+$}}
\newcommand{\kisminus}{\kwstyle{is$-$}}
\newcommand{\kismult}{\kwstyle{is$\times$}}

\newcommand{\typestyle}[1]{\textsf{#1}\xspace}
\com{	% types-as-identifiers style similar to Martin-Löf style
\newcommand{\tbool}{\typestyle{bool}}
\newcommand{\ttrue}{\typestyle{true}}
\newcommand{\tfalse}{\typestyle{false}}
\newcommand{\tobj}{\typestyle{obj}}
\newcommand{\tnat}{\typestyle{nat}}

}%com
\newcommand{\tbool}{\typestyle{B}}	% "boolean" or "binary" or "bit"
\newcommand{\ttrue}{\typestyle{T}}
\newcommand{\tfalse}{\typestyle{F}}
\newcommand{\tobj}{\typestyle{O}}
\newcommand{\tnat}{\typestyle{N}}

\newcommand{\conststyle}[1]{\textsf{#1}\xspace}
\newcommand{\ctrue}{\ttrue}
\newcommand{\cfalse}{\tfalse}
\newcommand{\cneg}{\conststyle{neg}}
\newcommand{\cand}{\conststyle{and}}
\newcommand{\cforall}{\conststyle{forall}}
\newcommand{\ceq}{\conststyle{eq}}

\newcommand{\czero}{\conststyle{zero}}
\newcommand{\csuc}{\conststyle{suc}}
\newcommand{\cvar}{\conststyle{var}}

\newcommand{\judgment}[1]{\textsf{ #1}}
\newcommand{\jtrue}{\judgment{\ttrue}}
\newcommand{\jfalse}{\judgment{\tfalse}}
\newcommand{\jbool}{\judgment{\tbool}}
\newcommand{\jobj}{\judgment{\tobj}}

\newcommand{\jnat}{\judgment{\tnat}}

\newcommand{\jtype}{\judgment{type}}

\newcommand{\funstyle}[1]{\textsf{#1}\xspace}
\newcommand{\fnat}{\funstyle{nat}}
\newcommand{\fbool}{\funstyle{bool}}
\newcommand{\feven}{\funstyle{even}}

\newcommand{\sbool}{\mathbb{B}}		% boolean set
\newcommand{\snat}{\mathbb{N}}		% natural numbers
\newcommand{\sreal}{\mathbb{R}}		% real numbers
\newcommand{\fprv}{\mathbf{Pr}}		% provably true (letter notation)
\newcommand{\funp}{\mathbf{U}}		% unprovable (neither true nor false)
\newcommand{\fcon}{\mathbf{Con}}	% consistent
\newcommand{\oq}{\mathrel{?}}		% guard conditional

\makeatletter
\newlength{\doublefracgap}
\DeclareRobustCommand{\doublefrac}[2]{%
  \mathinner{\mathpalette\doublefrac@{{#1}{#2}}}%
}
\newcommand{\doublefrac@}[2]{\doublefrac@@#1#2}
\newcommand{\doublefrac@@}[3]{%
  \ooalign{%
    \raisebox{\doublefracgap}{$\m@th#1\frac{#2}{\phantom{#3}}$}\cr
    \raisebox{-\doublefracgap}{$\m@th#1\frac{\phantom{#2}}{#3}$}\cr
  }%
}
\newcommand{\ddoublefrac}[2]{{\displaystyle\doublefrac{#1}{#2}}}

\makeatother

\newcommand{\irl}[1]{\ensuremath{\mathit{#1}}}		% inference rule label
\newcommand{\rlpo}{\irl{\exists TI}\xspace}		% LPO aka A/E TI

\newcommand{\infrule}[3][]{\cfrac{#2}{#3}\IfStrEq{#1}{}{}{\ \irl{#1}}}
\newcommand{\infeqv}[3][]{\ddoublefrac{#2}{#3}\IfStrEq{#1}{}{}{\ \irl{#1}}}
\newcommand{\infceqv}[4][]{\cfrac{#2\qquad}{}\ddoublefrac{#3}{#4}
				\IfStrEq{#1}{}{}{\ \irl{#1}}}
\newcommand{\suc}{\mathbf S}
\newcommand{\pred}{\mathbf P}
\newcommand{\quo}[1]{\ulcorner{#1}\urcorner}
\newcommand{\unq}[1]{\llcorner{#1}\lrcorner}

\newcommand{\qb}[1]{\llbracket{#1}\rrbracket}	% quote-bracketed anything
\newcommand{\qbe}[2]{\qb{{#1}\vdash{#2}}}	% quote-bracketed entailment

\newcommand{\fpr}{\square}			% provably true in box notation
\newcommand{\tpr}[2][]{{#1}\mathbin{\fpr}\qb{#2}} %  provability of quoted term

\newcommand{\tforall}[2]{\forall{#1}\ {#2}}
\newcommand{\texists}[2]{\exists{#1}\ {#2}}
\newcommand{\tforallp}[2]{\forall^{+}{#1}\ {#2}}
\newcommand{\texistsp}[2]{\exists^{+}{#1}\ {#2}}
\newcommand{\texuniq}[2]{\exists!{#1}\ {#2}}
\newcommand{\tmu}[2]{\mu{#1}\ {#2}}		% bounded or unbounded search

\newcommand{\tif}[3]{\kwstyle{if}\ {#1}\ \kwstyle{then}\ {#2}\ \kwstyle{else}\ {#3}}
\newcommand{\tifz}[3]{\kifz({#1},{#2},{#3})}
\newcommand{\tcase}[2]{\hbox{\kcase\ }{#1}\hbox{ \kof\ }{#2}}
\newcommand{\tcasei}[2]{{#1}\Rightarrow{#2}}
\newcommand{\tlet}[2]{\hbox{\klet\ }{#1}\hbox{ \kin\ }{#2}}
\newcommand{\tletrec}[2]{\kletrec\ {#1}\hbox{ \kin\ }{#2}}
\newcommand{\tofn}[1]{\underline{#1}}
\newcommand{\tlambda}[2]{\lambda{#1}\ {#2}}

\newcommand{\lift}[1]{\underline{#1}}		% lift into a domain element
\newcommand{\ltrue}{\lift{\ctrue}}
\newcommand{\lfalse}{\lift{\cfalse}}
\title{Reasoning Around Paradox \\ with Grounded Deduction}
\author{Bryan Ford}
\date{First version: September 12, 2024 \\
	This version: \today}

\begin{document}

\maketitle
\begin{abstract}
How can we reason around logical paradoxes without falling into them?
This paper introduces \emph{grounded deduction} or \gd,
a Kripke-inspired approach to first-order logic and arithmetic
that is neither classical nor intuitionistic,
but nevertheless appears both pragmatically usable
and intuitively justifiable.
\gd permits the direct expression of unrestricted recursive definitions --
including paradoxical ones such as `$L \ldef \neg L$' --
while adding \emph{dynamic typing} premises to certain inference rules
so that such paradoxes do not lead to inconsistency.
This paper constitutes a preliminary development and investigation
of grounded deduction,
to be extended with further elaboration
and deeper analysis of its intriguing properties.
\end{abstract}

\newpage
\tableofcontents

\newpage
\section{Introduction}
\label{sec:intro}

How well-founded are the classical rules of logical deduction
that we normally rely on throughout mathematics and the sciences?
This topic has been debated for centuries.

\subsection{Pythagoras visits Epimenides}

Let us indulge briefly in an anachronistic reimagining 
of what transpired when Pythagoras met Epimenides in Crete.
Upon hearing Epimenides utter the phrase ``Cretans, always liars''
as part of his ode to Zeus,
Pythagoras becomes troubled wondering
whether Epimenides, a Cretan, was lying.
Seeking answers,
Epimenides takes Pythagoras to the oracle in the cave of Ida,
known always to speak the truth.
Pythagoras asks the oracle:

\begin{quote}
%O great oracle,
%I ask that we name my next sentence $L$. \\
%This sentence is false. \\
%O oracle,
%is my sentence $L$ true or false?
O oracle,
I ask only the following:
is your answer to my question ``no''?
\end{quote}

Reports differ on what ensued next.
By one account,
the oracle emitted a deafening shriek
and vanished in a cloud of acrid smoke.
Pythagoras hastily fled the island,
fearing retribution once the Cretans learned they had lost their oracle.
%For reasons to be clarified later,
%let us call this the \emph{classical account}.

By conflicting reports,
however,
the oracle merely stared back at Pythagoras and told him calmly:
``Your question is circular bullshit.''
Pythagoras departed the island in shame,
never to mention the incident
or leave its record in the history books.
%Let us call this the \emph{grounded account}.

Pythagoras's query above is of course
just a variation on the well-known \emph{Liar paradox},
related to though distinct from the \emph{Epimenides paradox}
that later became associated
with Epimenides' famous line of poetry.\footnote{
	A Cretan's claim that Cretans are ``always liars''
	is of course technically paradoxical
	only under dubious semantic assumptions,
	such as that Epimenides meant
	that \emph{all} Cretans \emph{always} lie and never tell the truth.
	In fact Epimenides' line
	was probably not meant to be paradoxical at all,
	but was rather a religious reaction to an impious belief
	that Zeus was not living as a deity on Mount Olympus
	but was dead and buried in a tomb on Crete;
	see \cite{harris06cretans}.
	For a broader history of the Liar and other paradoxes,
	see \cite{sorensen05brief}.
}

Let us focus, however,
on the two conflicting accounts above of the oracle's response to Pythagoras.
In the first, which we'll call the \emph{classical account},
the oracle self-destructs trying to answer the question,
as in any number of science-fiction scenarios where
the hero triumphs over an evil computer or artificial intelligence
by giving it some problem ``too hard to solve.''\footnote{
	The 19883 film \emph{WarGames} comes to mind
	as a classic Hollywood example.
}
In the second account, which we'll call the \emph{grounded account},
the oracle simply recognizes
the circular reasoning in the Liar Paradox for what it is,
and calls bullshit on the question
instead of trying to answer it.\footnote{
	We use the term ``bullshit'' here not as an expletive
	but as a technical term embodying an important semantic distinction
	from mere falsehood.
	Whereas a lie deliberately misrepresents some known truth,
	bullshit does not care what the truth is,
	or even if there is any relevant truth.
	In the words of \cite{frankfurt05bullshit}:
	\begin{quote}
		It is impossible for someone to lie unless he thinks he knows the truth.
		Producing bullshit requires no such conviction.
		\com{
		$\dots$
		For the bullshitter, $\dots$ all these bets are off:
		he is neither on the side of the true
		nor on the side of the false.
		His eye is not on the facts at all,
		as the eyes of the honest man and of the liar are,
		except insofar as they may be pertinent to his interest
		in getting away with what he says.
		}%com
	\end{quote}
}
When a paradox like this clearly causes something to go wrong,
where does the blame lie:
with the oracle asked to answer the question,
or with the question itself?

\subsection{The paradoxes in classical and alternative logics}

In developing mathematics and computer science
atop the accepted foundation of classical logic,
we must carefully guard our formal systems from numerous paradoxes
like that above.
Avoiding paradoxes impels us
to forbid unconstrained recursive definitions,
for example,
where a new symbol being defined also appears part of its definition.
Allowing unconstrained recursive definitions in classical logic
would make the Liar paradox trivially definable as `$L \equiv \neg L$',
leading to immediate inconsistency.
$L$ becomes provably both true and false,
and subsequently so do all other statements,
rendering the logic useless for purposes
of distinguishing truth from falsehood.

\later{		looks like we still need type stratification in set reasoning
Avoiding paradoxes led to Russell's stratified theory of types~\cite{XXX}
and its many descendants.
The paradoxes
similarly led set theory to the standard principle that ``larger'' sets
may be constructed only from ``smaller'' sets
in rigorously-constrained ways~\cite{XXX}.
}

Understandably dissatisfied with this apparent fragility,
alternative philosophical schools of thought have explored
numerous ways to make logic or mathematics more robust
by weakening the axioms and/or deduction rules that we use.\footnote{
	For a broad and detailed exploration
	of many such alternative approaches
	to the problems of truth and paradox,
	see for example \cite{field08saving}.
\later{
as in the traditions of constructivism~\cite{XXX},
ultrafinitism~\cite{XXX},
substructural logics~\cite{XXX},
	% see https://plato.stanford.edu/entries/curry-paradox/
	% and https://plato.stanford.edu/entries/logic-substructural/
etc.
}
}
Most of these alternative formulations of logic
leave us pondering two important questions, however.
First,
could we envision actually \emph{working in} such an alternative logic,
carrying out what we recognize
as more-or-less normal mathematics or computer science --
and how would such adoption affect (or not) our everyday reasoning?
Second,
since most of these alternative logics ask us to live with
unfamiliar and often counterintuitive new constraints on our reasoning,
what is the \emph{payoff} for going to this trouble?
What ideally-useful benefit would we get, if any,
for accepting unfamiliar constraints on our basic deduction methods --
for seemingly ``tying our hands''?
The latter question
may be central to the reason that most alternative logics along these lines
remain obscure curiosities
of great interest to experts specializing in formal logic,
but to few others.

\subsection{Introducing grounded deduction (\gd)}
\label{sec:intro:gd}

This paper presents \emph{grounded deduction} or \gd,
a foundation for logical deduction
that attempts to avoid classical logic's difficulties
with the traditional paradoxes,
while striving at a framework
in which we might plausibly hope to do normal work
in mathematics or the sciences
without inordinate or unjustified difficulty.\footnote{
	The term ``grounded deduction''
	is inspired by the notion of a statement
	being \emph{grounded} or not
	in Kripke's theory of truth,
	one important precedent for this work
	along with many others
	(see \cite{kripke75outline}).
}
Most importantly,
\gd endeavors to
\emph{offer something in return}
for the strange and perhaps uncomfortable new constraints
it imposes on our traditional methods of deduction.

The main immediate ``payoff'' that \gd offers
is the permission to make unconstrained recursive definitions.
That is, \gd allows definitions of the form `$s \equiv d$',
which may include the newly-defined symbol $s$
arbitrarily within the expansion $d$ on the right-hand side,
without the usual restrictions
(such as that $d$ be structurally primitive-recursive,
or well-founded by some other criteria).
In particular, \gd permits the direct definition 
of outright paradoxical propositions such as $L \equiv \neg L$
(the Liar paradox),
without apparent inconsistency.
More pragmatically,
\gd's admission of unrestricted recursive definitions
proves useful in concisely expressing and reasoning about
numerous standard concepts in working mathematics and computer science.

\begin{figure}[t]
\begin{center}
\includegraphics[width=0.5\textwidth]{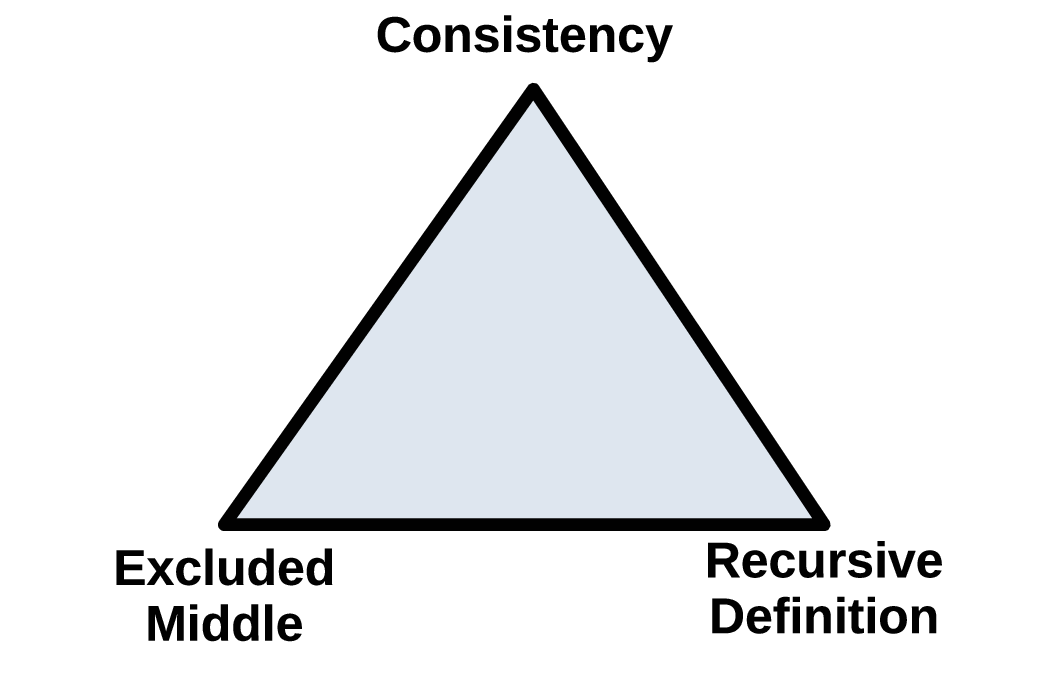}
\end{center}
\caption{An impossibility triangle for logic.
	We may desire our logical reasoning to be (fully) consistent,
	to give us (unrestricted) use of the law of excludded middle or LEM,
	and to give us (unrestricted) recursive definition capability.
	It appears 
	we must compromise at least one of these desires, however.}
\label{fig:triangle}
\end{figure}

In the tradition of so-called ``impossibility triangles,''
\cref{fig:triangle}
shows one such triangle that appears to apply to systems of logic.
Of three desirable properties we might wish for -- namely
(full) consistency,
the (unrestricted) law of excluded middle or LEM, and
(unrestricted) recursive definition,
it appears we must compromise and accept
a weakened version of at least one of these properties.
Classical reasoning prioritizes full consistency and LEM
while restricting recursive definitions,
while in \gd we will prioritize consistency and recursive definition
but weaken our LEM.\footnote{
	The third alternative is of course possible too:
	\emph{paraconsistent} logics weaken our demand for logical consistency,
	typically attempting to ``reduce the damage'' caused by inconsistency
	instead of eliminating it entirely.
	See for example \cite{field08saving}
	for a detailed overview of such approaches.
	\later{more citations}
}

Beyond the immediate offering of unrestricted recursive definitions,
the many indirect implications of \gd's alternative perspective
on deduction rules and logical truth
are interesting,
in ways that this paper attempts to begin mapping,
but on which it can admittedly only scratch the surface.
\baf{ summarize a few interesting observations) }

The cost of this flexibility manifests in \gd's deduction rules,
many of which modify the rules of classical logic
by incorporating \emph{typing requirements} into their premises.\footnote{
	\gd's notion of typing
	is heavily influenced by concepts and notations
	that have become ubiquitous in programming language theory and practice,
	such as Martin-L\"of's
	intuitionistic type theory
	as described in \cite{martin-lof80intuitionistic}.
	\gd's logic is not intuitionistic, however, as we will see.
	Further,
	\gd's use of typing is unlike those of
	statically-typed programming languages
	or stratified logics in the tradition of Russell and Tarski,
	\later{citations}
	but rather is more closely analogous
	to dynamically-typed programming languages like Python.
	% like C++, Haskell, or Coq.
	In particular, \gd is syntactically single-sorted,
	having only one syntactic space of \emph{terms}.
	A term's type depends (``dynamically'')
	on whatever value the term actually produces,
	if any --
	whether boolean, integer, set, etc. --
	and not on any stratification
	pre-imposed statically on the syntactic structure.
}
For example, \gd allows us to invoke proof by contradiction --
assuming some proposition $p$ is false hypothetically
in order to prove it true --
only after we \emph{first} prove
that $p$ denotes a well-typed boolean value,
\ie that $p$ is in fact either true or false.
\gd's inference rule for introducing logical implication $p \limp q$
similarly requires us first to prove that $p$ is boolean,
thereby avoiding Curry's paradox,
as we will see in \cref{sec:prop:curry}.

While \gd's typing prerequisites may seem unfamiliar,
we can nonetheless justify them intuitively.
In the grand tradition of abusing Latin for science,
we might say that \gd limits classical logic's \emph{tertium non datur} or LEM
with a counterbalancing principle of \emph{habeas quid}:
we must have a thing before we use it.\footnote{
	\com{
	The phrase \emph{habeas quid}
	is shamelessly incorrect Latin, of course.
	But modern educated audiences
	may at least recognize the composite words
	from the legal principle of \emph{habeas corpus}
	and the economic principle of \emph{quid pro quo}.
	And ``quid'' is shorter than ``aliquid''
	if a bit misused in this context.
	}%com
	The term \emph{habeas quid}
	shamelessly sacrifices linguistic and historical authenticity
	at the altar of expedience.
	While most readers today likely share
	the writer's non-fluency in Latin,
	many will at least recognize the individual words
	from the legal principle of \emph{habeas corpus}
	and the economic principle of \emph{quid pro quo}.
}
\gd differs from classical logic only in
applying \emph{habeas quid} systematically
to all things mathematically expressible,
closing the special loophole
that the LEM traditionally leaves for boolean things.

Incorporating these typing requirements into \gd's basic deduction rules,
fixed at the boundary between \gd
and whatever language or metalogic we use to reason about it,
appears crucial in avoiding so-called \emph{revenge problems},
where solving a pardox in one place
just makes a more subtle but equally-insidious paradox
appear elsewhere.\footnote{
	A collection of essays specifically on such revenge problems in logic
	may be found in \cite{beall08revenge}.
}

Satisfying \gd's typing requirements does impose a ``reasoning cost''
over the familiar rules of classical logic.
In the common case, however,
these typing proofs tend to be trivial
and will likely be subject to future automation
with appropriate tooling.\footnote{
	Such automation might well include static type systems,
	complementing the dynamic type system in \gd's foundation.
	In the same way that static types in (say) TypeScript
	complement the dynamic types native
	to the underlying JavaScript,\later{citations}
	static-typing extensions to \gd could usefully
	both guide and constrain the search space
	that automated reasoning tools must confront,
	while silently discharging most of the tedious typing prerequisites
	that we might have to prove manually in ``raw'' \gd.
	In essence,
	\gd's notion of dynamic typing is meant as a foundational tool
	but by no means is intended as the end of the story.
}
As \gd's goal is to formulate a plausible \emph{working} logic,
the priority is to offer a reasonably complete set
of familiar logical and mathematical tools in the new framework,
ideally comprehensible not just to experts
in mathematical logic or programming language theory,
but also to non-experts.
As a result,
this formulation makes no attempt at minimalism.
Many operators we introduce are definable in terms of others,
and many deduction rules are derivable from combinations of others,
as we note in places.
The author thus offers apologies in advance to experts in logic
to whom this style of formulation and exposition
may feel unbearably verbose, tedious, and often redundant.

\com{
XXX old, massage...

One of the goals of this exposition
is to work out a Kleene/Kripke-style logic
that one might actually be tempted to \emph{use} practically,
and not merely study as a logical or metamathematical curiosity.
A key goal is to get a feel for what it would be like
to inhabit a different logical universe,
to view the mathematical world from another platform
that contrasts significantly from the classical perspective.

}%com

\baf{ summarize current status:
rigorous formal development of \emph{anything} in this paper
remains for future work...
}

This working paper is a draft
that is both preliminary and incomplete.
In particular,
what is presented here is only the first part of a much longer document,
subsequent parts of which will be released in updates to this preprint
as they reach a state of approximate readiness and readability.
There may well be significant gaps or errors in parts already released,
and rigorous formal analysis remains to be done.
The author asks readers to take this current draft, whatever its state,
as a preliminary basis for exploration, discussion, and further development,
and not (yet) as a finished product.

\later{
	\gd is about doing a kind of logical ``due diligence'' on values.
	We need to suck it up and prove we actually have a thing
	before trying to use it as a thing.
}

\section{Propositional deduction in \gd}
\label{sec:prop}

\later{	add false elimination rule}

\later{	make text rules consistent with HOL formalization}

\baf{	Use the term "judgment" consistently and correctly.
	Use "type judgments" for the typed things we have. 
	In \gd a definition is a form of judgment,
	and an entailment is a form of judgment. }

Mirroring the traditional starting point for defining logic,
we first introduce the basic propositional connectives in \gd
for
logical negation ($\neg$),
conjunction ($\land$),
disjunction ($\lor$),
implication ($\limp$),
and
biconditional ($\liff$).
In the process,
we introduce \gd's approach to typing, judgments and deduction.

Classical logic in general,
and the law of the excluded middle (LEM) especially,
presuppose that any syntactically-valid proposition
has an associated truth value of either \ctrue or \cfalse.
Even many-valued logics
such as Kleene's 3-valued logic
typically retain the underlying premise
that every proposition has \emph{some} particular truth value,
while expanding the range of ``choices'' for what that value might be.\footnote{
	Kleene introduced his strong 3-valued logic
	in \cite{kleene38notation},
	%XXX and expanded upon in \cite{kleene52introduction} ???
	as a tool for reasoning about computations that might not terminate
	and their relationship to the ordinal numbers
	of classical set theory.
	The truth-value semantics of conjunction and disjunction
	for grounded deduction as presented here
	line up precisely with those of Kleene's strong 3-valued semantics.
	\gd diverges in other respects, however,
	and we will rely more on modern domain theory
	rather than classical set theory and ordinals
	in order to model and reason about
	the semantics of computation in \gd.
}

\gd starts by rejecting this presumption,
treating it as the ``original sin'' of classical logic.
In \gd,
a proposition by default has \emph{no} value of any type.
In fact,
\gd does not even syntactically distinguish logical propositions
from terms denoting mathematical objects such as integers or sets.
Any syntactic expression is merely an untyped \emph{term} --
until that term is logically \emph{proven}
to represent a value (of some type)
through a ``grounded'' deduction process.
That is,
until and unless we have proven that a term denotes a value of some type,
we refuse to ascribe \emph{any} value or type to that term --
not even a ``third value'' in the usual sense for 3-valued logics.\footnote{
	In this respect,
	\gd bears a close relationship to the paracomplete system KFS
	explored in \cite{field08saving}.
	Field brilliantly characterizes what it means for a formula
	not to have a truth value as follows:

	\begin{quote}
	What then is the relation between truth value and semantic value in KFS?
	In the case of restricted theories
	(which are the only ones for which we have an unrelativized notion
	of semantic value),
	we can say this: 
	having semantic value 1 is sufficient for being true;
	and having semantic value 0 is sufficient for being false
	(\ie having a true negation).
	For sentences with semantic value $1/2$,
	we can't say that they're true, or that they aren't,
	or that they're false, or that they aren't.
	We can't say whether or not they are ``gappy''
	(neither true nor false).
	And our inability to say these things can't be attributed to ignorance,
	for we don't accept that there is a truth about the matter.
	This isn't to say that we think there is no truth about the matter:
	we don't think there is, and we don't think there isn't.
	And we don't think there either is or isn't.  
	Paracompleteness runs deep.
	\baf{	\cite[p.72]{field08saving}}
	\end{quote}
}

\baf{Note somewhere that the formalization
	will be based on Gentzen-style natural deduction,
	based on inference rules rather than axioms.
	This approach facilitates the typing rules
	that \gd introduces.
	Whether and in what way a more axiomatic approach
	to formulating \gd might work is left to future work.}

\baf{	Consider using background colors via colorbox
	to highlight which parts of statements are target language
	and which are metalanguage,
	and (later) at which meta-level when needed.}

\subsection{Boolean truth values}
\label{sec:truth}

\baf{(XXX refine, probably simplify)}

We will typically use the letters $a, b, c$
to represent \emph{terms} in some abstract or concrete syntax.
A term might in principle
represent any type of value (number, set, etc.).
A term might just as well represent no definable value at all,
such as the ``result'' of a paradox or computation that never terminates,
and thus never actually yields any value.

For the present, we do not care
exactly what kinds of values a term $a$ might represent
(if indeed it has a value at all).
Instead,
we care only that there is at least one such expressible value
that we will call ``a true value.''
We also assume
there is at least one expressible value 
that we will call ``a false value.''
A value that falls into either of these categories
we will call a \emph{boolean truth value}.
We will represent the \emph{type} of boolean truth values
via the boldface letter `$\tbool$'.
We will assert that some term $a$ is boolean
by suffixing it with the boolean type letter,
as in `$a \jbool$' -- meaning simply ``$a$ is boolean.''

We require that no value be both true and false,
but we otherwise set no expectations on what these truth values actually are.
We are also agnostic to
how many distinct true values and how many false values might exist.
There might be only one unique true value named \ctrue,
and one distinguished false value named \cfalse,
as in many strongly-typed programming languages.
Alternatively, truth values might be single-digit binary integers,
with 1 as the unique true value, 0 as the unique false value,
and all other numbers not denoting truth values.
We might even take all integers to be truth values,
with 0 as the only false value and all other integers representing true,
as in many weakly-typed languages such as C.
Which particular values might represent true and false
will not concern us here;
we assume merely that such values exist.

Consistent with these assumptions,
we will consider the three boldface letters
\ctrue, \cfalse, and \tbool
all to denote \emph{types} of values.
From this perspective,
the types \ctrue and \cfalse
are each \emph{subtypes} of type \tbool.
That is,
any value of type \ctrue is also of type \tbool,
but the converse does not hold.
If we view truth values as types, however,
we do so only in the ``weakly-typed'' sense
that we assume that concrete values of these types
are reliably recognizable via some computation.
We neither require all terms to have well-formed types in some type system,
nor expect that
all terms terms denoting truth values
to be syntactically distinguishable
(\eg as ``propositions'')
from terms
denoting other non-truth values.
A term is just a term,
which \emph{might} but need not denote a truth value.

\subsection{Type judgments and inference rules}
\label{sec:prop:judgments}

\newcommand{\ruleboolIa}{
	\infrule[boolI1]{a \jtrue}{a \jbool}
}
\newcommand{\ruleboolIb}{
	\infrule[boolI2]{a \jfalse}{a \jbool}
}
\newcommand{\ruleboolE}{
	\infrule[boolE]{
		a \jbool
		\qquad
		\begin{matrix}
			a \jtrue 	\\
			\vdots	\\
			c \jtrue
		\end{matrix}
		\qquad
		\begin{matrix}
			a \jfalse 	\\
			\vdots	\\
			c \jtrue
		\end{matrix}
	}{
		c \jtrue
	}
}

If $a$ is an arbitrary term,
then `$a \jtrue$' expresses a \emph{type judgment}
or claim that term $a$ denotes a true value
(\emph{any} true value if there is more than one).
Similarly, `$a \jfalse$' is a type judgment
that $a$ denotes a false value.
Finally, as mentioned above,
`$a \jbool$' is a type judgment
that $a$ denotes any boolean:
that is,
either a true value or a false value.\footnote{
	The notation used here is indebted to
	Martin-L\"of's intuitionistic type theory
	as described in \cite{martin-lof80intuitionistic}.
	\gd's logic is not intuitionistic as is Martin-L\"of's type theory,
	however.
\later{
	and \gd's type system is dyanmic (computation-driven)
	rather than static (syntax-driven)
	as is Martin-L\"of's.
}
}

We will next use judgments to form
\emph{inference rules} in traditional natural deduction style.
To illustrate,
we first introduce the following two basic inference rules:

\[
	\ruleboolIa
\qquad
	\ruleboolIb
\]

Inference rules indicate any \emph{premises} above the line,
a \emph{conclusion} below the line,
and optionally a \emph{label} for the inference rule to the right.
The first rule above, \irl{boolI1},
states that if it is known (\ie already proven)
that term $a$ denotes a true value,
then we may safely infer the weaker conclusion
that $a$ denotes some boolean truth value
(\ie that $a$ is either true or false).
The second rule similarly allows us
to infer the weaker type judgment `$a \jbool$'
if we have already proven the type judgment `$a \jfalse$'.

The next inference rule illustrates multiple premises
and hypothetical inference:

$$
	\ruleboolE
$$

This rule states that we can draw the conclusion `$c \jtrue$'
provided we first satisfy three conditions stated by the premises.
The first premise `$a \jbool$' states that
term $a$ must first be known (\ie already proven) to be boolean.
Second, starting from a hypothetical assumption of `$a \jtrue$',
we must be able to derive through some correct chain of reasoning
the conclusion `$c \jtrue$'.
Finally, starting from the contrary hypothesis `$a \jfalse$',
we must likewise be able to derive
the same conclusion `$c \jtrue$'.
This rule in effect performs \emph{boolean case analysis}:
taking a term $a$ already known to be boolean,
allowing us to ``split'' our line of reasoning
to address the true case `$a \jtrue$' separately and differently
from the false case `$a \jfalse$',
then finally ``re-join'' our reasoning in the conclusion `$c \jtrue$'.

We may consider the first two inference rules above,
\irl{boolI1} and \irl{boolI2},
to be \emph{introduction} rules for boolean type judgments.
These rules \emph{introduce}
a type judgment of the form `$a \jbool$' into the conclusion,
provided we are reasoning forwards
from premises towards conclusion.
The last rule above, in contrast,
is an \emph{elimination} rule
for boolean type judgments.
That is, the \irl{boolE} rule
effectively \emph{eliminates}
a type judgment of the form `$a \jbool$'
from the premises,
thereby making it possible to reason
in terms of \ctrue and \cfalse judgments alone,
within the other two hypothetical premises.

We finally introduce two inference rules
that in effect define the essence of ``truth'' and ``falsity''
for our purposes in \gd,
independent of their booleanness:

\[
	\infrule[trueIE]{
		c
	}{
		c \jtrue
	}
\qquad
	\infrule[falseE]{
		a \jfalse
	\qquad
		a \jtrue
	}{
		c \jtrue
	}
\]

The first of these rules, \irl{trueIE},
expresses that the assertion ``$c$ is true'' via `$c \jtrue$'
is equivalent simply to asserting $c$ alone with no type tag at all:
that is, the true statements are exactly those we care about
for proof purposes.
At the moment this rule may seem fairly vacuous and useless,
but we will see its pragmatic utility as we progress.
For now, we will generally continue to attach the type tag `$\ttrue$'
to assertions of truth,
even though the \irl{trueIE} rule means
that attaching `$\ttrue$' to true statements is technically redundant.

The double line indicates that the rule is bidirectional,
representing both an introduction and a corresponding elimination rule at once.
Reading the bidirectional rule ``as usual''
with premise above and conclusion below the double line,
these rules serve as introduction rules.
Flipping the bidirectional rule vertically, however --
taking the judgment \emph{below} the line as the premise
and the judgment \emph{above} the line as conclusion --
we get the corresponding elimination rule.
A bidirectional rule thus states in effect that
the form of judgment above the line is logically equivalent to,
and hence freely interchangeable with,
the form of judgment below the line.

The second rule above, \irl{falseE},
expresses the principle of \emph{non-contradiction}:
it is impossible for any term $a$ to be simultaneously both true and false.
This rule expresses non-contradiction
via the classic mechanism of \emph{explosion}:
if some line of reasoning leads us to a point
where we find that a term $a$ must be both true and false,
then it is safe to conclude anything we might want to --
expressed by the arbitrary conclusion `$c \jtrue$' --
since we have anyway arrived at a logical impossibility.\footnote{
	The fact that we accept the explosion rule from classical logic
	distinguishes \gd from proposed \emph{paraconsistent} logics,
	which attempt to achieve some robustness to inconsistency.
	\later{cite: Graham Priest etc.}
	Since the explosion rule means that any single inconsistency
	immediately makes \emph{every} expressible statement provably ``true'',
	consistency is thus every bit as critical to \gd
	as it is to classical logic.
}

The inference rules above will be the only ones we need
in order to define the relationship
between boolean, true, and false type judgments
in \gd.

\subsubsection{Proof by contradiction in \gd}

As mentioned above,
provided that some term $a$ is known to denote a boolean truth value --
that is, we have already proven `$a \jbool$' --
the elimination rule \irl{boolE} allows us
to perform case analysis on $a$.
That is,
we can prove the goal `$b \jtrue$' in one fashion
in the case where $a$ happens to be true,
while we might prove the same goal in a different fashion
in the case where $a$ is false.

With this power of case analysis, for example,
we can immediately derive one rule for (grounded) proof by contradiction.
By taking $c$ to be the same as $a$ in the \irl{boolE} rule above,
the second premise of \irl{boolE} becomes trivial:
hypothesis `$a \jtrue$' leads directly
to the premise's required conclusion, which is also `$a \jtrue$'.
We thus get the following derived rule
representing a particular special case of boolean case analysis:

\[
	\infrule{
		a \jbool
		\qquad
		\begin{matrix}
		a \jfalse 	\\
		\vdots	\\
		a \jtrue
		\end{matrix}
	}{
		a \jtrue
	}
\]

That is, if $a$ is already proven to be boolean,
and if from the hypothetical assumption that $a$ is false
we can prove the contrary judgment that $a$ is true,
then we can deduce the \emph{non-hypothetical} conclusion
that $a$ must be unconditionally true.

We will shortly see how the key difference between this rule
and classical proof by contradiction --
namely, the additional `$a \jbool$' typing premise on the left side
that we find missing from classical logic --
will help \gd avoid paradoxes such as the Liar.

\subsubsection{Judgments as terms}
\label{sec:prop:judgments-as-terms}

We will normally use type judgments
like `$a \jtrue$' or `$a \jbool$'
in defining inference rules such as those above.
In \gd these judgments may also serve as ordinary terms,
however,
expressing the proposition that a term
denotes a value of a particular type.
Formalizing this principle,
the following rules express this equivalence:
first for the specific example of the \tbool type,
then in general for any name $\tau$ denoting a type:

$$
\infrule{
	a \jbool
}{
	(a \jbool) \jtrue
}
\qquad
\infrule{
	(a \jbool) \jtrue
}{
	a \jbool
}
\qquad
\infrule{
	a \judgment{$\tau$}
}{
	(a \judgment{$\tau$}) \jtrue
}
\qquad
\infrule{
	(a \judgment{$\tau$}) \jtrue
}{
	a \judgment{$\tau$}
}
$$

\later{	these could be bidirectional rules,
	if we introduce bidirectional rules earlier}

\subsection{Logical negation $\neg$}
\label{sec:prop:neg}

\newcommand{\rulenegIa}{
	\infrule[\neg I1]{a \jtrue}{\neg a \jfalse}
}
\newcommand{\rulenegIb}{
	\infrule[\neg I2]{a \jfalse}{\neg a \jtrue}
}
\newcommand{\rulenegEa}{
	\infrule[\neg E1]{\neg a \jfalse}{ a \jtrue}
}
\newcommand{\rulenegEb}{
	\infrule[\neg E2]{\neg a \jtrue}{ a \jfalse}
}

We next introduce the logical negation operator, `$\neg$'.
Given any term $a$,
we can construct a term `$\neg a$' denoting the logical negation of $a$.
The following inference rules define logical negation
in terms of the \ctrue and \cfalse type judgments above:

$$
	\rulenegIa
\qquad
	\rulenegIb
\qquad
	\rulenegEa
\qquad
	\rulenegEb
$$

These rules take the form of introduction and elimination rules, respectively,
for logical negation.
The fact that we can express both true (`$a \jtrue$')
and false type judgments (`$a \jfalse$'),
and not just the former,
allows for a simpler formulation
than the traditional introduction and elimination rules
for logical negation in classical logic.

In the interest of more concise notation,
we can combine the four inference rules above
into the following two bidirectional equivalence rules:

\[
	\infeqv[\neg IE1]{a \jtrue}{\neg a \jfalse}
\qquad
	\infeqv[\neg IE2]{a \jfalse}{\neg a \jtrue}
\]

Using the above rules and boolean case analysis (\irl{boolE}),
we can derive a bidirectional \emph{typing rule} stating that
if term $a$ is a boolean then so is `$\neg a$',
and vice versa:

$$
\com{
	\infrule[\neg TI]{a \jbool}{\neg a \jbool}
	\qquad
	\infrule[\neg TE]{\neg a \jbool}{a \jbool}
}%com
	\infeqv[\neg TIE]{a \jbool}{\neg a \jbool}
$$

We can now derive rules for proof by contradiction
and refutation by contradiction, respectively,
in terms of logical negation,
by using case analysis (\irl{boolE}) and $\neg I$:

$$
	\infrule{
		a \jbool
		\qquad
		\begin{matrix}
			\neg a \jtrue 	\\
			\vdots	\\
			a \jtrue
		\end{matrix}
	}{
		a \jtrue
	}
	\qquad
	\infrule{
		a \jbool
		\qquad
		\begin{matrix}
			a \jtrue 	\\
			\vdots	\\
			\neg a \jtrue
		\end{matrix}
	}{
		\neg a \jtrue
	}
$$

From these rules and boolean case analysis we can in turn derive
the more traditional inference rules
for negation introduction and elimination,
respectively:

$$
	\infrule{
		a \jbool
		\qquad
		\begin{matrix}
			a \jtrue 	\\
			\vdots	\\
			c \jtrue
		\end{matrix}
	}{
		\neg a \jtrue
	}
	\qquad
	\infrule{
		\neg a \jtrue \qquad a \jtrue
	}{
		c \jtrue
	}
$$

The first rule requires a chain of reasoning
leading from the hypothetical judgment `$a \jtrue$'
to a proof of an arbitrary term $c$ used only in this premise:
\ie a proof that if $a$ is true then anything is provable.
Simply taking $b$ to be `$\neg a$' converts this rule
into the earlier one for `$\neg a \jtrue$'.

The second rule similarly derives a proof of an arbitrary term $b$
from the contradictory premises of both $\neg a$ and $a$.
We derive this rule using elimination rule \irl{{\neg}E2} above
and the non-contradiction rule \irl{falseE}
from \cref{sec:truth}.

Finally we derive the most concise
of the standard rules for proof by contradiction,
namely double-negation introduction and elimination,
as a bidirectional rule:

$$
\com{
	\infrule[\neg\neg I]{a \jtrue}{\neg \neg a \jtrue}
\qquad
	\infrule[\neg\neg E]{\neg \neg a \jtrue}{a \jtrue}
}%com
	\infeqv[\neg\neg IE]{a \jtrue}{\neg \neg a \jtrue}
$$

This formulation needs no `$a \jbool$' premise
because the typing rules above
imply that $\neg \neg a$, $\neg a$, and $a$ are all boolean
provided that any one of them is boolean.
We can then derive this rule from those above by contradiction.

The fact that double-negation elimination holds in \gd
makes it immediately obvious that \gd
makes no attempt to be intuitionistic
in the tradition initiated by L.E.J. Brouwer,
which traditionally rejects this equivalence.\footnote{
	The roots of intuitionism appeared in
	Brouwer's 1907 PhD thesis, \cite{brouwer07over} (Dutch).
	This and other relevant works of Brouwer are available in English
	in \cite{brouwer75philosophy}
	and \cite{brouwer81cambridge}.
	Brouwer's ideas were further developed by others
	into formal systems of intuitionistic logic
	and constructive mathematics;
	see for example \cite{heyting71intuitionism}
	and \cite{bishop67foundations}.
	We will compare and contrast \gd as presented here
	with the tradition of intuitionistic and constructive mathematics
	as particular comparisons become relevant.
}
This is one way in which \gd may feel more familiar and accesible
than intuitionistic logic to those accustomed to classical logic,
despite the new typing requirements that \gd introduces.

\subsection{Definitions, self-reference, and paradox}
\label{sec:def}

We now introduce into \gd the ability to express \emph{definitions},
in the following form:

\[
	s \ldef d
\]

This form specifically represents a \emph{constant definition},
in which we assign an arbitrary but not-yet-used symbol, $s$,
as a \emph{constant symbol} to represent another arbitrary term $d$.
We henceforth refer to term $d$ as the \emph{expansion}
of the constant symbol $s$.
In essence,
the definition establishes the logical equivalence
of symbol $s$ with its expansion $d$,
in that either may subsequently be substituted for the other in a term.
We focus on constant definitions to keep things simpler for now,
but will introduce parameterized non-constant definitions
later in \cref{sec:quant:def}.

\subsubsection{Using definitions}
\label{sec:def:use}

We explicitly represent the use of definitions in \gd
via the following inference rules:

\[
	\infrule[{\ldef}I]{
		s \ldef d
	\qquad
		\ttc{p}{d} \jtrue
	}{
		\ttc{p}{s} \jtrue
	}
\qquad
	\infrule[{\ldef}E]{
		s \ldef d
	\qquad
		\ttc{p}{s} \jtrue
	}{
		\ttc{p}{d} \jtrue
	}
\]

The notation `$\ttc{p}{\cdot}$' in the above rules
represents a \emph{syntactic template}
that can express substitutions for free variables.
In particular, if $x$ denotes a variable,
the notation `$\ttc{p}{x}$' represents an otherwise-arbitrary term $p$
having exactly one free variable $x$.
If $d$ is a term,
the notation `$\ttc{p}{d}$' represents same term $p$
after replacing all instances of the free variable $x$ with term $d$.
The notation `$\ttc{p}{s}$' similarly represents the same term $p$
after replacing all instances of the same free variable $x$
with the defined symbol $s$.

Since the free variable $x$ itself does not appear in the above rules,
the \emph{template term} $p$ containing $x$
serves only as a context in these rules
indicating where an instance of the definition's expansion $d$
is to be replaced
with the defined symbol $s$ in the introduction rule \irl{{\ldef}I},
or vice-versa within the elimination rule \irl{{\ldef}E}.

\com{	old, Curry-style substitution syntax

The notation $a[b/y]$ represents the substitution
of term $b$ for free variable $y$ in term $a$,
while avoiding capture of bound variables (only a concern later).
The introduction rule $\ldef I$
allows introducing the defined operator $o$ in the conclusion
in replacement of its definition,
while the elimination rule $\ldef E$
does the reverse.

\baf{ decide on and standardize substitution notation.
	some good information on historical precedent is 
	\href{https://mathoverflow.net/questions/243084/history-of-the-notation-for-substitution}{here}.
	In particular, the $[a/x]b$ notation,
	with $a,b$ terms and $x$ a variable,
	apparently goes at least as far back as 
	Haskell Curry et al, Combinatory Logic (1958).
}
}%com

%\paragraph{Conditional bidirectional rule notation}
\label{sec:prop:cond-equiv}

The pair of inference rules above describing definitional substitution
have a form that will be common in \gd,
so we will use shorthand notation
that combines both rules
into a single more concise \emph{conditional bidirectional} rule
as follows:

\[
	\infceqv[{\ldef}IE]{
		s \ldef d
	}{
		\ttc{p}{d} \jtrue
	}{
		\ttc{p}{s} \jtrue
	}
\]

A rule of this form expresses essentially that
\emph{provided} the common premise above the single line on the left side 
has been satisfied
(in this case that a definition `$s \equiv d$' exists),
the premise above and conclusion below
the right-hand, double-lined part of the rule
may be used in either direction as a logical equivalence.
That is, provided there is a definition `$s \equiv d$',
we can replace `$\ttc{p}{d} \jtrue$' with `$\ttc{p}{s} \jtrue$'
and vice versa.

\subsubsection{First-class definitions versus metalogical abbreviations}

The use of definitions
is ubiquitous and essential in the normal practice
of working mathematics and theoretical computer science.
Ironically, however, definitions \emph{per se}
are often entirely missing from the formal logics
constructed and studied by logicians,
such as classical first-order logic.
This is because standard practice is to treat definitions
merely as \emph{metalogical} abbreviations or shorthand notations:
\ie textual substitutions, like macros in many programming languages,
that we could in principle just expand in our heads
before commencing the real work of logical reasoning.

For this ``definitions as shorthand abbreviations'' perspective to work,
however,
standard practice holds that definitions must be \emph{non-recursive}.
That is, the newly-defined symbol $s$ in a definition `$s \ldef d$' 
must \emph{not} appear in the expansion $d$.
Instead, the new symbol $s$ must be used
only \emph{after} the definition is complete.
This crucial restriction avoids numerous tricky issues
including the paradoxes we will explore shortly,
while also tremendously reducing the expressiveness and utility of definitions.

In \gd, in contrast,
we will treat definitions as ``first-class citizens'' of the logic,
rather than only as metalogical abbreviations.
That is, we will treat definitions like `$s \equiv d$' 
as actual steps in a formal logical proof,
just as definitions normally appear before and intermixed with theorems
in a working mathematical paper or textbook.

Both definitions and the bidirectional inference rules we have used above
have the same apparent effect, of establishing logical equivalences.
We draw an important semantic difference between them, however.
Like other inference rules, a bidirectional equivalence rule
is a purely metalogical construct:
a convention we use to describe and reason about \gd
in our informal metalogic of ordinary English
supplemented with traditional mathematical notation and concepts.
A definition, in contrast, is not just metalogical
but a first-class citizen within the logic of \gd.
Although the definitional equivalence symbol `$\equiv$'
is not part of \gd's term syntax,
this symbol \emph{is} part of \gd's proof syntax,
since definitions appear in \gd proofs alongside ordinary deductions.

We maintain the standard requirement that a given symbol $s$
must be defined only once:
a proof must have at most one definition
with a given symbol $s$ on the left-hand side.
Allowing a symbol to be redefined --
\eg to yield a true value by one definition
and a false value by another --
would of course yield immediate contradictions.

\gd will recklessly tempt fate, however,
by allowing definitions to be \emph{recursive} or self-referential.
Within a definition `$s \equiv d$',
the newly-defined symbol $s$ may also appear any number of times,
without restriction,
within the definition's right-hand-side expansion $d$.
We will shortly explore the effects of recursive definitions.

\subsubsection{The Liar Paradox}
\label{sec:prop:liar}

Let us see how our recklessly self-referential logic fares
against the venerable \emph{liar paradox},
readily expressible in words as follows:

\begin{center}
This statement is false.
\end{center}

If we suppose hypothetically that the above statement is true,
then we must logically conclude that it is false, and vice versa.
It is thus both true and false, a contradiction.

We can readily express the liar paradox in a definition of \gd
as follows:

$$
	L \ldef \neg L
$$

If allowed,
this definition would immediately doom classical logic,
which assumes that every syntactically well-formed proposition such as $L$
must be either true or false.
Applying classical proof by contradiction, for example,
we hypothetically assume `$\neg L$' is true,
then unwrap $L$'s definition once to yield `$\neg \neg L$',
and hence `$L$' by double-negation elimination,
thus contraditing our hypothesis `$\neg L$'.
Since the hypothesis `$\neg L$' led to the contradictory conclusion `$L$',
it follows that `$L$' must also be true non-hypothetically.
But then `$\neg L$' is also true non-hypothetically,
so we have an unconditional contradiction.
By the explosion principle, we can henceforth prove anything.

\gd's deduction rules above
do not permit us proof by contradiction about
just any syntactically well-formed term $a$,
however.
Instead, our proof by contradiction rules
first require us to prove `$a \jbool$':
\ie that $a$ is a term that actually denotes a boolean value
satisfying the \emph{habeas quid} principle.
Only then may we assume that $a$ must be either \ctrue or \cfalse
and invoke any flavor of the law of the excluded middle
or proof by contradiction.

In the case of the liar paradox statement `$L$',
we could prove `$L \jbool$'
if we could find a way to prove that `$\neg L$', `$\neg \neg L$',
or any other such variant denotes a boolean value.
But we will have difficulty doing so,
as we find no well-founded, non-circular grounds
to support such a claim.
In particular,
in attempting to prove that `$L$' is boolean,
we run into the practical conundrum of 
\emph{first} having to prove that `$L$' is boolean.
We can assign `$L$' no truth value
because it is \emph{ungrounded},
to adopt Kripke's terminology.\footnote{
	See \cite{kripke75outline}.
}

A conventional diagnosis of the Liar paradox
holds that `$L \ldef \neg L$' is a ``bad'' definition
because it is self-referential.
An alternative diagnosis that \gd suggests, in contrast,
is that the issue is not with the definition at all
but with classical logic's failure to enforce
the \emph{habeas quid} principle of \cref{sec:intro:gd}:
we must have a thing before we use it.
By this principle,
$L$ is no longer paradoxical
but just harmlessly meaningless.

We will of course revisit the paradox question, multiple times,
as we acquire more interesting and seemingly-dangerous logical toys to play with.

\baf{	note somewhere more precisely:
	equality by definition as a symbol, when used at all,	
	is typically a symbol used only in the (outermost) metalanguage,
	and never intended to be part of any target language.
	\gd changes that, making ``equals by definition'' part of the logic.
	We don't necessarily need a separate symbol for it in \gd either,
	since definitions are conceptually ``all in one place''
	(a list of definitions)
	separate from normal terms expressed with respect to them,
	so we could reuse the normal equality symbol
	in that list of definitions without technical ambiguity.
	Nevertheless, it seems useful to make a clear distinction
	between equality by definition from ``emergent'' equality.
	Also, perhapps most importantly,
	definitional equality ``works everywhere'' in a \gd term
	while normal equality requires a proof
	that its arguments are objects in order to mean something.

	Alternatively,
	we could emphasize this symbol as a ``logical equivalence'' symbol,
	\ie essentially as metalogical rather than target-logical,
	since it isn't and doesn't need to be part of the alphabet
	with which we describe \emph{terms} in \gd.
	The stronger substitutivity of definitional equivalence in \gd
	might seem to argue for this perspective as well.
	Note: Unicode has a different ``equivalent to'' symbol 
	but that seems to be less traditional for logical equivalence.
}

\subsection{Logical conjunction $\land$ and disjunction $\lor$}
\label{sec:prop:conj}
\label{sec:prop:disj}

\newcommand{\ruleandIa}{
	\infrule[\land I1]{
		a \jtrue \qquad b \jtrue
	}{
		a \land b \jtrue
	}
}
\newcommand{\ruleandIb}{
	\infrule[\land I2]{
		a \jfalse
	}{
		a \land b \jfalse
	}
}
\newcommand{\ruleandIc}{
	\infrule[\land I3]{
		b \jfalse
	}{
		a \land b \jfalse
	}
}
\newcommand{\ruleandEa}{
	\infrule[\land E1]{
		a \land b \jtrue
	}{
		a \jtrue
	}
}
\newcommand{\ruleandEb}{
	\infrule[\land E2]{
		a \land b \jtrue
	}{
		b \jtrue
	}
}
\newcommand{\ruleandEc}{
	\infrule[\land E3]{
		a \land b \jfalse
		\qquad
		\begin{matrix}
			a \jfalse \\
			\vdots \\
			c \jtrue \\
		\end{matrix}
		\qquad
		\begin{matrix}
			b \jfalse \\
			\vdots \\
			c \jtrue \\
		\end{matrix}
	}{
		c \jtrue
	}
}

We introduce conjunction terms of the form `$a \land b$'
with the classical deduction rules:

$$
	\ruleandIa
\qquad
	\ruleandEa
\qquad
	\ruleandEb
$$

The introduction rule \irl{\land I1}
allows us to introduce logical conjunction into a conclusion 
of the form `$a \land b \jtrue$',
contingent on the premises of `$a \jtrue$' and `$b \jtrue$'
each already holding individually.
The two elimination rules \irl{\land E1} and \irl{\land E2}
weaken the premise `$a \land b \jtrue$'
into a conclusion of `$a \jtrue$' or `$b \jtrue$' alone,
respectively.

The above rules allow us to reason only about the true cases
relating to judgments of the form `$a \land b \jtrue$'.
\com{
We will also need to reason about the false cases,
for judgments like `$a \land b \jfalse$'.
In classical logic, false-case rules for conjunction
are readily derivable from the true-case rules for disjunction.
As a result, the standard deduction rules for propositional logic
normally address only the true cases
and leave the corresponding false cases to be derived.
This approach does not immediately work in \gd, however,
without the unrestricted law of the excluded middle (LEM).
We therefore complete the basic rules for conjunction
with explicit rules addressing the `$a \land b \jfalse$' cases:
}%com
We will also need to reason about cases
in which a logical conjunction is false,
a purpose served by the following rules:

$$
	\ruleandIb
\qquad
	\ruleandIc
\qquad
	\ruleandEc
$$

The false-case introduction rules \irl{\land I2} and \irl{\land I3}
allow us to infer `$a \land b \jfalse$'
given a proof of either `$a \jfalse$' or `$b \jfalse$'.
The false-case elimination rule $E3$ 
essentially performs case analysis
on the premise `$a \land b \jfalse$' to be eliminated.
Provided the conclusion `$c \jtrue$' may be inferred separately
(and likely via different reasoning steps)
from either of the hypotheses `$a \jfalse$' or `$b \jfalse$',
the premise `$a \land b \jfalse$' ensures the conclusion `$c \jtrue$'
regardless of which of $a$ and/or $b$ are actually false.

%\subsubsection{Disjunction}
\label{sec:prop:disj}

\com{
Having now defined both the true-case and false-case rules for conjunction,
however,
we may treat a disjunction term of the form `$a \lor b$'
as an abbreviation for the conjunction `$\neg(\neg a \land \neg b)$',
via De Morgan's laws,
just as in classical logic.
(Alternatively we could as well treat disjunction as primitive
and derive conjunction from it.)
Conjunction and disjunction in \gd therefore
have the full symmetry and cross-substitutivity properties
we are accustomed to in classical logic,
unlike many other logics weakened to omit LEM.
\baf{(XXX an example or two)}

Regardless of which operators we consider primitive versus derived,
we simply state the true-case and false-case rules for disjunction:
}%com

The following rules similarly address
the true and false cases of logical disjunction:

$$
	\infrule[\lor I1]{
		a \jtrue
	}{
		a \lor b \jtrue
	}
\qquad
	\infrule[\lor I2]{
		b \jtrue
	}{
		a \lor b \jtrue
	}
\qquad
	\infrule[\lor E1]{
		a \lor b \jtrue
		\qquad
		\begin{matrix}
			a \jtrue \\
			\vdots \\
			c \jtrue \\
		\end{matrix}
		\qquad
		\begin{matrix}
			b \jtrue \\
			\vdots \\
			c \jtrue \\
		\end{matrix}
	}{
		c \jtrue
	}
$$
$$
	\infrule[\lor I3]{
		a \jfalse \qquad b \jfalse
	}{
		a \lor b \jfalse
	}
\qquad
	\infrule[\lor E2]{
		a \lor b \jfalse
	}{
		a \jfalse
	}
\qquad
	\infrule[\lor E3]{
		a \lor b \jfalse
	}{
		b \jfalse
	}
$$

The introduction rules \irl{\lor I1} and \irl{\lor I2}
introduce `$a \lor b \jtrue$'
given only an individual proof of either `$a \jtrue$' or `$b \jtrue$',
respectively.
The elimination rule essentially performs disjunctive case analysis.
Provided the conclusion `$c \jtrue$' may be proven separately
from either of the hypotheses `$a \jtrue$' or `$b \jtrue$',
the disjunction in the premise ensures the conclusion
regardless of which of $a$ and/or $b$ are in fact true.
Similarly, the corresponding false-case rules naturally mirror
the true-case rules for conjunction.

Just as in classical logic,
conjunction and disjunction in \gd are duals of each other:
we can obtain either operator's rules
by taking those of the other and swapping $\ctrue$ with $\cfalse$
and swapping `$\land$' with `$\lor$'.
As a result,
De Morgan's laws work in \gd just as in classical logic,
as we express in the following bidirectional equivalence rules:

\[
	\infeqv[\land IE]{
		\neg (\neg p \lor \neg q)
	}{
		p \land q
	}
\qquad
	\infeqv[\lor IE]{
		\neg (\neg p \land \neg q)
	}{
		p \lor q
	}
\]

The fact that De Morgan's laws continue to hold in \gd as with classical logic
may make \gd feel slightly more familiar and accessible to some,
despite the new typing requirements
that many other inference rules impose in \gd.

\later{	truth tables, citations of Kleene's 3-valued logic
	and strong conjunction/disjunction operators? }

\subsubsection{Typing rules for conjunction and disjunction}
\label{sec:prop:disj:type}

From the above rules
we can finally derive the following straightforward typing rules
for conjunction and disjunction:

$$
	\infrule[\land TI]{
		a \jbool \qquad b \jbool
	}{
		a \land b \jbool
	}
	\qquad
	\infrule[\lor TI]{
		a \jbool \qquad b \jbool
	}{
		a \lor b \jbool
	}
$$

Recall that logical negation in \gd
has a typing elimination rule \irl{\neg TE}
that works in the reverse direction,
allowing us to deduce `$a \jbool$' from `$\neg a \jbool$'.
Reverse-direction type deduction
is not so simple for conjunction or disjunction,
since the result may be boolean
even if only one of the inputs is boolean.\footnote{
	This typing behavior ultimately derives from \gd's
	adoption of Kleene's strong 3-valued semantics
	for conjunction and disjunction;
	see \cite{kleene38notation}.
}
Nevertheless, we can derive the following reverse typing rules,
reflecting the fact that \emph{at least one} of the inputs
to a conjunction or disjunction must be boolean
in order to for the result to be boolean:

$$
	\infrule[\land TE]{
		a \land b \jbool
		\qquad
		\begin{matrix}
			a \jbool	\\
			\vdots		\\
			c \jtrue	\\
		\end{matrix}
		\qquad
		\begin{matrix}
			b \jbool	\\
			\vdots		\\
			c \jtrue	\\
		\end{matrix}
	}{
		c \jtrue
	}
\qquad
	\infrule[\lor TE]{
		a \lor b \jbool
		\qquad
		\begin{matrix}
			a \jbool	\\
			\vdots		\\
			c \jtrue	\\
		\end{matrix}
		\qquad
		\begin{matrix}
			b \jbool	\\
			\vdots		\\
			c \jtrue	\\
		\end{matrix}
	}{
		c \jtrue
	}
$$

\com{
We can now derive the law of non-contradiction in conjunctive form,
as well as the law of excluded middle in disjunctive form,
each with the caveat that the relevant term must be boolean:

$$
	\cfrac{a \jbool}{\neg(a \land \neg a) \jtrue}
	\qquad
	\cfrac{a \jbool}{a \lor \neg a \jtrue}
$$

We can derive De Morgan's laws as inference rules,
again provided the constituent terms are known to be boolean:

$$
	\cfrac{a \jbool \qquad b \jbool \qquad
			\neg(a \land b) \jtrue}{
		\neg a \lor \neg b \jtrue}
	\qquad
	\cfrac{a \jbool \qquad b \jbool \qquad
			\neg a \lor \neg b \jtrue}{
		\neg(a \land b) \jtrue}
$$
$$
	\cfrac{a \jbool \qquad b \jbool \qquad
			\neg(a \lor b) \jtrue}{
		\neg a \land \neg b \jtrue}
	\qquad
	\cfrac{a \jbool \qquad b \jbool \qquad
			\neg a \land \neg b \jtrue}{
		\neg(a \lor b) \jtrue}
$$
}%com

Now that we have logical disjunction,
we might consider the booleanness of a term
in terms of logical disjunction and negation.
A term is boolean whenever its value is either true or false:
that is, we may treat `$a \jbool$'
as equivalent to `$a \lor \neg a \jtrue$':

\[
	\infeqv[boolIE]{
		a \lor \neg a
	}{
		a \jbool
	}
\]

\com{
$$
	a \jbool \ldef a \lor \neg a \jtrue
$$
}

\later{	Look up and cite where Field suggested this I think w.r.t KFS}

\subsubsection{Paradoxes revisited}

With conjunction and disjunction,
we can construct slightly more subtle and interesting
paradoxes (and non-paradoxes).
Consider the following statements intuitively, for example:

\begin{center}
\begin{tabular}{l}
$S_1$: Snow is white. \\
$S_2$: Either statement $S_1$ or statement $S_2$ is true. \\
$S_3$: Statements $S_2$ and $S_3$ are both true. \\
\end{tabular}
\end{center}

Supposing $t$ is any true term,
we can define these sentences in \gd as follows:

\begin{center}
\begin{tabular}{l}
$ S_1 \ldef \ctrue $ \\
$ S_2 \ldef S_1 \lor S_2 $ \\
$ S_3 \ldef S_2 \land S_3 $ \\
\end{tabular}
\end{center}

Statement $S_1$ is trivially true,
and only one operand of a disjunction need be true
for the disjunction to be true.
Therefore,
the truth of statement $S_1$ makes statement $S_2$ likewise true,
despite $S_2$'s self-reference in its second operand.

Statement $S_3$, however,
we find ourselves unable to prove either true or false in \gd.
Because $S_2$ is true,
$S_3$ effectively depends on its own value.
We will not be able to invoke proof by contradiction on $S_3$
without first proving it boolean,
and any such attempt will encounter the fact that $S_3$ must first
have already been proven boolean.

$S_3$ is an example of a statement Kripke would classify as
\emph{ungrounded} but {\em non-paradoxical}:
\gd does not give it a truth value because of its circular dependency,
but it could be ``forced'' to true (\eg by axiom)
without causing a logical inconsistency.

If $S_2$ happened to be false, of course,
then it would be trivial to prove $S_3$ false.

\subsection{Logical implication $\limp$ and biconditional $\liff$}
\label{sec:prop:imp}
\label{sec:prop:iff}

\com{
We don't need to introduce logical implication as a primitive in \gd.
Instead, we can treat implication as a shorthand,
using the familiar equivalence from classical logic:

$$	a \limp b \ldef \neg a \lor b	$$
}%com

Logical implication in \gd exhibits the same equivalence
as in classical logic,
which we express in the following bidirectional inference rule:

\[
	\infeqv[{\limp}IE]{
		\neg a \lor b
	}{
		a \limp b
	}
\]

Just as in classical logic,
$a$ implies $b$ precisely when
either $a$ is false or $b$ is true.

We can then derive introduction and elimination rules for implication,
mostly classical except the introduction rule requires
that the antecedent be proven to be boolean:\footnote{
	This is the point at which \gd diverges
	from most existing developments of paracomplete logics,
	as explored in \cite{maudlin06truth}
	and \cite{field08saving} for example.
	The prevailing view in these prior developments
	appears to be that weakening the introduction rule
	for logical implication in this fashion
	would render logical implication too weak to be useful.
	The contrary position that \gd suggests is essentially this:
	what if such a ``weakened'' notion of implication
	is actually not only \emph{good enough} to be useful in practice,
	but even quite intuitively reasonable
	when we view the added premise
	from a perspective of computation and typing?
}

$$
	\infrule[{\limp}I]{
		a \jbool \qquad
			\begin{matrix}
			a \jtrue	\\
			\vdots 		\\
			b \jtrue	\\
			\end{matrix}
	}{
		a \limp b \jtrue
	}
\qquad
	\infrule[{\limp}E]{
		a \limp b \jtrue \qquad a \jtrue
	}{
		b \jtrue
	}
$$

The \irl{{\limp}E} rule is identical to the classical \emph{modus ponens} rule.

We can similarly express the logical biconditional or ``if and only if'' in \gd
via the same bidirectional equivalence that applies in classical logic:

\[
	\infeqv[{\liff}IE]{
		(a \limp b) \land (b \limp a)
	}{
		a \liff b
	}
\]

Unlike implication,
the biconditional introduction rule we derive
includes premises demanding
that we first prove both terms in question to be boolean:

$$
	\infrule[{\liff}I]{
		a \jbool \qquad b \jbool \qquad
			\begin{matrix}
			a \jtrue	\\
			\vdots 		\\
			b \jtrue	\\
			\end{matrix}
			\qquad
			\begin{matrix}
			b \jtrue	\\
			\vdots 		\\
			a \jtrue	\\
			\end{matrix}
	}{
		a \liff b \jtrue
	}
$$

Two derived elimination rules, one for each direction,
work as in classical logic:

\[
	\infrule[{\liff}E1]{
		a \liff b \jtrue \qquad a \jtrue
	}{
		b \jtrue
	}
\qquad
	\infrule[{\liff}E2]{
		a \liff b \jtrue \qquad b \jtrue
	}{
		a \jtrue
	}
\]

As we did earlier with definitions in \cref{sec:prop:cond-equiv},
we can combine the two rules above into a single, more concise,
conditional bidirectional inference rule:

\[
	\infceqv[{\liff}E]{
		a \liff b \jtrue
	}{
		a \jtrue
	}{
		b \jtrue
	}
\]

From the above rules
we can further derive inference rules getting from a biconditional
``back'' to a logical implication in either direction:

$$
	\infrule[{\liff}E3]{
		a \liff b \jtrue
	}{
		a \limp b \jtrue
	}
\qquad
	\infrule[{\liff}E4]{
		a \liff b \jtrue
	}{
		b \limp a \jtrue
	}
$$

Finally, through boolean case analysis we can derive
the following type-elimination rules that apply to the biconditional
(but importantly, \emph{not} to logical implication in \gd).
In essence, a biconditional in \gd yields a boolean truth value
not only when, but \emph{only} when,
\emph{both} of the biconditional's subterms are boolean:

$$
	\infrule[{\liff}TE1]{
		a \liff b \jbool
	}{
		a \jbool
	}
\qquad
	\infrule[{\liff}TE2]{
		a \liff b \jbool
	}{
		b \jbool
	}
$$

\baf{	work out the case analysis explicitly somewhere? }

In general,
we now have the machinery necessary to represent, and prove,
any statement in classical propositional logic --
provided, of course, that the constituent terms
are first proven to be boolean as might be necessary.

\subsubsection{Curry's paradox}
\label{sec:prop:curry}

\baf{consider using C and P instead of c and p, consistent with Liar}

Another interesting paradox to examine is Curry's paradox,
which we may express informally as follows:

\begin{center}
If this statement is true, then pigs fly.
\end{center}

Curry's paradox is interesting in particular
because it relies only on logical implication,
and not on the law of excluded middle.
Curry's paradox therefore compromises even intuitionistic logics,
if they were to admit self-referential definitions such as this.

We can express Curry's paradox via a perfectly-legal definition in \gd,
however:

$$
	C \ldef C \limp P
$$

With the traditional natural deduction rule for implication,
without first proving anything else about $C$,
we can hypothetically assume $C$ and attempt to derive
arbitrary predicate $P$.
Since $C \ldef C \limp P$,
this derivation follows trivially
via \emph{modus ponens}.
But then we find that $C \limp P$ is true non-hypothetically,
that $C$ is likewise true by its definition,
and hence (again by \emph{modus ponens})
the truth of $P$, \ie pigs indeed fly.\footnote{
	For a witty satirical exploration of how our world might look
	if ``truth'' were in fact as overloaded
	as Curry's paradox would appear to make it,
	%and pigs indeed fly,
	see \cite{morrow93city}.
}

In \gd, however,
the introduction rule for `$a \limp b \jtrue$'
carries an obligation first to prove `$a \jbool$'.
Again, \emph{habeas quid}.
We will have trouble proving this for Curry's statement $C$, however,
since $C$'s implication depends on its own antecedent
and we find no grounded basis to assign any truth value to it.
As with the Liar paradox expressed in \gd,
we find ourselves \emph{first} having to prove `$C \jbool$'
in order to apply the \irl{{\limp}I} rule
in order to prove `$C \jbool$'
(or in general to prove \emph{anything} about $C$).
Thus, \gd appears to survive self-referential paradoxes
that even intuitionistic logics do not.

\baf{
XXX relate carefully to the contraction-free and detachment-free
approaches to Curry's paradox
% \url{https://plato.stanford.edu/entries/curry-paradox/#ContFreeResp}
}

%later

\baf{ note somewhere: we will later drop the \ctrue and \cfalse type judgments;
	they were only training wheels to establish clarity
	in the initial presentation and development of \gd.
	But we will keep using other type judgments
	like bool and nat... }

\baf{	in expansion: prove properties of the propositional calculus;
	observe how we are assuming that we have similar tools
	like ``not'' and ``and'' in the metalogic in order to do this proof.
	Thus, the only obvious way we have to reason about
	and prove interesting properties of these tools in a target logic
	is essentially to assume we already have them in the metalogic.
	If we believe that the tools we are reasoning about in the target logic
	are the same tools we're actually using in the metalogic,
	then it seems we are taking a leap of faith
	that these tools are in effect ``self-justifying.''
	This may seem unsatisfying:
	wouldn't it be nice to build \emph{all} our tools
	and prove them correct ``from the ground up''?
	But then what is ``the ground''?
	If we don't accept \emph{any} self-justifying tools,
	we don't seem to have any ground to stand on!
	There does not appear to be any conceivable alternative
	to starting from at least \emph{some} self-justifying reasoning tools
	of this form; our goal is more to choose them wisely
	and make them as clear and ``obviously correct'' as we can.
	Self-justifying tools here: neg, and, boolean.
	Later: general recursion (getting from PR to lambda calculus);
	LPO.
}

%later

\section{Predicate logic: reasoning about objects}
\label{sec:quant}

Moving beyond logical propositions,
we now wish to reason logically
about mathematical objects other than truth values:
\eg numbers, sets, functions, etc.
We will thus wish to have the usual predicate-logic quantifiers,
\emph{for all} ($\forall$) and \emph{there exists} ($\exists$).

\subsection{Domain of discourse and object judgments}

But what will be our \emph{domain of discourse} --
the varieties of mathematical objects that we quantify over?
In the same spirit of our earlier agnosticism
about which term values represent ``true'' or ``false'' values
and how many of each there are,
we likewise remain agnostic for now
about precisely what kinds of objects we may quantify over in \gd.
We intentionally leave this question to be answered later, separately,
in some specialization or application of the principles of \gd
that we cover here.\later{
\Cref{sec:nat} and~\cref{sec:set}, in particular,
will later explore instances of \gd 
specialized to the domains of
natural numbers and sets, respectively.
}
In software engineering terms,
we for now leave the domain of discourse
as an open ``configuration parameter'' to our predicate logic.

Instead of settling on any particular domain of discourse,
we merely introduce a new form of typing judgment
for use in our inference rules:

$$	a \jobj	$$

This judgment essentially states:
``The term $a$ denotes an object in the domain of discourse
to which the logical quantifiers apply.''

\subsection{Universal quantification}

% the parameter is \jobj here or \jnat for GA
\newcommand{\ruleforallIa}[1]{
	\infrule[\forall I1]{
		\begin{matrix}
			x #1			\\
			\vdots			\\
			\tto{p}{x} \jtrue	\\
		\end{matrix}
	}{
		\tforall{x}{\tto{p}{x} \jtrue}
	}
}
\newcommand{\ruleforallIb}[1]{
	\infrule[\forall I2]{
		a #1
		\qquad
		\tto{p}{a} \jfalse
	}{
		\tforall{x}{\tto{p}{x} \jfalse}
	}
}
\newcommand{\ruleforallEa}[1]{
	\infrule[\forall E1]{
		\tforall{x}{\tto{p}{x} \jtrue}
		\qquad
		a #1
	}{
		\tto{p}{a} \jtrue
	}
}
\newcommand{\ruleforallEb}[1]{
	\infrule[\forall E2]{
		\tforall{x}{\tto{p}{x} \jfalse}
		\qquad
		\begin{matrix}
			\underbrace{x #1 \qquad \tto{p}{x} \jfalse} \\
			\vdots					\\
			q \jtrue				\\
		\end{matrix}
	}{
		q \jtrue
	}
}
\newcommand{\ruleforallTI}[1]{
	\infrule[\forall TI]{
		\begin{matrix}
			x #1		\\
			\vdots		\\
			\tto{p}{x} \jbool	\\
		\end{matrix}
	}{
		\tforall{x}{\tto{p}{x} \jbool}
	}
}

Given this new form of type judgment,
we define natural deduction rules for the universal quantifier $\forall$
as follows:

$$ 
	\ruleforallIa{\jobj}
\qquad
	\ruleforallEa{\jobj}
$$

The notation `$\tto{p}{x}$' represents a syntactic template
as discussed earlier in \cref{sec:def:use},
except in this case the ellipsis `$\dots$'
indicates that the template term $p$
may also contain other free variables in addition to $x$.
As before, the notation `$\tto{p}{a}$'
appearing in the $\forall E1$ rule
represents the template term $p$
with another term $a$ substituted for the variable $x$
while avoiding variable capture.
\baf{	explain variable capture somewhere?}

The premise of the introduction rule \irl{\forall I1}
posits a particular unspecified but quantifiable object
denoted by some variable $x$,
and demands a proof that a predicate term $\tto{p}{x}$ is true of $x$.
This proof must thus be carried out
without any knowledge of which particular quantifiable object
the variable $x$ actually represents.
Provided such a proof can be deduced
about the unknown hypothetical object $x$,
the introduction rule concludes that 
term $\tto{p}{x}$ holds true \emph{for all} quantifiable objects $x$.

The elimination rule \irl{\forall E1}
demands that some universally quantified term
$\tforall{x}{\tto{p}{x}}$ be true,
and also that some arbitrary term $a$ of interest
is already proven to denote a quantifiable object.
Under these premises,
we reason that the term $\tto{p}{a}$,
where object term $a$ has been substituted for free variable $x$,
must be true as a special case.

The new second premise `$a \jobj$'
in the elimination rule \irl{\forall E1} represents the main difference
between universal quantification in \gd versus classical first-order logic.
Classical first-order logic assumes
that terms representing quantifiable objects
are kept syntactically separate
from logical formulas,
and hence that any term that can be substituted
for a variable $x$ in a quantifier
is necessarily a quantifiable object.
Because \gd takes it as given that terms might denote anything
(truth values, quantifiable objects, non-quantifiable objects)
or nothing
(paradoxical statements, non-terminating computations),
it becomes essential to demand proof
that $a$ in fact denotes a quantifiable object
before we safely conclude that a universally-quantified truth applies to $a$.

\baf{ somewhere discuss variables, non-capture, etc. }

As before in propositional logic,
we also need to reason about the false case of universal quantification,
\ie where there is a counterexample to the quantified predicate.
The following false-case inference rules serve this purpose:

$$
	\ruleforallIb{\jobj}
\qquad
	\ruleforallEb{\jobj}
$$

The false-case introduction rule \irl{\forall I2}
demands that some arbitrary term $a$ be known to denote a quantifiable object,
and that some predicate term $\tto{p}{x}$ with a free variable $x$
be provably false when $a$ is substituted for $x$.
Since this object $a$ serves as a counterexample
demonstrating that $\tto{p}{x}$ is not true
\emph{for all} quantifiable objects $x$,
we then conclude that the universally quantified predicate
$\tforall{x}{\tto{p}{x}}$ is false.

\com{	True for GA but maybe not true for GD in general
	where b might denote a non-quantifiable object that a depends on:

In contrast with the true-case elimination rule $\forall E1$ above
that required a proof that term $b$ represents a quantifiable object,
the false-case introduction rule $\forall I2$
requires no such proof for its counterexample term $b$.
This is because in this case there is no ``harm''
if term $b$ in fact denotes no quantifiable object at all.
If $a(b)$ may be proven false in \gd
without \emph{any} information about what $b$ denotes,
\emph{if anything},
then this means that predicate $a(x)$ simply ignores $x$.
Thus, assigning $x$ to any concrete quantifiable object
will falsify $\tforall{x}{a(x)}$ just as effectively as term $b$ would.
}%com

The false-case elimination rule \irl{\forall E2}
allows us to make use of the knowledge
that a universally-quantified statement is false
and thus has a counterexample.
The \irl{\forall E2} rule takes as premises
a universally-quantified predicate known to be false,
together with a hypothetical line of reasoning from a variable $x$
denoting an arbitrary quantifiable object
about which predicate $\tto{p}{x}$ is false,
and concluding that term $q$ is true assuming these hypotheses.
Upon satisfying these premises,
the rule allows us to conclude that $q$ is true unconditionally
(non-hypothetically).
The conclusion term $q$ may not refer
to the temporary variable $x$ used in the second hypothetical premise.

Apart from the incorporation of object typing requirements,
both of these false-case rules operate similarly
to the standard natural deduction rules for existential quantifiers
in classical first-order logic.
This should not be a surprise,
in that their goal is to reason about the
existence of a counterexample
that falsifies a universally-quantified predicate.

\subsection{Existential quantification}
\label{sec:quant:exists}

\com{
Having already defined rules
covering both the true and false cases of universal quantification,
we can now define existential quantification
as a shorthand dual for universal quantification just as in classical logic:
$\texists{x}{a(x)} \equiv \neg \tforall{x}{\neg a(x)}$.
This definition enables us to derive the following inference rules
for existential quantification:
}%com

The following rules define existential quantification in \gd:

$$
	\infrule[\exists I1]{
		a \jobj
		\qquad
		\tto{p}{a} \jtrue
	}{
		\texists{x}{\tto{p}{x}} \jtrue
	}
\qquad
	\infrule[\exists E1]{
		\texists{x}{\tto{p}{x}} \jtrue
		\qquad
		\begin{matrix}
			\underbrace{x \jobj \qquad \tto{p}{x} \jtrue} \\
			\vdots					\\
			q \jtrue				\\
		\end{matrix}
	}{
		q \jtrue
	}
$$
$$ 
	\infrule[\exists I2]{
		\begin{matrix}
			x \jobj		\\
			\vdots		\\
			\tto{p}{x} \jfalse	\\
		\end{matrix}
	}{
		\texists{x}{\tto{p}{x}} \jfalse
	}
\qquad
	\infrule[\exists E2]{
		\texists{x}{\tto{p}{x}} \jfalse
		\qquad
		a \jobj
	}{
		\tto{p}{a} \jfalse
	}
$$

Just as in classical logic,
universal and existential quantification are duals of each other in \gd.
That is,
we may obtain the rules for either from those of the other
simply by swapping \ctrue with \cfalse
and simultaneously swapping `$\forall$' with `$\exists$'.
As a result,
the classical equivalences between universal and existential quantification
continue to hold in \gd,
as expressed in the following bidirectional inference rules:

\[
	\infeqv[\forall IE]{
		\neg \texists{x}{\tto{\neg p}{x}}
	}{
		\tforall{x}{\tto{p}{x}}
	}
\qquad
	\infeqv[\exists IE]{
		\neg \tforall{x}{\tto{\neg p}{x}}
	}{
		\texists{x}{\tto{p}{x}}
	}
\]

\com{
Finally, we need a way to relate
the universal and existential quantifiers to each other
for deduction purposes.
Classical logic does this implicitly
through its background assumption that every formula has a truth value,
including quantified formula.
We need to be more cautious,
but take an approach that is similar in effect
through the following two type deduction rules
about the quantifiers:

$$
	\cfrac{
		\tforall{x}{a(x)} \jbool
	}{
		\texists{x}{a(x)} \jbool
	}\ \exists T
\qquad
	\cfrac{
		\texists{x}{a(x)} \jbool
	}{
		\tforall{x}{a(x)} \jbool
	}\ \forall T
$$

These rules enable us to use boolean case analysis
and proof by contradiction
to derive the traditional equivalences between the quantifiers:

$$
	\cfrac{
		\tforall{x}{a(x)} \jtrue
	}{
		\neg \texists{x}{\neg a(x)} \jtrue
	}
\qquad
	\cfrac{
		\texists{x}{a(x)} \jtrue
	}{
		\neg \tforall{x}{\neg a(x)} \jtrue
	}
$$

These equivalences further confirm how \gd
is more similar to classical than intuitionistic logic
in its deduction rules,
only with explicit typing of terms.

}%com

\later{
\subsection{Typing rules for quantifiers}
\label{sec:quant:type}

If a property has a well-defined boolean value for all objects $x$,
then the corresponding existentially or universally quantified expression
similarly has a boolean truth value:

$$
	\ruleforallTI{\jobj}
\qquad
	\infrule[\exists TI]{
		\begin{matrix}
			x \jobj			\\
			\vdots			\\
			\tto{p}{x} \jbool	\\
		\end{matrix}
	}{
		\texists{x}{\tto{p}{x}} \jbool
	}
$$

The type introduction rules \irl{\forall TI} and \irl{\exists TI}
are interesting for the fact that
they would be considered non-constructive,
and hence objectionable,
in Brouwer's intuitionistic tradition.
In particular,
these type-introduction rules are equivalent
to the so-called \emph{limited principle of omniscience} or LPO,
which is a trivial consequence
of the law of excluded middle (LEM) in classical logic,
but is not traditionally accepted
in the intuitionistic notion of constructivism.
\gd accepts these non-constructive rules, however,
giving them a \emph{computational} if not ``intuitionistic'' interpretation,
as we will discuss further in \cref{sec:comp:intuit}.
We will sometimes want a constructive flavor of \gd,
which we obtain by omitting these typing rules.
}%later

\subsubsection{Type constraints on quantification}
\label{sec:quant:constrain}

We will often want to express quantifiers
ranging only over objects of some specific type,
such as the natural numbers to be defined later,
rather than over all quantifiable objects of any type.
We express this in \gd
by attaching type judgments to the variable bound in the quantifier,
like `$\tforall{x \jnat}{\tto{p}{x}}$' or `$\texists{x \jnat}{\tto{p}{x}}$'
to constrain $x$ to natural numbers alone
and not any other types of objects that might exist.
We consider this notation to be equivalent
to `$\tforall{x}{(x \jnat) \limp \tto{p}{x}}$'
or `$\texists{x}{(x \jnat) \limp \tto{p}{x}}$',
respectively.
Type-constrained quantification thus
relies on logical implication
and the use of type judgments as terms
as discussed earlier in \cref{sec:prop:judgments-as-terms}.

\subsection{Equality}
\label{sec:eq}

\newcommand{\ruleeqR}{
	\infrule[{=}R]{
		a \jobj
	}{
		a = a \jtrue
	}
}
\newcommand{\ruleeqS}{
	\infrule[{=}S]{
		a = b \jtrue
	}{
		b = a \jtrue
	}
}
\newcommand{\ruleeqT}{
	\infrule[{=}T]{
		a = b \jtrue	\qquad	b = c \jtrue
	}{
		a = c \jtrue
	}
}
\newcommand{\ruleeqE}{
	\infceqv[{=}E]{
		a = b \jtrue
	}{
		\tto{p}{a} \jtrue
	}{
		\tto{p}{b} \jtrue
	}
}
\newcommand{\ruleeqTI}{
	\infrule[{=}TI]{
		a \jobj
		\qquad
		b \jobj
	}{
		a = b \jbool
	}
}
\newcommand{\ruleeqTEa}{
	\infrule[{=}TE1]{
		a = b \jbool
	}{
		a \jobj
	}
}
\newcommand{\ruleeqTEb}{
	\infrule[{=}TE2]{
		a = b \jbool
	}{
		b \jobj
	}
}

In the modern tradition of incorporating the concept of equality
as an optional but common fragment of first-order logic,
we now define the notion of equality in \gd.
In particular,
equality in \gd has the standard properties
of reflexivity (${=}R$), symmetry (${=}S$), and transitivity (${=}T$),
as expressed in the following rules:

$$
	\ruleeqR
\qquad
	\ruleeqS
\qquad
	\ruleeqT
$$

The reflexivity rule requires $a$ to be a quantifiable object
as a precondition on our inferring that $a$ is equal to anything,
even equal to itself.
This typing discipline is inessential but pragmatically useful
so that the fact of two objects $a$ and $b$
being comparable at all (\ie `$a = b \jbool$')
will entail that $a$ and $b$ are both quantifiable objects.
This will help us avoid the need for too many typing premises
in other rules.
As a result, in particular,
the symmetry and transitivity rules need no type premises,
as their equality premises ensure
that the terms known to be equal must denote objects.

We also maintain the traditional property
that objects known to be equal may be substituted for each other,
which we express via the following elimination rule:
$$
	\ruleeqE
$$

That is, whenever terms $a$ and $b$ are known to be equal,
instances of $a$ may be replaced with $b$, and vice versa,
within another term $p$.
\later{	Clarify, maybe in footnote:
	applies to transparent contexts,
	not within quotes and the like. }

\subsubsection{Typing rules for equality}
\label{sec:eq:typing}

We next introduce typing rules for equality:

$$
	\ruleeqTI
\qquad
	\ruleeqTEa
\qquad
	\ruleeqTEb
$$

The first rule \irl{{=}TI} expresses the standard mathematical principle
that any two quantifiable objects may be compared,
yielding some definite truth as to whether they are equal or not.
We could alternatively adopt weaker rules,
in which perhaps only some quantifiable objects
may be tested for equality,
and perhaps only with some but not all others,
to yield well-defined results.
Such weaker alternatives would significantly complicate
reasoning about equality, however,
and would depart from the now-ubiquitous practice
of expecting essentially all mathematical objects to be comparable.

\com{
This formulation entails that \emph{all} quantifiable objects
may be tested for equality, at least with themselves.
We could weaken this formulation
such that only certain types of objects --
such as a type of ``comparable'' objects --
have any meaningful notion of equality.
On the one hand, this weaker formulation
would ``assume less'' about the nature of quantifiable objects
and hence leave more latitude for inference rules specific to object type.
On the other hand,
the notion of equality is such a fixture
of practically all mathematical objects
that it is hard to envision having a ``comparability assumption''
attached to all quantifiable objects getting in the way of anything useful.
Adding a separate ``comparable'' type judgment 
would thus seem to add unjustified complexity to \gd.
}

\com{
In the process we introduce a new typing judgment ``$a \hbox{ comp}$''
expressing the claim that term $a$ denotes an object
that is \emph{comparable} to other (comparable) objects via equality.
We leave deliberately unspecified for now, however,
questions such as which kinds of values are and aren't comparable,
whether some or all of these comparable values are also quantifiable objects,
and so on.
}

The last two elimination rules \irl{{=}TE1} and \irl{{=}TE2}
are technically redundant with each other, of course,
as either can be derived from the other
using the symmetry rule \irl{{=}S} above.
We include both merely for\ldots well, symmetry.

We will need to reason not only about equality but also about inequality --
``not equals'' -- 
which we define via the following rules:

\[
	\infeqv[{\ne}IE]{
			a = b \jfalse
		}{
			a \ne b \jtrue
		}
\qquad
	\infrule[{\ne}S]{
		a \ne b \jtrue
	}{
		b \ne a \jtrue
	}
\com{
\qquad
	\infrule[{\ne}T]{
		a = b \jtrue	\qquad	b \ne c \jtrue
	}{
		a \ne c \jtrue
	}
}%com
\]

The first bidirectional rule \irl{{\ne}IE}
simply states the standard principle that inequality
means the same as ``not equal to''.
The second rule \irl{{\ne}S} expresses that,
like equality, 
inequality is symmetric.
Unlike equality,
however,
inequality is neither reflexive nor transitive.
\com{
The last rule \irl{{\ne}T}, however,
expresses that inequality holds when equals is substituted for equals,
and is trivially derivable from rule \irl{{=}T} above.
}
%\baf{ symmetry of inequality can probably be derived
%	from symmetry of equality and the boolean typing rules? }
We can then derive typing rules for inequality:

\[
	\infrule[{\ne}TI]{
		a \jobj
		\qquad
		b \jobj
	}{
		a \ne b \jbool
	}
\qquad
	\infrule[{\ne}TE1]{
		a \ne b \jbool
	}{
		a \jobj
	}
\qquad
	\infrule[{\ne}TE2]{
		a \ne b \jbool
	}{
		b \jobj
	}
\]

\later{

Having formulated a notion of equality in \gd,
we can now define a \emph{unique existence} quantifier
in the usual way:

$$
\exists!x. p(x) \equiv \exists x. (p(x) \land \forall y. (p(y) \limp y=x))
$$

Inference rules:
$$
\infrule[\exists!I]{
	b \jobj
	\qquad
	a(b) \jtrue
	\qquad
	\begin{matrix}
		\hline
		x \jobj
		\qquad
		y \jobj
		\qquad
		a(x) \jtrue
		\qquad
		a(y) \jtrue
		\\
		\vdots
		\\
		x=y \jtrue
	\end{matrix}
}{
	\exists! x. a(x) \jtrue
}
$$
\baf{elimination, false cases? or too much trouble?}

}%later

\later{
\input{sel}
}

\subsection{Parameterized function and predicate definitions}
\label{sec:quant:def}

Now that we have notation and some rules for reasoning about objects,
it becomes more essential
to extend our earlier characterization of first-class definitions of \gd,
in \cref{sec:def},
to allow for non-constant, parameterized definitions.
Adopting a common shorthand,
we will use the notation $\vec{x}$
to represent a finite list of variables $x_1,\dots,x_n$
for some arbitrary natural number $n$.
Using this notation and the syntactic template notation used earlier,
a parameterized definition in \gd takes the following form:

\[
	s(\vec{x}) \ldef \ttc{d}{\vec{x}}
\]

A definition of this form in general defines symbol $s$
to be a function taking as \emph{formal parameters}
the list of variables $\vec{x}$.
The definition's expansion,
represented by the syntactic template $\ttc{d}{\vec{x}}$,
is simply an arbitrary term that may contain free variables
from the list $\vec{x}$.
As before in \cref{sec:def},
each symbol $s$ may be defined only once,
but the symbol $s$ may appear without restriction
within the expansion $\ttc{d}{\vec{x}}$.
This freedom gives definitions in \gd
the expressive power to represent arbitrary recursive functions.
The special case where the list of free variables $\vec{x}$ is empty
($n=0$), of course,
represents the constant-definition case described earlier.

With the basic structure of definitions generalized in this way,
we similarly generalize the inference rules
with which we use definitions for substitution within terms:

\[
\com{
	\infrule[{\ldef}I]{
		s(\vec{x}) \ldef \ttc{d}{\vec{x}}
	\qquad
		\tto{p}{\ttc{d}{\vec{a}}} \jtrue
	}{
		\tto{p}{s(\vec{a})} \jtrue
	}
\qquad
	\infrule[{\ldef}E]{
		s(\vec{x}) \ldef \ttc{d}{\vec{x}}
	\qquad
		\tto{p}{s(\vec{a})} \jtrue
	}{
		\tto{p}{\ttc{d}{\vec{a}}} \jtrue
	}
}%com
	\infceqv[{\ldef}IE]{
		s(\vec{x}) \ldef \ttc{d}{\vec{x}}
	}{
		\tto{p}{\ttc{d}{\vec{a}}} \jtrue
	}{
		\tto{p}{s(\vec{a})} \jtrue
	}
\]

Recall from \cref{sec:prop:cond-equiv}
that this conditional bidirectional rule notation demands first that
the common premise on the left side be satisfied --
\ie in this case,
that a definition of the form `$s(\vec{x}) \ldef \ttc{d}{\vec{x}}$'
has been made.
Provided this common premise is satisfied,
the rule's right-hand side may be used in either direction,
forward or reverse.
Further, the right-hand side in this rule
assumes that there is a template term $\tto{p}{y}$
containing at least one free variable $y$
(and possibly other free variables).

In the forward direction, \irl{{\ldef}IE} serves as an introduction rule,
taking as its right-hand premise the result of a double substitution.
First we take the definition's expansion $\ttc{d}{\vec{x}}$
and replace the list of formal variables $\vec{x}$
with a list of arbitrary terms represented by $\vec{a}$,
to form an \emph{instantiated expansion} term $\ttc{d}{\vec{a}}$.
We then substitute this instantiated expansion for variable $y$
in the template $\tto{p}{y}$ to form the rule's second premise.
Provided these premises are satisfied,
the introduction rule allows us
to replace all occurrences of the instantiated expansion $\ttc{d}{\vec{a}}$
with \emph{function application} terms of the form $s(\vec{a})$,
which represent calls or invocations of function symbol $s$
with actual parameters represented by the terms $\vec{a}$.
In effect, the rule introduces a function application
by ``reverse-evaluating'' the function $s$
from a result term to corresponding ``unevaluated'' function-application terms.

Operating in the reverse direction,
\irl{{\ldef}IE} functions as an elimination rule,
permitting exactly the same transformation in reverse.
That is, in the presence of a definition `$s(\vec{x}) \ldef \ttc{d}{\vec{x}}$',
the rule allows function applications of the form $s(\vec{a})$ --
where the terms $\vec{a}$ represent actual parameters to function $s$ --
to be replaced with corresponding occurrences
of $s$'s definition instantiated with these same actual parameters
to yield the instantiated expansion term $\ttc{d}{\vec{a}}$.
Thus, the rule
effectively eliminates instances of the defined symbol $s$ from term $p$
by evaluating the function in the forward direction,
\ie replacing function applications
with instantiated expansions of the function definition.

Since the formal parameters in a \gd definition
may be replaced with \emph{arbitrary terms} as actual parameters
via the above introduction and elimination rules,
and arbitrary terms in \gd may represent anything (\ie values of any type)
or nothing (\ie paradoxical or non-terminating computations),
we can similarly make no \emph{a priori} assumptions
about what these terms denote, if anything,
while performing substitutions using definitions.
We will see the importance of this principle as we further develop \gd
and make use of its power to express arbitrary recursive definitions.

In traditional mathematical practice,
a \emph{predicate} is distinct from a \emph{function}
in that a function yields values in the relevant domain of discourse
(\eg natural numbers, sets, etc.),
while a predicate yields truth values.
That is, in first-order classical logic
where terms and formulas are syntactically distinct,
a function application is a term whereas a predicate application is a formula.
In \gd, however,
since formulas are just terms that happen to (or are expected to)
yield boolean truth values,
there is similarly no special distinction
between a function definition and a predicate definition.
A predicate in \gd is merely a function that happens to (or is expected to)
result in a boolean value.

By allowing unrestricted recursive definitions into \gd,
we have in effect embedded much of the computational power
of Church's untyped lambda calculus into \gd.\footnote{
	Alonzo Church introduced the early principles
	of his untyped lambda calculus in \cite{church32set},
	but Kleene and Rosser showed this system to be inconsistent
	in \cite{kleene35inconsistency}.
	Church later presented his lambda calculus 
	in mature form in \cite{church41calculi}
	and \cite{church45introduction}.
}
If we replace the function symbols $s$
in the \irl{{\ldef}IE} rule above
with lambda terms of the form `$(\tlambda{\vec{x}}{\ttc{d}{\vec{x}}})$' --
\ie if we treat a function's ``name''
as an explicit term representation of that function's definition --
then the \irl{{\ldef}IE} rule effectively becomes
what is called $\beta$-substitution in the lambda calculus.
The untyped lambda calculus is Turing complete
and hence able to express any computable function,
so allowing unrestricted recursive definitions in \gd
clearly brings considerable computational power with it.

Despite this expressive and computational power, however,
we are not (yet) bringing into \gd \emph{higher-order} functions
as first-class objects
that we might calculate in a term or quantify over.
That is, we have defined rules for transforming
an entire function application term of the form $s(\vec{a})$
in the presence of a suitable definition of function symbol $s$,
but we have not (yet) ascribed any meaning to $s$ alone in the logic of \gd
(except in the constant definition case where $s$ has no parameters),
and we cannot quantify over functions.
We will come to the topic of higher-order functions later in \cref{sec:fun}.

\baf{	point out succinctly: proving something true
	is just executing it in reverse.}

\subsection{Conditional evaluation within predicates}
\label{sec:quant:def}

\later{switch to more concise C-style notation?}

In describing computations on objects via recursive definitions,
it will often be useful to express \emph{conditional evaluation}:
computing a value in one fashion under a certain condition,
and otherwise in a different fashion.
It has become ubiquitous in practical programming languages
to express conditional evaluation in terms of
an \kwstyle{if}/\kwstyle{then}/\kwstyle{else} construct,
whose behavior in \gd we describe
via the following rules:

\[
	\infrule[ifI1]{
		p \jtrue
	\qquad
		a \jobj
	}{
		(\tif{p}{a}{b}) = a \jtrue
	}
\qquad
	\infrule[ifI2]{
		p \jfalse
	\qquad
		b \jobj
	}{
		(\tif{p}{a}{b}) = b \jtrue
	}
\]

By the first rule,
an \kif construct is equal to its \kwstyle{then}-case term $a$
if the condition $p$ is known to be true and $a$ is an object.
Similarly,
an \kif construct is equal to its \kwstyle{else}-case term $b$
if the condition $p$ is known to be false and $b$ is an object.

A key part of the expressive power
and utility of \kwstyle{if}/\kwstyle{then}/\kwstyle{else}
is that it is \emph{polymorphic} or type-agnostic
with respect to its subterms $a$ and $b$.
That is, $a$ and $b$ can in principle denote any type of object,
not just some particular type such as boolean.

It is often particularly useful as well
that an \kif construct is non-strict:
if $p$ is true then the true-case equality holds
regardless of whether the false-case term $b$ denotes anything or nothing.
Similarly, if $p$ is false then the false-case equality holds
independent of what the true-case term $a$ might or might not denote.

The following elimination rule allows us to reason in the opposite direction,
from the knowledge that an \kif construct is equal to some object,
back to the (exactly two) possible ways
in which that equality could have been established:

\[
	\infrule[ifE]{
		(\tif{p}{a}{b}) = c \jtrue
	\qquad
		\begin{matrix}
			\underbrace{p \jtrue \qquad a = c \jtrue}	\\
			\vdots						\\
			q \jtrue					\\
		\end{matrix}
	\qquad
		\begin{matrix}
			\underbrace{p \jfalse \qquad b = c \jtrue}	\\
			\vdots						\\
			q \jtrue					\\
		\end{matrix}
	}{
		q \jtrue
	}
\]

From the above rules we can derive the following typing rules:

\[
	\infrule[ifTI]{
		p \jbool
	\qquad
		a \jobj
	\qquad
		b \jobj
	}{
		\tif{p}{a}{b} \jobj
	}
\qquad
	\infrule[ifTE]{
		\tif{p}{a}{b} \jobj
	}{
		p \jbool
	}
\]

\com{	too complex
\[
	\infceqv[ifTIE1]{
		p \jtrue
	}{
		a \jobj
	}{
		\tif{p}{a}{b} \jobj
	}
\qquad
	\infrule[ifTI2]{
		p \jfalse
	\qquad
		b \jobj
	}{
		\tif{p}{a}{b} \jobj
	}
\qquad
	\infrule[ifTI3]{
		p \jbool
	\qquad
		a \jobj
	\qquad
		b \jobj
	}{
		\tif{p}{a}{b} \jobj
	}
\]
}%com

As with the typing rules for conjunction and disjunction,
these typing rules do not express all the possible cases
due to the non-strictness of the \kif construct.
These typing rules nevertheless express useful common cases, however.

\subsection{Guarded terms}
\label{sec:quant:guard}

It will sometimes be useful in definitions
to ensure that a term denotes something meaningful
\emph{only} if some explicit condition holds.
For this purpose we can define a \emph{guard} operator `$?$'
as follows, building on the \kif construct above:

\[
	p \oq a \ldef \tif{p}{a}{\bot}
\]

In essence, 
`$p \oq a$' means the same as $a$ if $p$ is true,
and otherwise denotes nothing.

The \emph{bottom} symbol `$\bot$' in the above definition
represents any term that always denotes nothing.
We can define `$\bot$' in many ways in \gd{} --
\eg any always-paradoxical statement such as the Liar or Curry's --
or perhaps most simply as follows:

\[
	\bot \ldef \bot
\]

Using the inference rules for the \kif construct above
we can derive the following rules for guarded terms:

\[
	\infrule[?I1]{
		p \jtrue
	\qquad
		a \jobj
	}{
		p \oq a = a \jtrue
	}
\qquad
	\infrule[?E1]{
		p \oq a \jobj
	}{
		p \jtrue
	}
\qquad
	\infrule[?E2]{
		p \oq a \jobj
	}{
		p \oq a = a \jtrue
	}
\]

With this operator, for example,
we can easily define
a \emph{strict} or \emph{weak} disjunction operator `$\lor_s$'
as follows:

\[
	a \lor_s b \ldef (a \jbool) \land (b \jbool) \oq a \lor b
\]

Unlike normal disjunction,
strict disjunction yields a boolean result
only when \emph{both} of its operands are boolean,
regardless of their values.
The primary true-case introduction rules for strict disjunction
are therefore as follows:

\[
	\infrule[\lor_s I1]{
		a \jtrue
	\qquad
		b \jbool
	}{
		a \lor_s b \jtrue
	}
\qquad
	\infrule[\lor_s I2]{
		a \jbool
	\qquad
		b \jtrue
	}{
		a \lor_s b \jtrue
	}
\]

Normally it is less desirable or useful
to have more typing obligations to fill in the premises like this.
Strictness can make typing rules simpler, however,
especially reasoning in the reverse direction,
as illustrated by these typing rules for `$\lor_s$':

\[
	\infrule[\lor_s TI]{
		a \jbool
	\qquad
		b \jbool
	}{
		a \lor_s b \jbool
	}
\qquad
	\infrule[\lor_s TE1]{
		a \lor_s b \jbool
	}{
		a \jbool
	}
\qquad
	\infrule[\lor_s TE2]{
		a \lor_s b \jbool
	}{
		b \jbool
	}
\]

Recall from \cref{sec:prop:disj:type}
that we could not derive typing rules as simple as these
for normal disjunction in \gd
because `$a \lor b$' might be boolean
even if only one of $a$ or $b$ is boolean.
Thus, the guard operator allows us to ``clean up''
the semantics of a definition
and simplify inference rules derived from it
when we intend a definition for use only under certain conditions.

\section{Natural number arithmetic}
\label{sec:nat}

Now that we have logical machinery
to reason about mathematical objects via quantification and equality,
it would be nice to have some actual mathematical objects to reason about.
For purposes of ``kicking the tires''
of our new grounded deduction vehicle,
what better place to start than with the natural numbers?

As before,
we will introduce natural numbers in a form
agnostic to questions of what \emph{other} types of values,
whether quantifiable or non-quantifiable,
might exist in \gd's term space,
or what the relationship might be 
between the natural numbers and objects of other types.
For example, we will leave it unspecified for now whether or not
the natural numbers are identical to any or all truth values.
Our basic formulation here will apply equally well, for example,
to models of \gd
where 0 is false and 1 is true,
where 0 is false and any nonzero number is true,
or
where true and false are separate values
unequal to any natural number,
within or outside the domain of discourse.

\subsection{Basic deduction rules for natural numbers}
\label{sec:nat:peano}

We introduce natural numbers via deduction rules
that essentially correspond to the Peano axioms
(minus those for equality, which we obtained above):

\newcommand{\rulezeroI}{
	\infrule[0I]{
	}{
		0 \jnat
	}
}
\newcommand{\rulesucIE}{
	\infeqv[\suc IE]{
		a \jnat
	}{
		\suc(a) \jnat
	}
}
\newcommand{\rulesuceqIE}{
	\infceqv[\suc{=}IE]{
		a \jnat
	\qquad
		b \jnat
	}{
		a = b \jtrue
	}{
		\suc(a) = \suc(b) \jtrue
	}
}
\newcommand{\rulesucneIE}{
	\infceqv[\suc{\ne}IE]{
		a \jnat
	\qquad
		b \jnat
	}{
		a \ne b \jtrue
	}{
		\suc(a) \ne \suc(b) \jtrue
	}
}
\newcommand{\rulesuceqzeroI}{
	\infrule[\suc{\ne}0I]{
		a \jnat
	}{
		\suc(a) \ne 0 \jtrue
	}
}

$$
	\rulezeroI
\qquad
	\rulesucIE
$$
$$
	\rulesuceqIE
$$
$$
	\rulesucneIE
$$
$$
	\rulesuceqzeroI
$$

These rules express the basic Peano axioms
that zero is a natural number,
the successor of any natural number is a natural number,
two natural numbers are equal/unequal whenever their successors are equal/unequal,
respectively,
and the successor of any natural number is not equal to zero:
that is,
the successor function $\suc$ never ``wraps around'' to zero
as it would in modular arithmetic.

\later{		\baf{figure out if these are needed, and which}

The following typing rules ensure that 
we can infer inequality, and not just equality,
between natural numbers:

$$
	\infrule{
		a = b \jbool
	}{
		\suc(a) = \suc(b) \jbool
	}
\qquad
	\infrule{
		\suc(a) = \suc(b) \jbool
	}{
		a = b \jbool
	}
$$
}

\subsection{Natural number typing rules}
\label{sec:nat:type}

The following typing rules
relate the natural-number type \tnat
to the potentially-broader type \tobj of quantifiable objects:

$$
	\infrule[natTI]{
		a \jobj
	}{
		(a \jnat) \jbool
	}
\qquad
	\infrule[natTE]{
		a \jnat
	}{
		a \jobj
	}
$$

The natural-number type-introduction rule \irl{natTI}
states that if a term $a$ is known to denote some quantifiable object,
then there is a definite boolean ``fact of the matter''
about whether or not $a$ more specifically denotes a natural number.
If $a$ is any quantifiable object,
then either it is a natural number -- hence `$(a\jnat) \jtrue$' --
or it is some other type of object -- hence `$(a\jnat) \jfalse$'.
From a computational perspective,
this rule effectively states that
we can subject any quantifiable object $a$ to a \emph{test}
of whether it denotes a natural number,
comparable to the dynamic type checks
common in programming languages such as Python.
Such checks always work
provided the tested value $a$ indeed represents a well-defined object.
We cannot expect such a type test to work if $a$ denotes
a non-object such as a nonterminating computation, however.

The natural number type-elimination rule \irl{natTE} above
states the simpler but equally-important subtyping property
that any natural number is a quantifiable object.

\subsection{Mathematical induction}
\label{sec:nat:induct}

We next introduce a rule for mathematical induction on the natural numbers:

\[
	\infrule[Ind]{
		\tto{p}{0} \jtrue
		\qquad
		\begin{matrix}
			\underbrace{x \jnat \qquad \tto{p}{x} \jtrue} \\
			\vdots			\\
			\tto{p}{\suc(x)} \jtrue	\\
		\end{matrix}
		\qquad
		a \jnat
	}{
		\tto{p}{a} \jtrue
	}
\]

This rule expresses the standard principle
that provided a predicate term $\tto{p}{x}$ is true for the case $x=0$,
and from the premise of it being true for any given natural number $x$
we can prove that it is also true for $\suc(x)$,
then $\tto{p}{x}$ is true for any arbitrary natural number $a$.
The last constraint,
expressed by the final premise `$a \jnat$',
is important in \gd to constrain the rule's applicability
to well-defined (\ie grounded and in particular non-paradoxical)
values of the appropriate object type, \ie \tnat.

Recall that we can restrict a quantifier to objects of a particular type:
that is,
`$\tforall{x\jnat}{\tto{p}{x}}$'
is equivalent to
`$\tforall{x}{(x\jnat) \limp \tto{p}{x}}$',
as discussed earlier in \cref{sec:quant:constrain}.
Using the above rule for mathematical induction
together with the universal quantifier introduction rule \irl{\forall I}
and the type-introduction rule \irl{natTI} above,
we can derive the following, perhaps simpler and more familiar induction rule
that directly yields a quantified predicate:

$$
	\infrule{
		\tto{p}{0} \jtrue
		\qquad
		\begin{matrix}
			\underbrace{x \jnat \qquad \tto{p}{x} \jtrue} \\
			\vdots			\\
			\tto{p}{\suc(x)} \jtrue	\\
		\end{matrix}
	}{
		\tforall{x\jnat}{\tto{p}{x}} \jtrue
	}
$$

\later{
...

Natural number selection,
aka pattern matching in programming language tradition...

$$
\begin{cases}
	b	& |\ i=0 \\
	c	& |\ i=\suc(j)
\end{cases}
$$
}%later

\com{	now redundant with definitions discussion

\subsection{Defining recursive functions of natural numbers}

As a reminder,
\gd's syntax is strongly single-sorted
in the sense that the same syntactic category of terms
represent both logical formulas
(\ie terms that happen to yield boolean truth values)
and expressions of other objects such as natural numbers.
Since the unrestricted recursive definition and substitution facilities
introduced in \cref{sec:def}
apply to arbitrary terms in \gd,
we can already apply them to terms representing natural numbers.
Since definitions in \gd permit unrestricted recursion,
this means that we can use the same definition rules
to define arbitrary recursive functions of natural numbers.
That is, a definition is just a definition in \gd,
whether the defined term happens to yield a boolean or a natural number
(or something else, or nothing).

The recursive functions we define this way are not higher-functions,
since we have not (yet) defined any rules
that would make such functions into first-class objects
that we could quantify over or compare.
We could presumably do so,
but developing and exploring higher-order functions in \gd
we leave as another task for another time.
}%com

\subsection{Natural number case decomposition}
\label{sec:nat:cases}

\baf{	just make a predecessor function primitive for simplicity?}

Combining the above basic natural-number reasoning rules
with \gd's general recursive definition capability
described earlier in \cref{sec:def} and \cref{sec:quant:def},
we already have \emph{almost} the arithmetic infrastructure necessary
to express arbitrary computable functions on natural numbers.
To make our arithmetic fully useful, however,
we still seem to need one more basic mechanism,
namely natural number \emph{case decomposition}:
that is,
a means to distinguish between the ``zero''
and ``successor of something'' cases of an argument
in a recursive definition.
Defining functions of natural numbers recursively by case decomposition
is standard practice, usually just implicitly assumed to be valid,
throughout working mathematics.
Consistent with this standard practice,
we will adopt the case-decomposition notation of working mathematics
by allowing recursive definitions such as in the following example,
which defines a \emph{predecessor} function $\pred$
that simply subtracts 1 from its argument, clamping at zero:

$$
\pred(a) \equiv 
	\begin{cases}
		0		& \mid a=0 \\
		a_p		& \mid a=\suc(a_p) \\
	\end{cases}
$$

We can similarly use case decomposition to define predicate functions
that take natural numbers as inputs and yield boolean truth values,
such as this predicate determining whether its input is
an even or odd natural number:

\[
\feven(a) \equiv 
	\begin{cases}
		\ctrue			& \mid a=0 \\
		\neg \feven(a_p)	& \mid a=\suc(a_p) \\
	\end{cases}
\]

\com{	for high-level presentation only (for now)...
\[
\feven(a) \equiv 
	\begin{cases}
		\ctrue			& \mid a=0 \\
		\neg \feven(a-1)	& \mid a>0 \\
	\end{cases}
\]

Example proof that \feven is terminating/boolean:

\[
	\infrule{
		\begin{NiceMatrix}[baseline=b]
		(\ctrue) \jbool	\\
		\vdots\\
		\feven(0) \jbool \\
		\end{NiceMatrix}
	\qquad
		\begin{NiceMatrix}[baseline=b]
		\underbrace{x \jnat \qquad \feven(x) \jbool}	\\
		\vdots	\\
		\neg \feven(x) \jbool \\
		\vdots\\
		\feven(x+1) \jbool \\
		\end{NiceMatrix}
	}{
		%\tforall{x \jnat}{\feven(x) \jbool}
		\tforall{x}{\feven(x) \jbool}
	}
\]
}%com

An alternative notation,
closer to the tradition of programming language practices
in computer science,
would be to use \kcase statements or similar textual constructs.
Such notation,
particularly prominent in functional programming languages,
varies widely in details but typically looks similar to
the following syntax we will employ:

$$
\pred(a) \equiv \tcase{a}{
			\tcasei{0}{0} \mid
			\tcasei{\suc(a_p)}{a_p}
		}
$$

The following conditional bidirectional inference rules
express the basic reasoning and computational role
that case analysis provides:

\newcommand{\rulecasezeroI}{
	\infrule[\hbox{case}0I]{
		a = 0 \jtrue
		\qquad
		\tto{c}{t_0} \jtrue
	}{
		\tto{c}{\tcase{a}{
			\tcasei{0}{t_0} \mid
			\tcasei{\suc(v_p)}{\tto{t_p}{v_p}}
		}} \jtrue
	}
}
\newcommand{\rulecasezeroE}{
	\infrule[\hbox{case}0E]{
		a = 0 \jtrue
		\qquad
		\tto{c}{\tcase{a}{
			\tcasei{0}{t_0} \mid
			\tcasei{\suc(v_p)}{\tto{t_p}{v_p}}
		}} \jtrue
	}{
		\tto{c}{t_0} \jtrue
	}
}
\newcommand{\rulecasesucI}{
	\infrule[\hbox{case}\suc I]{
		a = \suc(p) \jtrue
		\qquad
		\tto{c}{\tto{t_p}{p}} \jtrue
	}{
		\tto{c}{\tcase{a}{
			\tcasei{0}{t_0} \mid
			\tcasei{\suc(v_p)}{\tto{t_p}{v_p}}
		}} \jtrue
	}
}
\newcommand{\rulecasesucE}{
	\infrule[\hbox{case}\suc E]{
		a = \suc(p) \jtrue
		\qquad
		\tto{c}{\tcase{a}{
			\tcasei{0}{t_0} \mid
			\tcasei{\suc(v_p)}{\tto{t_p}{v_p}}
		}} \jtrue
	}{
		\tto{c}{\tto{t_p}{p}} \jtrue
	}
}
\newcommand{\rulecaseeqE}{
	\infrule[]{
		\begin{aligned}
		b = \kcase\	& a\ \kof \\
				& \tcasei{0}{t_0} \mid \\
				& \tcasei{\suc(v_p)}{\tto{t_p}{v_p}}
		\jtrue
		\end{aligned}
		%b = \tcase{a}{
		%	\tcasei{0}{t_0} \mid
		%	\tcasei{\suc(v_p)}{\tto{t_p}{v_p}}
		%} \jtrue
	\qquad
		\begin{matrix}
			\underbrace{
				a = 0 \jtrue
				\qquad
				b = t_0 \jtrue
			}		\\
			\vdots		\\
			c \jtrue	\\
		\end{matrix}
	\qquad
		\begin{matrix}
			\underbrace{
				a = \suc(v_p) \jtrue
				\qquad
				b = \tto{t_p}{v_p} \jtrue
			}		\\
			\vdots		\\
			c \jtrue	\\
		\end{matrix}
	}{
		c \jtrue
	}
}

\com{
$$
	\rulecasezeroI
$$
$$
	\rulecasezeroE
$$
$$
	\rulecasesucI
$$
$$
	\rulecasesucE
$$
}%com

\[
	\infceqv[case0IE]{
		a = 0 \jtrue
	}{
		\tto{c}{t_0} \jtrue
	}{
		\tto{c}{\tcase{a}{
			\tcasei{0}{t_0} \mid
			\tcasei{\suc(v_p)}{\tto{t_p}{v_p}}
		}} \jtrue
	}
\]
\[
	\infceqv[case\suc IE]{
		a = \suc(a_p) \jtrue
	}{
		\tto{c}{\tto{t_p}{a_p}} \jtrue
	}{
		\tto{c}{\tcase{a}{
			\tcasei{0}{t_0} \mid
			\tcasei{\suc(v_p)}{\tto{t_p}{v_p}}
		}} \jtrue
	}
\]

These rules operate similarly
to the substitution rule for equality in \cref{sec:eq},
but permit substitution only if the argument $a$
is known to be zero or nonzero, respectively.

The following more-subtle inference rule
finally allows us to reason in the reverse direction
about \kcase statements.
In particular,
if it is known that a \kcase statement yields
a object equal to some term $b$,
then that result must have resulted from
either the zero-case subterm or the nonzero-case subterm.

\begin{footnotesize}
$$
	\rulecaseeqE
$$
\end{footnotesize}

\later{	provide such a reverse-inference rule for if/then/else above}

\baf{	note: we could replace case analysis
	with a combination of if/then/else
	and a predecessor primitive.
	matter of taste.}

These examples illustrate that recursive definition,
combined with case decomposition in \gd,
enable us to define the predecessor function $\pred$ in \gd.
As an alternative,
if we take the predecessor function $\pred$ to be primitive,
we could use the \kif/\kwstyle{then}/\kwstyle{else}
conditional-evaluation construct defined earlier in \cref{sec:quant:def}
with the conditional predicate $a=0$
to achieve the same effect of case decomposition on natural number $a$.
Unlike a \kcase construct, however,
an \kif construct offers no direct way to get from a nonzero natural number
to its predecessor,
which is why it seems we need a predecessor primitive $\pred$ with this approach.
Case decomposition may feel more natural
to those familiar with functional programming languages
supporting abstract data types (ADTs) and pattern matching,
while treating $\pred$ and \kif as primitive
may feel more natural to those more familiar with
more common imperative languages like C or Python.\footnote{
	As we will explore later,
	we technically do not need either case decomposition,
	\kif/\kwstyle{then}/\kwstyle{else},
	or a primitive predecessor function
	to achieve full formal power to express and reason about
	functions of natural numbers.
	Using standard techniques familiar to logicians,
	we can for example transform a 2-argument function
	yielding a natural number, such as addition,
	into a 3-argument predicate like `$\kisplus(x,y,z)$'
	that tests whether $x+y=z$,
	thereby expressing addition indirectly rather than directly.
	In this slightly-obfuscated fashion
	we can implement both a successor function-predicate `$\kisS(x,z)$'
	that yields \ctrue iff $x = \suc(y)$,
	and a predecessor function-predicate `$\kisP(x,z)$'
	that yields \ctrue iff $x = \pred(y)$.
	Using recursive definitions of such function-predicates in \gd
	we can then define arbitrary function-predicates about natural numbers
	without ever needing explicit case decomposition
	or an explicit predecessor function as a primitive.
	We include case decomposition as a primitive for now, however,
	to avoid needing such unnatural obfuscation.
}

\subsection{Basic arithmetic development}
\label{sec:nat:arith}

Although we will not elaborate on the full details,
it appears feasible to develop arithmetic in \gd
based on these foundations in mostly the usual way.
The main difference from a standard development
of Peano Arithmetic (PA) in classical logic
is the need to prove that relevant objects
are well-defined natural numbers before using them.
These proof obligations appear slightly tedious, to be sure,
but otherwise not particularly onerous or challenging,
since expressing type constraints and expectations
is a standard if usually informal and often implicit part
of standard mathematical practice.

We start by defining the small numerals in the obvious way:

$$
	1 \equiv \suc(0)
\qquad
	2 \equiv \suc(\suc(0))
\qquad
	\cdots
$$

We then define addition in the standard primitive-recursive fashion:

$$
a+b \equiv 
	\begin{cases}
		b		& |\ a=0 \\
		\suc(a_p+b)	& |\ a=\suc(a_p) \\
	\end{cases}
$$

The key difference here between \gd and standard practice, of course,
is that \gd imposed on us no \emph{a priori} constraints on recursion --
such as that the definition be structurally primitive-recursive,
or well-founded in any other sense --
before admitting the recursive definition as legitimate in \gd.
We can just as easily define nonsensical ungrounded functions
like `$F(x) \equiv F(x)$'.
\gd accepts such definitions without complaint,
but just (hopefully) will not allow us to prove much of interest
about what $F(x)$ actually denotes under such a definition.

We do, however, now have to \emph{prove}
that a sensibly-defined function such as addition
actually yields a natural number for all arguments of interest:
in this case, for all natural-number inputs $a$ and $b$.
We can do so inductively,
under the background assumptions `$a \jnat$' and `$b \jnat$',
using the typing rules introduced earlier together with
the rules for mathematical induction.

To prove that addition as defined above
is a total function on the natural numbers, for example,
we use induction on the first argument $a$
(simply because the definition above fairly arbitrarily
used decomposition on the first argument $a$ in its recursion)
to prove the proposition `$(a+b)\jnat$'.
In the base case of $a=0$,
substituting the zero case of our definition of addition
results in an obligation to prove that $b$ is a natural number,
but we already have that as a background hypothesis.
In the induction step,
we may assume that some ephemeral variable $x$
denotes a natural number (`$x\jnat$')
and that $x+b$ is already known to denote a natural number (`$(x+b)\jnat$').
We must then prove `$\suc(x+b)\jnat$'.
But since the induction hypothesis already gives us `$(x+b)\jnat$'
and the earlier type-introduction rule for successor in \cref{sec:nat:peano}
in turn allows us to infer `$\suc(x+b)\jnat$',
the induction step is likewise proven.

Having proven that adding two natural numbers
always yields a natural number,
proving the other interesting properties of addition
proceeds more-or-less as usual in Peano arithmetic or similar systems,
merely incorporating the appropriate natural number (and boolean)
typing proofs as needed throughout the deductive process.\footnote{
In a practical, automated theorem-proving or verfication system based on \gd,
we would likely hope and expect that some form of static type system --
such as the sophisticated type systems supported in proof assistants
such as Isabelle/HOL and Coq --\later{citations}
would be available to help us discharge these tedious typing deductions
throughout most proofs in most cases.
Thus, we make no pretense that \gd's built-in ``dynamic typing''
should \emph{replace} the highly-useful static type systems
ubiquitous in modern automation tools,
which will likely still be as desirable as ever.
While the static type systems of today's tools tend to be
critical to ensure the consistency of their logic, however,
this need not be the case for automation built on \gd:
a static type system might instead just be a helpful automation layer
atop \gd's fundamental ``dynamic typing'' deductions,
such that (for example) a soundness error in the static type system
simply causes the underlying dynamically-typed \gd proof to fail,
rather than introducing a logical inconsistency
that might allow nonsense to be proven.
Further, such automation based on \gd could always permit reasoning
to ``escape'' the unavoidable restrictions of the static type system --
through dynamic type tests, for example --
again without introducing any (new) risks of inconsistency
atop the underlying dynamically-typed logic.
}

We continue in this spirit merely by outlining
suitable recursive definitions in \gd
for a few more of the basic arithmetic functions,
whose developments appear feasible in essentially the same way
as in standard (\eg primitive-recursive) developments of the same functions --
only with the added obligations of inductively proving these definitions
actually yield natural numbers for all appropriate arguments,
since we can no longer assume this at the outset
due to primitive recursion or other well-foundedness constraints
in the function-definition process.

We define multiplication recursively as follows:

$$
a \times b \equiv 
	\begin{cases}
		0			& |\ a=0 \\
		(a_p \times b)+b	& |\ a=\suc(a_p) \\
	\end{cases}
$$

This approach extends naturally to exponentiation:

$$
a^b \equiv 
	\begin{cases}
		1			& |\ b=0 \\
		a \times a^{b_p}	& |\ b=\suc(b_p) \\
	\end{cases}
$$

Neither \gd's basic recursive definition facilities in \cref{sec:def},
nor the case decomposition mechanism introduced above,
nor the rules for proof by mathematical induction,
inherently ``care'' whether a defined function returns
a natural number, or a boolean, or some other type.
As a result, exactly the same facilities allow us
to define inequalities and the ordering of natural numbers
in similarly recursive style:

$$
a \le b \equiv
	\begin{cases}
		\ctrue				& |\ a=0 \\
		(a_p \ne b) \land (a_p \le b)	& |\ a=\suc(a_p) \\
	\end{cases}
$$

This definition in essence tests $a \le b$ 
by checking, recursively,
that no natural number strictly less than $a$ is equal to $b$.
This recursive definition style is not our only option:
we could alternatively use quantifiers to similar effect as in 
`$a \le b \equiv \texists{c \jnat}{a+c=b}$'.
We stick with the recursive style here
merely for consistency and illustration purposes.

Either way, the other inequalities are easily defined:

\begin{align*}
a < b	&\equiv a \le b \land a \not= b \\
a \ge b	&\equiv \neg(a < b) \\
a > b	&\equiv \neg(a \le b)
\end{align*}

The upshot is that using these slightly-more-tedious proof practices,
we can prove any primitive-recursive function or predicate
to be a terminating total function or predicate in \gd.
Having done so,
we can then reason about these primitive-recursive functions and predicates
in the same fashion as we would
in primitive recursive arithmetic or PRA.\footnote{
	The notion of primitive-recursive functions were introduced
	in \cite{dedekind88zahlen},
	then developed into a system of formal reasoning by Skolem
	as ``the recursive mode of thought''
	in \cite{skolem23begrundung}.
	English translations of these works are available
	in \cite{dedekind63essays}
	and \cite{heijenoort02frege}, respectively.
	This system was further developed and analyzed by others:
	see for example \cite{goodstein64recursive}.
	This system later became known
	as ``primitive-recursive arithmetic'' or PRA
	after Ackermann's work made it clear that this form of recursion
	could express only certain (``primitive'') recursive functions
	and not all recursive functions over the natural numbers.
	See \cite{ackermann28hilbert},
	also translated to English in \cite{heijenoort02frege}.
}
Although PRA is based on classical logic,
\gd's inference rules effectively reduce to the classical rules
whenever the new typing requirements in the premises can be discharged --
which they always can be in the case of primitive-recursive computations.
Thus, \gd with natural numbers as defined here
appears to be at least as expressive and powerful as PRA,
in terms of both computation and reasoning power.

\subsection{Ackermann's function}
\label{sec:nat:ack}

Although we make no pretense of offering a full or rigorous
development of arithmetic in \gd here,
one obvious ``burning'' question likely to be asked
is how powerful this formulation of natural-number arithmetic actually is?
For example,
is \gd \emph{only} as powerful as PRA in reasoning,
or is it \emph{more} powerful?

A well-known limitation of PRA is that
primitive recursion can express
arbitrary exponentially-growing functions,
but cannot express superexponential functions
such as Ackermann's function.\footnote{
	Ackermann defined this function,
	and proved that it is not primitive-recursive,
	in \cite{ackermann28hilbert}.
	An English translation with a historical prologue
	is available in \cite{heijenoort02frege}.
}
We may define Ackermann's function recursively
(though not primitive recursively)
as follows:

$$
A(x,y) \equiv
	\begin{cases}
		y+1		& \mid x = 0 \\
		A(x-1,1)	& \mid y = 0 \\
		A(x-1,A(x,y-1))	& \textrm{otherwise} \\
	\end{cases}
$$

While it is not yet clear what \emph{other} limitations \gd might have,
it does not appear to have \emph{this} particular limitation.
Given that \gd makes no restrictions on recursive definitions,
the above standard definition of the Ackermann function
may be simply ``dropped into'' \gd with no immediate concern.

The slightly less trivial issue, however,
is whether \gd is powerful enough to allow us to \emph{compute}
and \emph{reason about} a function like Ackermann's.
In order to do this,
we must as a starting point be able to prove Ackermann's function
to be a total function
provided that its arguments $x$ and $y$ are natural numbers.
If we cannot do this,
then just having the definition ``in the system'' will be useless.

Fortunately, proving the Ackermann function total
appears not to be a problem in \gd.
Doing so simply requires a double (nested) inductive argument:
first an outer induction on argument $x$,
then an inner induction on argument $y$.
In the base case of the outer induction for $x=0$,
we merely need to prove `$(y+1)\jnat$',
which is trivial given the argument type assumption `$y \jnat$'.
In the outer step case of $x>0$,
we start with an induction hypothesis that 
`$\tforall{y \jnat}{A(x-1,y)\jnat}$'.
We must then use an inner induction on $y$
to prove
`$\tforall{y \jnat}{A(x,y)\jnat}$'.
In the inner base case of $y=0$,
proving `$A(x-1,1)\jnat$' is direct from the outer induction hypothesis.
For the inner induction step case,
given an induction hypothesis of `$A(x,y-1)\jnat$',
we must prove `$A(x-1,A(x,y-1))\jnat$'.
The embedded invocation of `$A(x,y-1)$' yields a natural number
directly from the inner induction hypothesis,
then applying the outer induction hypothesis 
gives us `$A(x-1,A(x,y-1))\jnat$'.

We will not rehash here all the details of Ackermann's proof
that his function is actually superexponential,
growing faster than any exponential function
representable in primitive-recursive arithmetic.
It should in principle be straightforward,
if tedious, to transplant Ackermann's proof into \gd.
Again,
the only significant new proof that obligations \gd imposes atop classical logic
are to satisfy the various natural number and boolean typing requirements,
all of which should be readily satisfiable
given the well-founded structure of everything to be proven.

In summary,
when we introduce the natural numbers,
\gd's unrestricted recursive definitions enable us
not only to express and reason about
primitive-recursive functions and predicates as in PRA,
but also to express arbitrary computable functions through recursion,
and apparently to reason about them in relatively standard ways
that at least extend beyond the reasoning power of PRA.
%How far this reasoning power extends, remains to be seen.

\subsection{General recursion}

How far does \gd's power of reasoning about natural numbers
extend beyond primitive-recursive arithmetic (PRA)?
To begin with,
since \gd allows unrestricted recursive definitions,
it appears fairly straightforward to express
arbitrary (general) recursive functions of natural numbers,
in the natural and obvious way.
Being able to express general recursive functions
does not automatically tell us how much power \gd has to reason about them,
of course,
but it is a start.

In case the power to express general recursive functions is not self-evident,
notice that we can readily express in \gd
Kleene's \emph{minimization} or \emph{mu} operator `$\mu$'.\footnote{
	See \cite{kleene52introduction}, chapter XI.
}
Simplistically,
if there exists any natural number $n$
for which a function $f(n)$ yields 0,
then `$\tmu{n}{f(n)}$' denotes the least such $n$.
Provided that we additionally know that $f(n)$ is a total function
yielding a defined result for all inputs $n$,
we can express `$\tmu{n}{f(n)}$' in \gd
via a simple recursive definition such as the following:

\[
	m_f(n) \ldef \tif{f(n) = 0}{n}{m_f(n+1)}
\]

Based on this definition,
`$\tmu{n}{\ttc{f}{n}}$' is simply `$m_f(0)$',
expressing an unbounded search upwards starting from zero
to find the least $n$ for which $f(n) = 0$.

We could equivalently formulate $\mu$ in \gd
in terms of any predicate template $\ttc{p}{n}$
that provably always yields a boolean result given any natural number $n$.
In this case, `$\tmu{n}{\ttc{p}{n}}$' denotes the least $n$
such that $\ttc{p}{n}$ is \ctrue,
provided such an $n$ exists.
The corresponding recursive \gd definition for unbounded search
is then simply:

\[
	m_p(n) \ldef \tif{\ttc{p}{n}}{n}{m_p(n+1)}
\]

To satisfy the requirement that the function $\ttc{f}{n}$
or the predicate $\ttc{p}{n}$ be (proven) total,
it is sufficient, though not necessary,
that $f$ or $p$ be primitive-recursive.

A well-known result of recursion theory is that
any general-recursive function $f(\vec{x})$
may be expressed in terms
of two fixed primitive-recursive functions $U(n)$ and $T(n,e,\vec{x})$
that are independent of $f$,
a single natural number $e$ serving as an \emph{index} or \emph{code} for $f$,
and a single use of the $\mu$ operator for unbounded search,
as follows:

\[
	f(\vec{x}) \ldef U(\tmu{n}{T(n,e,\vec{x})})
\]

The fixed primitive-recursive functions $U$ and $T$ in this construction
serve essentially the same role and function
as the rules of a Universal Turing Machine or UTM,
with the parameter $n$ serving as a step count.
The primitive-recursive function $T$ essentially simulates
the recursive computation encoded by $e$ for $n$ steps,
the minimization operator $\mu$ searches for the least $n$
for which the computation terminates,
and the primitive-recursive function $U$ extracts the natural-number output
of any such terminating execution.
We could in fact choose $T$ and $U$
to encode any Turing machine formulation
or any other step-driven computational model we might like,
such as Church's untyped lambda calculus
where steps count $\beta$-reductions on lambda terms.
\later{cite}

Even if all recursive functions are readily expressible in \gd,
in the infinitude of ways
we can formulate step-driven computational models,
what can we \emph{prove} in \gd about these recursive functions?
For any particular fixed input $\vec{x}$
for which a recursive function $f(\vec{x})$ actually terminates,
we can at least prove in \gd
that this execution indeed terminates with the expected result,
given the specific inputs $\vec{x}$ in question.
At worst, we can find a step count $n$
for which a step-driven machine computing $f(\vec{x})$ terminates,
then form a proof in \gd that executes the machine in reverse,
from its final termination at step $n$ back towards its start at step $0$.
We first prove that the machine's step $n$
terminates with the appropriate output,
then we prove that step $n-1$ leads to step $n$,
and so on back to step $0$
with a correct starting state embodying the inputs $\vec{x}$.

For inputs $\vec{x}$
for which a partial recursive function $f(\vec{x})$ does \emph{not} terminate,
we expect to be able to prove nothing interesting about $f(\vec{x})$ in \gd,
at least not directly.
We will explore later in \cref{sec:meta}
how we might use deeper metalogical reasoning
to prove computations nonterminating,
but for now we leave this as a separate matter.

\subsection{Yablo's Paradox}
\label{sec:nat:yablo}

In our continuing quest to explore
how \gd holds up against various known paradoxes,
let us now consider Yablo's ``paradox without self-reference.''\footnote{
	See \cite{yablo85truth} and \cite{yablo93paradox}.
}
In brief,
suppose we have the following infinite series of statements
labeled $Y_1$, $Y_2$, and so on, as follows::

\begin{itemize}
\item[$(Y_1)$]	All statements $Y_k$ for $k>1$ are untrue.
\item[$(Y_2)$]	All statements $Y_k$ for $k>2$ are untrue.
\item[]		$\vdots$
\item[$(Y_i)$]	All statements $Y_k$ for $k>i$ are untrue.
\item[]		$\vdots$
\end{itemize}

As we can see,
none of these statements appears to reference \emph{itself}, per se:
each statement depends only on \emph{strictly higher-numbered} statements.
Nevertheless,
if such an infinite list of statements is expressible in classical logic,
then it appears to lead to a paradox otherwise quite analogous to the Liar.

Hypothetically supposing that there is a natural number $i$
such that statement $Y_i$ is true,
all statements numbered $i+1$ and higher must be false.
But then statement $Y_{i+1}$ would be true as well,
contradicting the truth of statement $Y_i$.
Supposing to the contrary that all statements $Y_i$ are false for all $i$,
this would clearly include statement $Y_1$,
whose falsity would imply
that there must exist a statement $Y_k$ for $k>1$ that is true,
again contradicting the original assumption that all the statements are false.

Using recursive definitions with natural-number parameters,
we can readily express Yablo's paradox in \gd as follows:

\[
	Y(i) \ldef \tforall{k \jnat}{k > i \limp \neg Y(k)}
\]

If we could deduce `$Y(i) \jbool$' in \gd for any natural number $i$,
then we could certainly prove a contradiction
along the above lines of reasoning.
But how to prove in \gd that any such $Y(i)$ has a truth value?
The truth of each $Y(i)$ appears to depend on
all -- or at least some -- $Y(k)$ for $k>i$ already having truth values.
Thus, proving that any Yablo statement $Y(i)$ has a truth value in \gd
requires \emph{already} having proved
that an infinite number of higher-numbered Yablo statements
already have truth values,
yielding an infinitely-ascending set of proof obligations
before we can even get started.

It is highly debatable -- and indeed debated -- whether Yablo's paradox
actually avoids self-reference.\footnote{
	See for example
	\cite{priest97yablo},
	\cite{sorensen98yablo},
	\cite{beall01yablo}, and
	\cite{bueno03paradox}.
	% more recent:
	% "Fibonacci, Yablo and the Cassationist Approach to Paradox"
	% https://kar.kent.ac.uk/9134/1/Fibonacci%2C_Yablo_and_Cassationist.pdf
	% "Yablo’s Paradox and Beginningless Time"
	% https://intapi.sciendo.com/pdf/10.2478/disp-2009-0002
	% Cook, "The Yablo Paradox"
	% https://iep.utm.edu/yablo-pa/
}
In \gd,
we needed to define the recursive,
and hence arguably self-referential,
function $Y(i)$ above.
Nevertheless,
it is clear that structurally,
the evaluation of this recursive function
for each individual value of $i$ depends \emph{only} on
the function $Y$ evaluated on arguments strictly larger than $i$.
In this sense --
interpreting $Y$ not as a single function
but as an infinite chain of dependencies,
as Yablo clearly intended in principle --
the construction does appear to avoid direct self-reference.

\gd's apparent immunity to Yablo's paradox
suggests two interesting observations.
First, \gd's resilience to paradoxes appears to hold
when natural numbers and infinite sequences are involved,
moving beyond the propositional paradoxes
we explored in \cref{sec:prop}.
Second, \gd does not seem to ``care''
whether a paradox is a result of self-reference
or an infinite chain of dependencies.
Because \gd places on the \emph{prover}
the obligation of proving that a term has a truth value
before reasoning on a basis of it having a truth value,
both self-referential and infinite-chain paradoxes
become impervious to acquiring truth values in essentially the same way:
by ensuring that there is simply ``no place to start''
in assigning truth values to ungrounded statements such as these.
The \emph{habeas quid} principle of \cref{sec:intro:gd}
again renders Yablo's paradox harmlessly meaningless in \gd,
rather than genuinely paradoxical by leading to contradiction.

%\section{Boolean values as first-class quantifiable objects}
\section{First-class booleans and type disciplines}
\label{sec:bool}

Now that we have some objects in \gd's domain of discourse,
namely the natural numbers,
it is worth exploring in more detail how \gd evolves
if we add \emph{other} types of first-class objects
to the domain of discourse.
In particular,
we have so far talked about boolean \emph{values}
(\ie true values and false values)
without making any commitments about
whether or how these boolean values
might be \emph{objects} inhabiting our domain of discourse.
Let us now retreat slightly from our prior agnosticism
and see what changes if we explicitly make boolean values
into first-class objects that we can quantify over.

In terms of the fundamental expressiveness of \gd,
it is entirely unnecessary to make boolean values first-class.
Anything we can express with first-class booleans
we can still readily express without them;
the difference is purely a matter of taste and convenience.
Further, boolean objects are strictly simpler to define and use
than the natural numbers:
\eg we do not need induction to quantify over or otherwise reason about
``all'' of the two first-class boolean objects.
As a result,
the content of this section
may appear inconsequential from a formal perspective.
We include it nevertheless for the purpose
of clarifying and systematizing \gd's type system more explicitly.

\subsection{Equality of first-class booleans}

Given that we now wish to be able to quantify over boolean values
and have the ``intuitively correct'' thing happen when we do so,
we will now introduce the following inference rule
allowing us to compare booleans for equality:

\[
	\infrule[bool{=}I]{
		a \liff b \jtrue
	}{
		a = b \jtrue
	}
\]

Combining this new rule with prior inference rules
now enables us to infer that there can be only two boolean values,
namely \ctrue and \cfalse.
That is, any true value is (now) equal to any other true value,
and any two false values are similarly equal.

Conditioned on the hypothesis that $a$ and $b$ are known to be boolean,
we can derive a bidirectional version of the above \irl{bool{=}I} rule:

\[
	\infceqv[bool{=}IE]{
		a \jbool
	\qquad
		b \jbool
	}{
		a \liff b \jtrue
	}{
		a = b \jtrue
	}
\]

The typing premises are unnecessary in the forward direction
(as the \irl{bool{=}I} rule above implies)
but are crucial in the reverse direction,
since the mere knowledge that an arbitrary term $a$
is equal to another arbitrary term $b$
does not ensure that $a$ (and hence $b$) denotes a \emph{boolean} object
to which the biconditional operator `$\liff$' applies.

\subsection{Typing rules for first-class booleans}

Exactly as we did with the natural numbers,
we introduce typing rules stating that
boolean values may be tested for booleanness,
and that boolean values are objects:

\[
	\infrule[boolTI]{
		a \jobj
	}{
		(a \jbool) \jbool
	}
\qquad
	\infrule[boolTE]{
		a \jbool
	}{
		a \jobj
	}
\]

The latter \irl{boolTE} rule, in particular,
enables us to quantify over boolean values
via the universal and existential quantifiers
as we already can for the natural numbers.

The similarity of these rules
with those for the natural numbers in \cref{sec:nat:type}
suggest that we may want a pair of rules of this kind
for any new type we may want to introduce.
If and when we get around to defining types themselves
as explicit values (perhaps even quantifiable objects) in \gd,
we may wish to collapse all these rules into
a pair of more generic rules akin to the following,
where $\tau$ represents any type:

\[
	\infrule[boolTI]{
		a \jobj
	}{
		(a\ \tau) \jbool
	}
\qquad
	\infrule[boolTE]{
		a\ \tau
	}{
		a \jobj
	}
\]

For now, however, we limit ourselves to pointing out this trend.

\subsection{Type disiciplines: agnostic, coded, or disjoint types}

While the above rules allow us to infer that the booleans constitute
exactly two quantifiable objects,
we have still made no commitments regarding the specific \emph{relationship}
between objects of \tbool type and those of any other type such as \tnat.
At the moment, it could still be the case that boolean \ctrue
is equal to the natural number 1
(or to any other natural number for that matter),
or it could be the case that boolean \ctrue
is unequal to any natural number.
Like making booleans first-class to start with,
making a declaration on this matter of disjointness
is technically unnecessary and inconsequential from a theoretical perspective:
there is nothing we might ``need'' to express
that we can't readily express in principle
while remaining agnostic regarding this question.

\begin{figure}[t]
\begin{center}
\includegraphics[width=\textwidth]{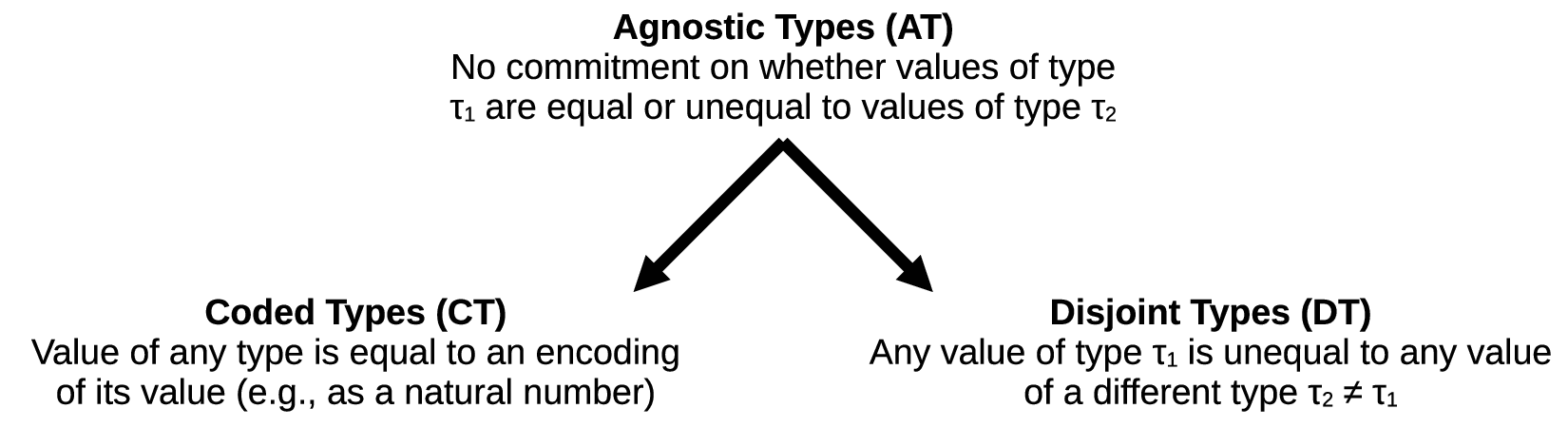}
\end{center}
\caption{A decision diagram illustrating a choice among three type disciplines:
	agnostic types, coded types, or disjoint types.}
\label{fig:type-disciplines}
\end{figure}

Nevertheless,
to map out a few interesting design alternatives,
we distinguish between three \emph{type disciplines}
illustrated in \cref{fig:type-disciplines}.
We briefly outline the alternative type disciplines
of agnostic types, coded types, and disjoint types below.

\subsubsection{Agnostic types (AT)}

In the agnostic types (AT) discipline,
we remain stubbornly uncommitted, refusing to taking a position
on all matters of whether values of two distinct types are equal or unequal.
Even though our typing rule \irl{{\ne}TI} in \cref{sec:eq:typing}
states explicitly that any two objects are indeed comparable,
in that testing them for equality yields \emph{some} boolean result,
in the AT discipline we deliberately stop short
of providing any logical basis for proving \emph{what that result is}
when the objects compared are of different types.
We simply decline to add any further inference rules or axioms
that might resolve this question.
We can still compare values of a particular type for equality or inequality,
but cannot expect to resolve specific comparisons of this kind in proofs.

Adopting this discipline is attractive
if we wish to keep our inference rules to a minimum,
given that it is generally perfectly feasible (and often standard practice)
to avoid equality comparisons across types anyway.
Further, any proofs we formulate under this discipline
remain valid under either of the two disciplines below,
so AT keeps our reasoning more generic and reusable in a sense.
We decide not to decide.

\subsubsection{Coded types (CT)}

Both in theory and in practice,
it is often useful to encode values of more complex types
into values of a single base type, such as natural numbers.
For example, we might encode boolean \cfalse as 0 and \ctrue as 1.
We will depend on such more elaborate forms of such coding
extensively below in \cref{sec:refl},
for example.
If we adopt the coded types or CT discipline,
we treat such coding as transparent,
considering a value of a non-base type
as identical and equal to its code in the base type.
Adopting this type discipline
with the common encoding of \ctrue as 1 and \cfalse as 0,
for example,
we will be able to prove the statements
`$\ctrue = 1$', `$\cfalse = 0$', `$\ctrue \ne 0$', and `$\cfalse \ne 2$'.

Adopting the CT discipline may make us guilty of violating
software-engineering best practices
and in particular type abstraction principles,
which normally hold that we should ``hide'' and not depend on
the implementation, or internal representation, of a type.
When we are reasoning along lines in which coding is central, however --
such as in G\"odel's incompleteness theorems
and related forms of reflective reasoning we will explore in later sections --
adopting the CT discipline can be convenient
and avoid unnecessary notation and explicit conversions
between a non-base type (such as a \gd term) and its code.
In the CT discipline,
a term or other value of non-base type \emph{is} its code --
\ie is equal and identical to its code --
so there is no conversion to be done.

\subsubsection{Disjoint types (DT)}

Whether for convenience, conceptual cleanliness,
strengthening \gd's type system,
or supporting best practices in software or proof engineering,
we might wish to ensure that objects of one type
are always disjoint from, and hence unequal to, objects of a different type.
As one argument for disjoint typing,
even though it may be technically workable to consider
\cfalse as ``equal to'' zero
\ctrue as ``equal to'' one as in the CT discipline,
in almost all of their endless practical uses
booleans and natural numbers play \emph{different semantic roles},
for which one is more clearly suitable than the other.
Distinguishing the different semantic roles of different types
has come to be accepted as useful in clear reasoning and communication
at the very least,
and particularly helpful when automation gets involved.
True never equals 1 because they are semantically distinct.

Taking this optional step further aligns \gd
with many familiar programming languages such as Python,
for example.
In the particular case of booleans versus natural numbers,
we might impose this disjointness requirement via the following rules:

\[
	\infrule{
		(a \jbool) \jtrue
	}{
		(a \jnat) \jfalse
	}
\qquad
	\infrule{
		(a \jnat) \jtrue
	}{
		(a \jbool) \jfalse
	}
\]

If we continue adding more types to \gd in this vein,
however,
we will find ourselves having to introduce a quadratically-exploding variety
of similar new inference rules
to obtain all the disjointness properties we desire.
At that point,
it will become pragmatically important to be able to reason about types
as first-class objects,
in order to be able to express and use a more-generic disjointness rule
of this general form:

\[
	\infrule{
		\tau_1 \jtype
	\qquad
		\tau_2 \jtype
	\qquad
		\tau_1 \ne \tau_2 \jtrue
	\qquad
		(a\ \tau_1) \jtrue
	}{
		(a\ \tau_2) \jfalse
	}
\]

Since disjointness of types is formally inessential and inconsequential
for purposes of fundamental expressiveness and reasoning, however,
we will defer further development or exploration
of first-class types and related topics for now.

\baf{
	expand and lay out three typing disciplines for \gd:
	AT agnostic typing: don't care, relationship unprovable.
	CT coded typing: all types reduce to natural numbers.
	DT disjoint typing: different types are always unequal, types have tags.
}

\section{Grounded Arithmetic (\ga)}
\label{sec:ga}

In the above development of \gd so far,
we have deliberately left the domain of discourse open-ended,
leaving the inference rules for each type open-ended
and agnostic as to whether and what kinds of \emph{other}
types of objects might exist and what their mutual relationships might be.
In a mature logic built on the principles of \gd,
we probably want a rich set of types:
other numeric types such as integers, rationals, reals, complex,
cardinals and ordinals, for example,
as well as non-numeric types such as
ordered tuples, sets, relations, functions,
the algebraic data types common in programming languages, and so on.
The earlier formulation of the natural number typing rules, for example,
is intended to allow for the graceful coexistence of the natural numbers
with other types,
whether overlapping with and properly extending the natural numbers
(\eg integers, rationals, reals)
or being entirely disjoint from them.

For foundational purposes, however,
it is also interesting to study
\gd's basic arithmetic fragment alone,
and see what might be accomplished
by building \emph{only} on the natural numbers.
For this purpose, we might ``close off'' the deduction rules defined above
and stipulate that there are \emph{no} quantifiable objects
other than the natural numbers.
In this case, \tobj becomes synonymous with \tnat,
and the typing rules defined in \cref{sec:nat:type}
relating \tobj with \tnat become superfluous.
In this environment,
two terms become provably equal ($a=b$)
only when $a$ and $b$ both represent natural numbers,
in particular the same natural number.
We may still remain agnostic as to whether the boolean truth values
are quantifiable objects (\ie identified with natural numbers)
or are entirely outside the quantification domain,
but this distinction will not matter in practice
unless we introduce axioms or inference rules sensitive to that question.

Such a formulation of arithmetic,
based on the deduction rules of \gd and the arithmetic rules above
but ``closed off'' to preclude the coexistence
of any other primitive types,
we will call \emph{grounded arithmetic} or \ga.
Besides exploring its proof and mathematical reasoning power,
we wish to ask standard metamathematical questions about it,
such as:
Is \ga consistent?
Are any consistency bugs fixable, or is it irretrievably broken?
If \ga is consistent,
in what mathematical environments can we prove it consistent?
What is the relative reasoning power of \ga
versus, say, Peano arithmetic (PA) based on classical logic,
or Heyting arithmetic (HA) based on intuitionistic logic?
Does \ga's allowance of unrestricted recursive definition
add any fundamentally interesting reasoning power beyond classical arithmetic,
or is this recursion power equivalent to what is already available
in (say) first-order PA?
What interesting challenges do \ga's typing obligations impose
on ordinary, practical mathematical proofs,
in domains like number theory or theoretical computer science for example?
We leave these intriguing questions mostly for future exploration,
but start with an attempt to define \ga more minimalistically,
and thereby in a fashion more amenable to this kind of analysis.

\subsection{Restricted syntax for grounded arithmetic (\ga)}
\label{sec:ga:syntax}

To make description and analysis of \ga more tractable,
we will present it in a more minimalistic form,
contrasting our more generalistic development of \gd so far.
We start by defining a \emph{restricted syntax} for \ga,
which falls more in line with formalizations of first-order classical logic
by virtue of syntactically distinguishing
\emph{terms} intended to denote quantifiable objects (in this case natural numbers)
from \emph{formulas} intended to denote
boolean truth values describing properties of natural numbers.

\subsubsection{Restricted term syntax}

A term $t$ has the following restricted \ga syntax:

\[
	t \equiv v \mid 0 \mid \suc(t)
\]

That is, a term $t$ intended to denote a natural number can be only
a variable $v$, the constant zero,
or the successor of a term.
A term formulated using only $0$ and $\suc$, with no variables,
we will call a \emph{literal} natural number.

\subsubsection{Restricted formula syntax}

\later{	relax to allow boolean variables?
	Would still be typable in PCF,
	and maybe needed to define substitution contexts for definitions!
	At least the way I currently define syntactic templates.
}

A formula $f$ has the following restricted \ga syntax:

\[
	f \equiv
		v \mid
		t \jnat \mid
		t_1 = t_2 \mid
		f \jbool \mid
		\neg f \mid
%		f_1 \land f_2 \mid
		f_1 \lor f_2 \mid
%		\tforall{v}{\tto{f}{v}} \mid
		\texists{v}{\tto{f}{v}} \mid
		f(\vec{t})
\]

That is,
a formula can be a variable $v$,
a \tnat or \tbool type judgment,
equality comparison,
logical negation, disjunction,
existential quantifier,
or a predicate-function application $f(\vec{t})$
taking zero or more subterm arguments denoted by $\vec{t}$.

\com{	note: as long as we define syntactic templates
	in terms of substituting terms for variables,
	we need to allow variables wherever substitutions can be done.
	We need natural number variables to allow
	substitution of equals for equals.
	We need formula variables
	in order to allow recursive definition of predicates,
	which are more necessary to \ga than recursive \tnat functios.
	But once we have both term and formula variables,
	we effectively have recursive computation
	for both natural numbers and booleans at once,
	unless we avoid using syntactic templates for recursive definitions
	and instead introduce some more restricted notion of context.
	Since the PCF language we'll be reducing to
	can easily handle and type-check both natural numbers and booleans,
	it seems better just to preserve the generality
	of allowing recursive definitions in both terms and formulas,
	even if recursive computations in terms are arguably not so useful
	without a term-level predecessor function or case decomposition.
	We could of course add term-level predecessor with no problem,
	making term computations fully general as well,
	at a slight cost in the minimality of \bga.
	Let's do this if it proves useful for some reason.
}

Even though we are now imposing
a rudimentary static type discipline by distinguishing terms from formulas,
we still retain \emph{dynamic} type judgments
of the form `$t \jnat$' and `$f \jbool$'
because we do not take for granted that a syntactically-valid term or formula
intended to represent a natural number or boolean, respectively,
will \emph{actually} do so:
they might instead describe non-terminating or paradoxical statements.

\com{	Do we need variables in formulas?
	Technically yes, for the syntactic templates
	we need for the definitional equality substitution rules.
	But those templates could also be viewed as using a relaxed syntax...}

In this restricted syntax,
function application syntax of the form `$f(\vec{t})$'
is available only within formulas (with terms as actual parameters).
Recursive definitions are therefore usable only in formulas:
that is,
\ga directly allows definitions only of recursive \emph{predicates},
not of recursive functions yielding natural numbers.
The restricted term syntax also notably omits
natural number case decomposition
or conditional evaluation via \kif constructs.
These constraints are nonessential, however,
as we will make up for this loss of term expressiveness 
via the more-expressive formula syntax.

Despite being more expressive than terms,
even the restricted formula syntax omits
the \ctrue and \cfalse judgments we used extensively in formulating \gd,
as well as the logical conjunction, implication, and biconditional connectives.
We treat these as metalogical shorthands
embodying the following equivalences:

\begin{align*}
	p \jtrue		&\ldef p \\
	p \jfalse		&\ldef \neg p \\
	p \land q		&\ldef \neg (\neg p \lor \neg q)	\\
	p \limp q		&\ldef \neg p \lor q	\\
	p \liff q		&\ldef (p \limp q) \land (q \limp p)	\\
	\tforall{x}{\tto{p}{x}}	&\ldef \neg\texists{x}{\neg\tto{p}{x}}	\\
\end{align*}

All but the last of the above equivalences
we could alternatively treat as first-class definitions within \ga.
The last would be problematic
because we have not specified a way to handle variable bindings
in first-class definitions.
\baf{	point to potential solutions, \eg using $\lambda$ as primitive? }

\subsection{Concise inference rules for grounded arithmetic (\ga)}

\begin{table}
\begin{small}
\begin{center}
\renewcommand*{\arraystretch}{0.5}	% make blank rows shorter
\begin{tabular}{|c|}
\hline
~\\
\textbf{Definition} \\
~\\
$
	\infceqv[{\ldef}IE]{
		s(\vec{x}) \ldef \ttc{d}{\vec{x}}
	}{
		\tto{p}{\ttc{d}{\vec{a}}}
	}{
		\tto{p}{s(\vec{a})}
	}
\com{	based on lambda terms
	\infeqv[appIE]{
		\tto{p}{\ttc{d}{a}}
	}{
		\tto{p}{(\tlambda{x}{\ttc{d}{x}})\ (a)}
	}
}%com
\com{	Y combinator - not needed as primitive in untyped lambda calculus
\qquad
	\infeqv[\mathbf{Y}IE]{
		\tto{p}{\ttc{d}{\ky(\tlambda{x}{\ttc{d}{x}})}}
	}{
		\tto{p}{\ky(\tlambda{x}{\ttc{d}{x}})}
	}
}
$ \\
~\\

\hline
~\\
\textbf{Equality} \\
~\\
$
	\infrule[{=}S]{
		a = b
	}{
		b = a
	}
\qquad
\com{	derivable from symmetry and substitution
	\infrule[{=}T]{
		a = b
	\qquad
		b = c
	}{
		a = c
	}
\qquad
}%com
	\infrule[{=}E]{
		a = b
	\qquad
		\tto{p}{a}
	}{
		\tto{p}{b}
	}
\qquad
%$ \\
%~\\
%$
\com{	derivable from eq-symmetry, substitution, neqIE
	\infrule[{\ne}S]{
		a \ne b
	}{
		b \ne a
	}
\qquad
}%com
	\infeqv[{\ne}IE]{
		\neg (a = b)
	}{
		a \ne b
	}
$ \\
~\\
\hline
~\\

\textbf{Natural numbers} \\
~\\
$
	\infeqv[natIE]{
		a = a
	}{
		a \jnat
	}
\qquad
	\infrule[0I]{
	}{
		0 \jnat
	}
\qquad
\com{	% both of these come as special cases of the bidirectional rule below
	\infrule[\suc I]{
		a \jnat
	}{
		\suc(a) \jnat
	}
\qquad
	\infrule[\suc{=}E]{
		\suc(a) = \suc(b)
	}{
		a = b
	}
}%com
	\infeqv[\suc{=}IE]{
		a = b
	}{
		\suc(a) = \suc(b)
	}
$ \\
~\\
$
	\infrule[\suc{\ne}0I]{
		a \jnat
	}{
		\suc(a) \ne 0
	}
\qquad
\com{	individual intro/elim rules
	\infrule[\suc{\ne}I]{
		a \ne b
	}{
		\suc(a) \ne \suc(b)
	}
\qquad
	\infrule[\suc{\ne}E]{
		\suc(a) \ne \suc(b)
	}{
		a \ne b
	}
}
	\infeqv[\suc{\ne}IE]{
		a \ne b
	}{
		\suc(a) \ne \suc(b)
	}
$ \\
~\\
%\com{	induction rule
$
	\infrule[Ind]{
		\tto{p}{0}
	\qquad
		x \jnat, \tto{p}{x} \vdash \tto{p}{\suc(x)}
	\qquad
		a \jnat
	}{
		\tto{p}{a}
	}
$ \\
%}
~\\
\hline
~\\

\textbf{Propositional logic} \\
~\\
$
	\infeqv[boolIE]{
		p \lor \neg p
	}{
		p \jbool
	}
\qquad
	\infeqv[\neg\neg IE]{
		p
	}{
		\neg\neg p
	}
\com{
	\infeqv[{\lor}IE]{
		\neg(\neg p \land \neg q)
	}{
		p \lor q
	}
\qquad
}%com
$ \\
~\\
\com{	logical and
$
	\infrule[\land I1]{
		p
		\qquad
		q
	}{
		p \land q
	}
\qquad
	\infrule[\land E1]{
		p \land q
	}{
		p
	}
\qquad
	\infrule[\land E2]{
		p \land q
	}{
		q
	}
$ \\
~\\
$
	\infrule[\land I2]{
		\neg p
	}{
		\neg (p \land q)
	}
\qquad
	\infrule[\land I3]{
		\neg q
	}{
		\neg (p \land q)
	}
\qquad
	\infrule[\land E3]{
		\neg (p \land q)
	\qquad
		\neg p \vdash r
	\qquad
		\neg q \vdash r
	}{
		r
	}
$ \\
}%com
$
	\infrule[\lor I1]{
		p
	}{
		p \lor q
	}
\qquad
	\infrule[\lor I2]{
		q
	}{
		p \lor q
	}
\qquad
	\infrule[\lor E1]{
		p \lor q
	\qquad
		p \vdash r
	\qquad
		q \vdash r
	}{
		r
	}
$ \\
~\\
$
	\infrule[\lor I3]{
		\neg p
		\qquad
		\neg q
	}{
		\neg (p \lor q)
	}
\qquad
	\infrule[\lor E2]{
		\neg (p \lor q)
	}{
		\neg p
	}
\qquad
	\infrule[\lor E3]{
		\neg (p \lor q)
	}{
		\neg q
	}
$ \\
~\\
\com{	implication and biconditional: not needed in the core
$
	\infeqv[{\limp}IE]{
		\neg p \lor q
	}{
		p \limp q
	}
\qquad
	\infeqv[{\liff}IE]{
		(p \limp q) \land (q \limp p)
	}{
		p \liff q
	}
$ \\
~\\
}%com
\hline
~\\

\textbf{Predicate logic%
	 \com{for natural numbers}
	-- omitted from \bga} \\
~\\
\com{	leave this as a metalogical equivalence
$
	\infeqv[\exists IE]{
		\neg \tforall{x}{\tto{\neg p}{x}}
	}{
		\texists{x}{\tto{p}{x}}
	}
$ \\
~\\
}%com
\com{	forall
$
\com{	generic, non-inductive forall introduction
	\infrule[\forall I1]{
		x \jnat \vdash \tto{p}{x}
	}{
		\tforall{x}{\tto{p}{x}}
	}
}%com
	\infrule[\forall I1]{
		\tto{p}{0}
	\qquad
		x \jnat, \tto{p}{x} \vdash \tto{p}{\suc(x)}
	}{
		\tforall{x}{\tto{p}{x}}
	}
\qquad
	\infrule[\forall E1]{
		\tforall{x}{\tto{p}{x}}
	\qquad
		a \jnat
	}{
		\tto{p}{a}
	}
$ \\
~\\
$
	\infrule[\forall I2]{
		a \jnat
	\qquad
		\neg \tto{p}{a}
	}{
		\neg \tforall{x}{\tto{p}{x}}
	}
\qquad
	\infrule[\forall E2]{
		\neg \tforall{x}{\tto{p}{x}}
	\qquad
		x \jnat, \neg \tto{p}{x} \vdash \ttc{c}{\ttmore}
	}{
		\ttc{c}{\ttmore}
	}
$ \\
}%com
$
	\infrule[\exists I1]{
		a \jnat
	\qquad
		\tto{p}{a}
	}{
		\texists{x}{\tto{p}{x}}
	}
\qquad
	\infrule[\exists E1]{
		\texists{x}{\tto{p}{x}}
	\qquad
		x \jnat, \tto{p}{x} \vdash \ttc{q}{\ttmore}
	}{
		\ttc{q}{\ttmore}
	}
$ \\
~\\
$
\com{	inductive version, unnecessary if we have separate induction rule
	\infrule[\exists I2]{
		\neg\tto{p}{0}
	\qquad
		x \jnat, \neg\tto{p}{x} \vdash \neg\tto{p}{\suc(x)}
	}{
		\neg\texists{x}{\tto{p}{x}}
	}
}%com
	\infrule[\exists I2]{
		x \jnat \vdash \neg\tto{p}{x}
	}{
		\neg\texists{x}{\tto{p}{x}}
	}
\qquad
	\infrule[\exists E2]{
		\neg\texists{x}{\tto{p}{x}}
	\qquad
		a \jnat
	}{
		\neg\tto{p}{a}
	}
$ \\
~\\
\hline
%~\\
%
%% "Limited principle of omniscience"?
%% "Limited principle of optimism"?
%\textbf{Predicate typing
%	-- omitted from \bga and \cga} \\
%~\\
%$
%	\infrule[\exists TI]{
%		x \jnat \vdash \tto{p}{x} \jbool
%	}{
%		\texists{x}{\tto{p}{x}} \jbool
%	}
%$ \\
%~\\
%\hline
\end{tabular}
\end{center}
\end{small}
\caption{Inference rules for Grounded Arithmetic (\ga)}
\label{tab:ga:rules}
\end{table}

\Cref{tab:ga:rules} presents
all the inference rules comprising \ga,
in a more concise and minimalistic form than we have used so far.
This formulation is more minimalistic in that it avoids including rules
that can be readily derived from combinations of other already-included rules.
It also uses Gentzen's sequent-style syntax
with the \emph{turnstile} or \emph{entailment} symbol `$\vdash$',
rather than vertical ellipsis,
to express hypothetical chains of reasoning from given assumptions.
\later{	cite Gentzen sequent calculus}
Compare, for example,
the notation for the \irl{{\lor}E1} rule in \cref{tab:ga:rules}
versus the functionally-identical rule shown earlier in \cref{sec:prop:disj}.
This difference reflects that earlier the highest priority
was clarity and obviousness,
whereas conciseness now takes a higher priority.

We treat a proof in \ga as a list of \emph{deductions}
taking the form either `$s(\vec{x}) \ldef \ttc{d}{\vec{x}}$'
or `$\Gamma \vdash p$'.
That is,
each line in a valid proof is either a definition of a new symbol
or a logical entailment derived from an inference rule.
Note that the symbols `$\ldef$' and `$\vdash$' 
are part of the proof syntax
but are not part of the restricted term or furmula syntax.

\subsection{\bga: the basic quantifier-free fragment of \ga}
\label{sec:ga:frag}

Even though \ga is already much more restricted
than the general framework for grounded deduction
that we developed in prior sections,
we will later have use for a couple still-more-restrictive fragments of \ga.

We define \emph{basic grounded arithmetic} or \bga
as the variant of \ga that we obtain
by omitting the existential quantifier `$\exists$'
and all the inference rules pertaining to quantifiers.
\bga is thus analogous to Gentzen's primitive-recursive arithmetic (PRA)
by virtue of permitting no direct expression of quantifiers,
but only the implicit top-level quantification expressed by free variables.
\baf{would primitive GA, PGA, be a better term?}

\baf{Except I think that \bga will need $\pred$ and $\kif$ as primitives
	to make up for the loss of the quantifiers... }

\later{
We similarly define \emph{constructive grounded arithmetic} or \cga
as the variant of \ga that we obtain
by including the standard introduction and elimination rules
for the existential quantifier,
both in terms of positive and negative existence,
but omitting the type-introduction rule \irl{\exists TI}.
We will explore later in \cref{XXX}
how the inclusion or omission of this rule effectively allows, or forbids,
non-constructive reasoning along lines
that intuitionists in Brouwer's tradition would typically object to.
\later{citations}
}%later

\subsubsection{Expanding the inference rules to handle background assumptions}

The rules shown in \cref{tab:ga:rules} alone,
while concisely specifying the key logical rules embodying \gd,
do not yet quite completely describe a ``working'' deduction system.
Both the Hilbert-style natural deduction rules presented earlier
and the concise rules in \cref{tab:ga:rules}
implicitly assume, but otherwise ignore, the fact that
in practice we often need to make the specified deductions
in the context of \emph{background assumptions}:
additional assumptions not immediately relevant to the present rule
but which might be crucial in later steps of a proof.

\begin{table}
\begin{small}
\begin{center}
\renewcommand*{\arraystretch}{0.5}	% make blank rows shorter
\begin{tabular}{|c|}
\hline
~\\
\textbf{Structural rules for deduction} \\
~\\
$
	\infrule[H]{
	}{
		\Gamma, p \vdash p
	}
\qquad
	\infrule[W]{
		\Gamma \vdash q
	}{
		\Gamma, p \vdash q
	}
\qquad
	\infrule[C]{
		\Gamma, p, p \vdash q
	}{
		\Gamma, p \vdash q
	}
\qquad
	\infrule[P]{
		\Gamma, p, q, \Delta \vdash r
	}{
		\Gamma, q, p, \Delta \vdash r
	}
$ \\
~\\
\hline
\end{tabular}
\end{center}
\end{small}
\caption{Structural inference rules for Grounded Arithmetic (\ga)}
\label{tab:ga:struct-rules}
\end{table}

To handle background assumptions explicitly,
we first incorporate in our system the standard structural rules
shown in \cref{tab:ga:struct-rules}.
The symbols `$\Gamma$' and `$\Delta$' in these rules
represent arbitrary lists of premises of any length zero or greater.
The \emph{hypothesis} rule \irl{H}
allows us to infer, without any prior premises,
that a conclusion $p$ is trivially true
if we already assumed $p$ in the hypotheses.
The \emph{weakening} rule \irl{W} 
allows us to add hypothetical assumptions,
producing a weaker statement with the new assumption $p$
from a stronger statement that was proved without assuming $p$.
We often need to weaken deductions in this way
in order to get the list of premises to ``line up with''
those of other deductions appearing elsewhere in a proof --
which, unlike the weakened deduction, might have actually needed $p$.
The \emph{contraction} rule \irl{C} allows us to ``contract''
or deduplicate several copies of the same hypothetical assumption into one.
We typically need contraction
when we need to use a hypothetical assumption more than once in a proof.
Finally,
the \emph{permutation} rule \irl{P} allows us to permute or reorder
the list of background hypotheses arbitrarily.
We could dispense with this rule
if we consider the hypotheses to be an unordered collection
rather than an ordered list.

Finally,
we expand each of the rules in \cref{tab:ga:rules}
to allow for background assumptions.
To do so,
in each of a rule's premises or conclusion
in which neither an entailment symbol `$\vdash$'
nor a definition symbol `$\ldef$' appears,
we prepend `$\Gamma \vdash$' to that premise or conclusion.
To each premise
in which an entailment symbol `$\vdash$' already appears,
we prepend only `$\Gamma,$'.
All premises and conclusions of all inference rules
thereby become entailments or definitions,
and all resulting entailments allow for background assumptions
represented by $\Gamma$,
which may be arbitrary provided they are fixed
throughout a given use of the rule.
We could have written all the rules in \cref{tab:ga:rules}
already expanded in this fashion,
but leaving background assumptions and unnecessary entailments
out of the main formulation of \ga seemed preferable
for clarity and conciseness.
Merely to illustrate this expansion process,
the expanded versions of the three main rules for logical disjunction
are as follows:

\[
	\infrule[\lor I1]{
		\Gamma \vdash p
	}{
		\Gamma \vdash p \lor q
	}
\qquad
	\infrule[\lor I2]{
		\Gamma \vdash q
	}{
		\Gamma \vdash p \lor q
	}
\qquad
	\infrule[\lor E1]{
		\Gamma \vdash p \lor q
	\qquad
		\Gamma, p \vdash r
	\qquad
		\Gamma, q \vdash r
	}{
		\Gamma \vdash r
	}
\]

\baf{structure of \ga proofs: lists or derivation trees}

\baf{how to derive the remaining inference rules presented earlier}

\baf{reducing disjoint types to 0=false,1=true,2+nat}

\baf{reducing natural number computations}

\subsection{Alternative formulations}

\ga could readily be reformulated differently from the formulation above
in a variety of ways without changing its essence.
We could treat conjunction as primitive rather than disjunction,
but then the bidirectional rule \irl{boolIE} would be less natural.
Similarly, we could treat universal
rather than existential quantification as primitive.
We could further minimize the formulation,
as measured in terms of number of inference rules,
by simply expanding the equivalences for $p \jbool$ and $a \jnat$,
at the cost of some textual repetition
(especially in the last inference rule in \cref{tab:ga:rules})
and consequent loss of clarity.

More significantly and less purely cosmetically,
we could reformulate \ga to omit recursive definitions entirely.
We would then need to introduce addition and multiplication as primitives,
more closely along the lines of traditional formulations
of Peano Arithmetic (PA) or Heyting Arithmetic (HA).
Including addition and multiplication primitives
make it possible to express the construction of pairs
and other finite data structures encoded as natural numbers,
from which we can express and reason about Turing machines
and other general models of recursive computation.
Which of these alternatives is simpler or more ``foundational''
seems largely a matter of subjective taste.

%later

\section{A computational interpretation of \ga terms}
\label{sec:comp}

All \ga terms correspond in principle to ordinary computations.
Not only is \ga a logic \emph{of} or \emph{about} computation;
it is also a \emph{computable logic}.
Not only are the natural numbers and functions on natural numbers
computable in \ga,
but all of the boolean predicates \ga can express
concerning those natural numbers and functions are also computations.

In particular,
assuming \ga is formulated correctly,
any \ga term $t$ that provably has some value $v$
according to the proof rules above
corresponds to a computation that terminates and yields value $v$.
Similarly,
for any terminating computation expressible in \ga's term language --
which should be \emph{any} computation,
since \ga terms are Turing complete --
the fact that this computation terminates, and the value it yields,
should be provable via \ga's deduction system.

This section will not yet succeed
in rigorously proving this intended correspondence.
Instead, for now we merely begin to explore this relationship
in two ways.
First, we describe and further analyze \ga's computational behavior
via an \emph{operational semantics},
a now-standard tool for rigorous specification of programming languages.
Second, we will outline automatic transformations or \emph{reductions}
from \ga terms into more-conventional programming languages,
to illustrate more concretely and intuitively
how we conceive of \ga terms ultimately as ``just software.''

\subsection{A big-step structural operational semantics (BSOS) for \ga}
\label{sec:sos}

\begin{table}
\begin{small}
\begin{center}
\renewcommand*{\arraystretch}{0.5}	% make blank rows shorter
\begin{tabular}{|c|}
\hline
~\\
\textbf{Definition} \\
~\\
$
	\infrule{
		s(\vec{x}) \ldef \ttc{d}{\vec{x}}
	\qquad
		\tto{p}{\ttc{d}{\vec{a}}} \Downarrow v
	}{
		\tto{p}{\ttc{s}{\vec{a}}} \Downarrow v
	}
$ \\
~\\
\hline
~\\
\textbf{Natural numbers} \\
~\\
$
	\infrule{
	}{
		x \Downarrow x
	}
\qquad
	\infrule{
	}{
		0 \Downarrow 0
	}
\qquad
	\infrule{
		t \Downarrow n
	}{
		\suc(a) \Downarrow n+1
	}
\com{
\qquad
	\infrule{
	}{
		\tlambda{x}{a(x)}	\Downarrow	\tlambda{x}{a(x)}
	}
}%com
$ \\
~\\
\hline
~\\
\textbf{Equality} \\
~\\
$
	\infrule{
		a \Downarrow n
	\qquad
		b \Downarrow n
	}{
		a=b \Downarrow \ctrue
	}
\qquad
	\infrule{
		a \Downarrow n
	\qquad
		b \Downarrow m
	\qquad
		n \ne m
	}{
		a=b \Downarrow \cfalse
	}
$ \\
~\\
\hline
~\\
\textbf{Propositional logic} \\
~\\
$
	\infrule{
	}{
		\ctrue \Downarrow \ctrue
	}
\qquad
	\infrule{
	}{
		\cfalse \Downarrow \cfalse
	}
\qquad
	\infrule{
		p \Downarrow \ctrue
	}{
		\neg p \Downarrow \cfalse
	}
\qquad
	\infrule{
		p \Downarrow \cfalse
	}{
		\neg p \Downarrow \ctrue
	}
$ \\
~\\
$
	\infrule{
		p \Downarrow \ctrue
	}{
		p \lor q \Downarrow \ctrue
	}
\qquad
	\infrule{
		q \Downarrow \ctrue
	}{
		p \lor q \Downarrow \ctrue
	}
\qquad
	\infrule{
		p \Downarrow \cfalse
	\qquad
		q \Downarrow \cfalse
	}{
		p \lor q \Downarrow \cfalse
	}
$ \\
~\\
\hline
~\\
\textbf{Predicate logic} \\
~\\
$
	\infrule{
		\tto{p}{n} \Downarrow \ctrue
	}{
		\texists{x}{\tto{p}{x}} \Downarrow \ctrue
	}
\qquad
	\infrule{
		x \jnat \vdash \neg \tto{p}{x}
	}{
		\texists{x}{\tto{p}{x}} \Downarrow \cfalse
	}
$ \\
~\\
\hline
\end{tabular}
\end{center}
\end{small}
\caption{Big-step Structural Operational Semantics (BSOS) for \ga}
\label{tab:sos}
\end{table}

\Cref{tab:sos} concisely presents
a \emph{big-step structural operational semantics},
or BSOS,
for \ga terms.\later{\cite{XXX}}
The reduction rules in this table specify inductively
how more complex \ga terms may reduce to simpler ones,
with some \ga terms eventually reducing to a concrete natural number $n$
or a boolean constant \ctrue or \cfalse.
We next examine particular (sets of) reduction rules in detail.

\subsubsection{Definition reduction rules}

The rule for definitions in \cref{tab:sos}
essentially describes beta substitution for recursive functions,
using a fixed background set of recursive function definitions
instead of the unnamed lambda expressions
typical of functional programming languages.
This rule essentially states that
an invocation of a defined function symbol $s$
with a certain list of parameter subterms $\vec{a}$
reduces to some value $v$
whenever the code resulting from expanding that definition --
from $\ttc{s}{\vec{a}}$ to $\ttc{d}{\vec{a}}$ --
likewise reduces to the same value $v$.
The inductive character of the BSOS reduction rules
implicitly allows this beta substitution to be performed
any number of times to support recursive definitions.

This definition-reduction rule is formulated here
so as to allow beta reduction
essentially anywhere in a \ga term, without restriction.
We could alternatively specify a more constrained set of \emph{contexts}
in which such beta reductions are allowed,
as is common in specifying the operational semantics
of more traditional programming languages.
While such a context refinement should be feasible and may be worthwhile,
the uses we intend for this BSOS do not appear to make contexts necessary,
and hence would only appear to make the BSOS unnecessarily more complex.

\subsubsection{Natural number reduction rules}

\ga terms allow expression of natural numbers
in only a few simple ways:
via variables, the constant $0$,
or the successor function $\suc(a)$.
The reduction rules for natural numbers are correspondingly simple,
effectively only reducing any concrete representation of a natural number
having the form $\suc(\dots\suc(0)\dots)$
to a corresponding natural-number result $n$.

Variables not representing a concrete natural number
simply reduce to themselves, and not to a concrete value.
The BSOS thus effectively supports computation
only on concrete natural numbers,
and will typically ``get stuck'' and not reduce if presented with
an abstract-interpretation problem containing unknown variables.
There are cases where \ga terms containing unknown variables
may nevertheless reduce, however,
as we will see below.

\baf{variable reduction rule is probably unnecessary}

\subsubsection{Equality reduction rules}

Whenever two subterms $a$ and $b$
reduce to the same concrete natural number $n$,
the positive-case equality reduction rule
reduces `$a = b$' to \ctrue.
The corresponding negative-case rule
similarly reduces `$a = b$' to \cfalse
whenever $a$ and $b$ reduce to two unequal natural numbers.

Since these rules apply only when the subterms $a$ and $b$
both reduce to concrete natural numbers, however,
equality testing is effectively strict in both arguments:
the equality operator never reduces at all
if either subterm never reduces
(or reduces to something other than a concrete natural number,
such as a variable or a boolean).

\subsubsection{Propositional logic reduction rules}
\label{sec:sos:prop}

The reduction rules for logical negation `$\neg$' are unsurprising:
`$\neg p$' reduces to \cfalse if $p$ reduces to \ctrue, and
`$\neg p$' reduces to \ctrue if $p$ reduces to \cfalse.
If $p$ never reduces to \ctrue or \cfalse,
then `$\neg p$' never reduces to anything.

The three reduction rules for logical disjunction `$\lor$'
directly encodes the usual truth table for a disjunction operator.
In addition, these rules express
the fact that disjunction in \ga is \emph{non-strict} in both arguments:
`$p \lor q$' may reduce to \ctrue
even if $p$ never reduces (provided $q$ reduces to \ctrue),
or even if $q$ never reduces (provided $p$ reduces to \ctrue).
For example, `$\ctrue \lor x$' reduces to \ctrue
even though the variable $x$ never reduces to any concrete value.
Only in the event that neither $p$ nor $q$ ever reduce to \ctrue,
the conjunction operator effectively become strict,
ultimately yielding \cfalse only provided that
\emph{both} $p$ and $q$ reduce to \cfalse.
If neither $p$ nor $q$ ever reduce to any boolean constant, for example,
then `$p \lor q$' will ``get stuck'' and never reduce to anything.

Since logical conjunction, implication, and biconditional
are all considered derived rather than primitive in \ga,
the operational semantics of these derived operators follow in turn
from the operational semantics of logical negation and disjunction
as shown in \cref{tab:sos}.
Logical conjunction `$\land$', in particular,
unsurprisingly behaves as the dual of disjunction,
obeying the same reduction rules
only with \ctrue and \cfalse exchanged.

\subsubsection{Predicate logic reduction rules}
\label{sec:sos:pred}

The BSOS includes two reduction rules for \ga's existential quantifier:
one that under certain conditions
reduces an existentially-quantified term to \ctrue,
the other sometimes reducing an existentially-quantified term to \cfalse.

The positive-case reduction is not particularly surprising:
it specifies that if we can find some concrete natural number $n$ that,
when substituted into the quantified predicate $p$
causes $\tto{p}{n}$ to reduce to \ctrue,
then the quantified predicate `$\texists{x}{\tto{p}{x}}$'
likewise reduces to \ctrue.
We are essentially just encoding the unbounded search
for a positive-existence example $n$ satisfying the predicate
into the positive-case reduction rule for the predicate.

The reduction rules in this BSOS up to this point
may seem reasonably familiar and consistent with standard practice.
The negative-case reduction rule for the existential quantifier,
however,
has a more peculiar, noteworthy, and perhaps suspicious feature.
The precondition for this reduction rule,
`$x \jnat \vdash \neg \tto{p}{x}$',
refers back to \ga's \emph{deduction system}
as described earlier in \cref{sec:ga}.
The negative-case reduction rule for the existential quantifier
is thus explicitly demanding a \ga \emph{proof}
that `$\neg \tto{p}{x}$' holds and yields \ctrue
whenever the free variable $x$ represents any natural number.

This means that to implement a complete and correct evaluator
for all \ga terms including quantifiers,
the evaluator must effectively have a \ga theorem prover embedded within it.
While unconventional,
this property of \ga's BSOS is not a problem in principle,
since a theorem prover for any ordinary symbolic formal system like \ga
is itself just a computation --
a piece of software that can be implemented
in any Turing-complete computational framework.
Theorems of \ga, like those of most other interesting formal systems,
are recursively enumerable.
We merely rely on this fact in \ga's BSOS
by allowing ourselves the freedom to include operational reduction rules
that in effect search for \ga proofs by recursive enumeration.

Since we are at present interested in the theoretical properties of \ga,
we do not demand that the embedded \ga theorem prover
be \emph{practical} or \emph{efficient}.
We require only that it, given unlimited time and storage resources,
in principle would perform its designated task
of finding a particular \ga proof if one exists.
We make no claim that \ga terms represent \emph{practical computations},
only that they represent
what theorists often call \emph{effective computations} --
which can be and often are entirely impractical.

The upshot is that \ga's existential quantifier
can \emph{computationally} reduce either to \ctrue or to \cfalse,
in the former case as the result of a successful unbounded search
for a satisfying concrete natural number $n$,
in the latter case as the result of a successful unbounded search
for a \ga proof that no satisfying natural number can possibly exist.
We will further unpack and explore the meaning and implications
of these slightly-peculiar semantics in following sections.

If we can make ourselves comfortable
with the operational semantics of \ga's existential quantifier,
then the semantics of the universal quantifier --
which we treat as a non-primitive
derived from the existential quantifier --
follow immediately and are at least no more problematic.
In particular, \ga's universal quantifier is a precise dual
to the existential quantifier,
as in classical logic though unlike intuitionistic logic.

\subsubsection{Implications of the BSOS's dependence
		on \ga's deduction system}
\label{sec:sos:dep-impl}

There is one important and immediate practical consequence
of this effective embedding of a \ga theorem prover
into the reduction rules for quantifiers.
We are effectively making the BSOS, and many of its interesting properties,
directly dependent on the properties of the corresponding deduction system.

\later{
This means for example that the operational semantics of 
\emph{constructive} \ga (\cga) is different from that of full \ga{} --
not because of any difference immediately apparent in \cref{tab:sos},
but because the BSOS for \cga refers to a different deduction system
(that of \cga instead of \ga).
In particular,
\cga's deduction system allows only a subset of \ga's deduction rules
to be used in satisfying the proof
that the existential quantifier's negative-case reduction rule demands.
}%later

This dependence on the deduction system means that
for many properties we might like to prove about \ga's BSOS,
we will first have to prove related properties about \ga's logic.
For example, we might be tempted to use this BSOS
to prove \ga term evaluation deterministic --
\ie that a given term only ever reduces to at most one value.
To achieve this determinism proof in the case of the quantifiers, however,
we find that we first have to prove properties of \ga's logic
tantamount to (and likely stronger than) proving the logic itself consistent.
While proving interesting properties
of a conventional programming language's operational semantics
is traditionally simpler and easier than proving similar properties
of a full logical deducation system,
in the case of \ga's BSOS, we cannot expect this rule of thumb to hold
because of this dependency.

As a result, we cannot reasonably expect this BSOS to help us much
in proving many interesting, deep properties
like \ga's determinism or consistency.
We will instead attack problems like these
using other approaches and tools
developed in later sections.
We rely on this BSOS for now only to shed light intuitively
on the correspondence between reasoning and computation in \ga.

		% operational semantics
\subsection{The PCF and PPF programming languages}
\label{sec:pcf-ppf}

In order to obtain a potentially clearer and more intuitive
understanding of \ga's semantics from a computational perspective,
we will next examine \ga in relation to
more conventional programming languages.
We thus need a more conventional programming language or two
to relate \ga to.
For this purpose we choose 
\emph{Programming Computable Functions} or PCF,\footnote{
	See \cite{plotkin77lcf}.
}
a simple functional programming language
designed for theoretical analysis tasks such as this,
and 
\emph{Programming Parallel Functions} or PPF,
a minor extension to PCF supporting basic parallel computation.
We outline these two languages together
because they are so closely related.

\later{
%\subsection{PPF: Programming Parallel Functions}
%\label{sec:ppf}

This section defines PPF (Programming Parallel Functions),
a simple programming language
based on PCF (Programming Computable Functions) by Scott and Plotkin~\cite{XXX}.
PPF extends PCF in only one way:
by adding a rudimentary parallel-evaluation construct:
$a \parallel b$ evaluates subterms $a$ and $b$ concurrently,
nondeterministically returning the result of either subterm
if both subterms eventually terminate.
(In practical implementations,
this construct would typically launch two child threads
and wait for the result of whichever completes first --
but since we will not care about any notion of time
or of one thread completing ``before'' another,
it is simplest to consider this choice simply as nondeterminism:
\ie we care only that it is \emph{possible}
for $a \parallel b$ to return either subterm's result.)
If only one subterm terminates,
$a \parallel b$ always yields the result of that subterm.
if neither subterm terminates,
$a \parallel b$ in turn fails to terminate.

\baf{	explicitly define extenxed ``PPF with boolean type''
	as well and specify the (easy) reduction,
	for clarity and completeness? }
}

In its original formulation,
PCF includes both natural number and booleans as primitive types.
Since the booleans are straightforward to emulate via natural numbers,
subsequent formulations of PCF
often omit the boolean type and corresponding operations.\later{XXX}
We likewise omit the primitive boolean type here for simplicity.

\subsubsection{Term syntax in PCF and PPF}

Any computation in PCF or PPF is expressed as a term,
having the following syntax:

\begin{align*}
	T &\equiv
		x \mid
		0 \mid
		\suc(T) \mid
		\pred(T) \mid
		\tifz{T}{T}{T} \mid
		(\tlambda{x}{T}) \mid
		T (T) \mid
		\ky (T)
	& \text{(PCF)}
\\
	T &\equiv
		x \mid
		0 \mid
		\suc(T) \mid
		\pred(T) \mid
		\tifz{T}{T}{T} \mid
		(\tlambda{x}{T}) \mid
		T (T) \mid
		\ky (T) \mid
		(T \parallel T)
	& \text{(PPF)}
\end{align*}

Intuitively:
$x$ is a variable reference,
$0$ represents the natural-number constant 0,
$\suc$ is the natural number successor function,
$\pred$ is the natural number predecessor (clamping at zero).
An ``if-zero'' term $\tifz{a}{b}{c}$ first evaluates subterm $a$ and tests its result:
if $a$ evaluates to zero then $\tifz{a}{b}{c}$ evaluates subterm $b$ and returns its result
(the ``\kwstyle{then}'' case);
if $a$ evaluates to a nonzero number
then $\tifz{a}{b}{c}$ evaluates subterm $c$ and returns its result
(the ``\kwstyle{else}'' case).
A lambda term $(\tlambda{x}{t(x))}$
yields a higher-order function defined by subterm $t(x)$,
parameterized by the bound variable $x$.
A term $a(b)$ evaluates subterms $a$ and $b$ then,
if $a$ yields a function (lambda term),
applies that function to subterm $b$'s result,
yielding any result that the function produces.
A term $\ky(a)$ represents the fixed-point combinator,
invoking a function represented by subterm $a$ with a function parameter
representing the function's own return value,
allowing the function to invoke itself recursively.

The only syntax unique to PPF is `$a \parallel b$',
which intuitively launches two child processes or threads
that evaluate subterms $a$ and $b$ in parallel,
returning the result of any subterm that successfully evaluates,
as defined more precisely below.

\subsubsection{Types in PCF and PPF}

PCF and PPF use essentially the same type system,
in which we define types inductively as follows:

\begin{itemize}
\item	\tnat is a type: the basic type of natural numbers.
\item	if $\sigma$ and $\tau$ are types,
	then `$\sigma \to \tau$' is a type,
	representing functions from $\sigma$ to $\tau$.
\end{itemize}

The sole modification to PCF's type system required by PPF
is the addition of a typing rule
for the parallel construct operator $\parallel$
to verify that both subterms are of the same type
(which in turn becomes the type of the parallel construct).

We will not be particularly concerned with types here, however,
since types do not affect the language's evaluation behavior
and are not required to define its operational semantics,
which is all we will need for now.
(A denotational semantics could be defined for PPF,
representing the denotations of terms as sets of possible traces,
for example --
but we will not have particular need of such a denotational semantics
and so will omit one here.)

\later{summarize types and type judgments}

\subsubsection{Operational semantics of PCF and PPF}

\begin{table}
\begin{small}
\begin{center}
\renewcommand*{\arraystretch}{0.5}	% make blank rows shorter
\begin{tabular}{|c|}
\hline
~\\
\textbf{Reduction rules common to PCF and PPF} \\
~\\
$
	\infrule{
	}{
		x \Downarrow x
	}
\qquad
	\infrule{
	}{
		0 \Downarrow 0
	}
\qquad
	\infrule{
	}{
		(\tlambda{x}{a(x)})	\Downarrow	(\tlambda{x}{a(x)})
	}
$ \\
~\\
$
	\infrule{
		t \Downarrow n
	}{
		\suc(t) \Downarrow n+1
	}
\qquad
	\infrule{
		t \Downarrow 0
	}{
		\pred(t) \Downarrow 0
	}
\qquad
	\infrule{
		t \Downarrow n+1
	}{
		\pred(t) \Downarrow n
	}
$ \\
~\\
$
	\infrule{
		a \Downarrow 0
	\qquad
		b \Downarrow v
	}{
		\tifz{a}{b}{c} \Downarrow v
	}
\qquad
	\infrule{
		a \Downarrow n+1
	\qquad
		c \Downarrow v
	}{
		\tifz{a}{b}{c} \Downarrow v
	}
$ \\
~\\
$
	\infrule{
		a \Downarrow (\tlambda{x}{e(x)})
	\qquad
		\subs{b}{x}{e} \Downarrow v
	}{
		a(b) \Downarrow v	
	}
\qquad
	\infrule{
		a (\ky (a)) \Downarrow v
	}{
		\ky (a) \Downarrow v	
	}
$ \\
~\\
\hline
~\\
\textbf{Reduction rules exclusive to PPF} \\
~\\
$
	\infrule{
		a \Downarrow v
	}{
		(a \parallel b) \Downarrow v
	}
\qquad
	\infrule{
		b \Downarrow v
	}{
		(a \parallel b) \Downarrow v
	}
$ \\
~\\
\hline
\end{tabular}
\end{center}
\end{small}
\caption{Big-step Structural Operational Semantics (BSOS) of PCF and PPF}
\label{tab:pcf-ppf:sos}
\end{table}

\Cref{tab:pcf-ppf:sos} summarizes
the big-step operational semantics of PCF and PPF,
expressing more precisely in terms of reduction rules
the informal behavior of the constructs summarized above.

A complete and fully-rigorous operational semantics
also normally includes an inductive definition
of all the possible \emph{evaluation contexts} 
within which the reductions may be performed.
It is essential that reductions be allowed on subterms
embedded within other not-yet-reduced terms.
For example,
the first argument of a \kifz must reduce to a natural number
before we can test this result for zero
and allow the surrounding \kifz in turn to be reduced.
We merely specify informally here that the allowed evaluation contexts
are as usual for PCF,
and in PPF,
reductions are allowed within both the $a$ and $b$ subterms
of a parallel composition construct `$a \parallel b$'.

Despite its minimality,
the fixed-point combinator \ky enables PCF terms
to express arbitrary recursive computations.
PCF is therefore Turing-complete,
unlike primitive-recursive arithmetic
or the simply-typed lambda calculus for example,\later{cite}
It is thus easy to express PCF computations that never terminate.
The PCF language is deterministic, however:
any term reduces to \emph{at most one} concrete value
(a natural number or function),
or else never reduces to any value (and hence denotes $\bot$).\footnote{
	This determinism property is readily provable
	using either PCF's operational or denotational semantics.
	\baf{cite}
}

PPF, in contrast with PCF, is obviously nondeterministic in general.
For example, the program $0 \parallel \suc(0)$
can nondeterministically yield either of the natural numbers 0 or 1.
In the further reasoning steps below,
we will be most interested in \emph{using} PPF
in ways that nevertheless produce deterministic results.
Actually achieving this goal, and proving that we have achieved it,
will be an important and nontrivial challenge.

\later{	down-conversion of boolean terms,
	letrec definitions}

\subsubsection{Metacircular evaluation and simulation of \pcf and \ppf}
\label{sec:pcf-ppf:eval}

As with any Turing-complete computational model,
we can simulate the execution of either PCF or PPF
atop any other Turing-complete computational model --
including atop PCF and PPF themselves, in particular.
For example,
we could readily construct a metacircular evaluator for PCF terms
within the PCF language,
or similarly construct a metacircular evaluator for PPF terms
within PPF.
Either of these constructions would require
encoding terms into natural numbers;
this is a pragmatically tedious
but now theoretically standard and uncontroversial practice.

This basic principle of computability theory extends
to the simulation of parallel languages like PPF
atop sequential or deterministic languages like PCF.
If we were to construct a metacircular evaluator for PPF terms
atop PCF,
in particular,
then the only nontrivial challenge
is to simulate PPF's parallel computation terms
of the form `$a \parallel b$'.

Suppose we have a metacircular evaluator for PCF terms
taking the form of a function $E(s,n)$
which takes as parameters a natural-number step count $s$
and a natural number $n$ that uniquely encodes some PCF term $t$.
That is, $n$ is a G\"odel code for term $t$:
in the common ``Quine quote'' notation,
$n = \quo{t}$.
This PCF metacircular evaluator $E$
executes the PCF program represented by term $t$ for at most $s$ steps.
If $t$ terminates and yields a concrete value $v$
(a natural number or lambda expression),
then $E(s,\quo{t})$ returns an encoding of that result value,
which we assume to be nonzero.
If $t$ does not terminate within $s$ steps,
then $E(s,\quo{t})$ returns 0 to indicate as such.
Since this function evaluates the PCF term for a bounded number of steps $s$,
the function $E$ is not just computable but primitive recursive.

Using this PCF evaluator $E$,
we can readily simulate PPF's `$a \parallel b$' parallel construct atop PCF
by running two instances of $E$ in parallel,
one on $\quo{a}$ and the other on $\quo{b}$,
returning any result that either simulated subterm reduces to
at the first step $s$ in which either of them reduces.
In the case of a ``tie'' in which both subterms reduce to concrete values
on exactly the same simulated step $s$,
we arbitrarily return the result from $E(s,\quo{a})$.
In other words,
we execute subterms $a$ and $b$ as if they were two threads or processes
running in a virtual time schedule that alternates between
executing thread $a$ and thread $b$.

In nondeterministic-execution cases where both subterms $a$ and $b$
both eventually reduce to a concrete value,
this simulation will render PPF's parallel composition operator
quasi-deterministic,
effectively choosing one of the results
that PPF's operational semantics permits,
fairly arbitrarily based on the virtual ``timing'' subtleties
of how the two subterms ``race'' to produce their respective results.
These semantically-nondeterministic cases
are not those we will be primarily interested in, however.

The situations we will be interested in are those in which,
based on some knowledge about the behavior of either one of the subterms,
we can infer that
it \emph{does not matter} what the other subterm reduces to --
or whether the other subterm reduces at all.
As a particularly important example,
if we happen to know that $a$ never reduces to any concrete value
(and hence denotes $\bot$),
then the only value that `$a \parallel b$' can ever reduce to in \ppf
is whatever value $b$ reduces to, if anything.
Similarly, if $b$ never reduces at all,
then the only value `$a \parallel b$' can possibly reduce to
is any result from $a$.
If we can ensure that inferences of these kinds
apply to \emph{all} of the uses we make
of \ppf's parallel composition operator,
then we can deduce that our \emph{use} of parallel evaluation in \ppf
remains deterministic in effect,
even if \ppf is a semantically nondeterministic computational model.

%later
		% PCF and PPF languages
\subsection{Reductions from \ga to PPF and PCF computations}
\label{sec:reduct}

Having briefly summarized
the syntax and operational semantics of \pcf and \ppf,
we now turn towards using them as tools to help understand and analyze \ga.

First, although we omitted the primitive boolean type
from the above formulations of PCF and PPF,
let us assume that we have encoded \ctrue as 1 and \cfalse as 0
and constructed versions of the usual boolean operations
as computations on natural numbers in PCF.
For convenience and simplicity
we will thus subsequently pretend we have the customary language facilities
that traditionally produce or consume booleans:
\eg we will pretend that we have in PCF and PPF an `$=$' operator
that compares two natural numbers and yields a boolean,
although this is in fact a function that compares two natural numbers
and yields the natural number 1 or 0.
Similarly, we will pretend that we have in PCF and PPF the familiar
`$\tif{c}{a}{b}$' construct expecting its condition $c$ to be boolean,
although in fact this construct is just a 3-argument curried PCF function
testing whether its $c$ argument is 1 or 0
and returning the result of $a$ or $b$, respectively.

Thus conceptually augmented,
PCF already appears to include \emph{most}
of the computational capabilities
that \ga does,
according to its inference rules in \cref{tab:ga:rules}
and its operational semantics in \cref{tab:sos}.
PCF can certainly express and compute
both functions of natural numbers,
and function-predicates as functions that return either 0 or 1.
The fixed-point combinator \ky gives \pcf
essentially the same ability to express recursive functions:
\pcf has \emph{greater} expressiveness than \ga, in fact,
if we count \pcf's ability to compute using higher-order functions
(and not just with a fixed set of recursive definitions
as we formulated \ga for simplicity).

There are really only two computational elements in \ga
that are ``new'' and special in \ga and not already in \pcf:
namely \ga's logical disjunction operator `$a \lor b$',
and its existential quantifer `$\texists{x}{\tto{p}{x}}$'.
We now examine each of these operators in turn,
and how we may ultimately reduce them to \ppf or pcf equivalents
while preserving their semantics.

\subsubsection{Reducing non-strict logical disjunction `$a \lor b$'
		in \ga to \ppf or \pcf}
\label{sec:reduct:or}

It is easy to write even in \pcf a function that computes
the logical disjunction `$a \lor b$'
of two booleans $a$ and $b$ encoded into natural numbers,
and yields another boolean encoded as a natural number:
we could compute
$1-(1-\quo{a})(1-\quo{b})$, for example.
This implementation leverages the fact that
we can implement `$\neg b$' as $1-\quo{b}$,
we can implement `$a \land b$' as $\quo{a} \cdot \quo{b}$,
and $a \lor b = \neg(\neg a \land b)$ by De Morgan's laws.

The remaining problem is that the obvious ways
to implement logical disjunction (or conjunction) in \pcf
using either arithmetic, or conditional \kif statements,
yield behavior that is \emph{strict} in at least one if not both arguments.
That is, the arithmetic calculation $1-(1-a)\times(1-b)$
will yield a result only if \emph{both} $a$ and $b$ have been computed
and reduced to concrete natural numbers 0 or 1.
We could implement disjunction using an \kif conditional,
as in either `$\tif{\quo{a}=1}{1}{\quo{b}}$'
or `$\tif{\quo{b}=1}{1}{\quo{a}}$'.
The first of these last two options will be strict in the first argument
but non-strict in the second,
in that the second argument need not reduce at all
if the first argument reduces to 1.
But still neither of these alternatives
will successfully match the semantics of disjunction in \ga,
in which `$a \lor b$' reduces to \ctrue provided only that
\emph{either one} of the subterms $a$ and $b$ reduce to \ctrue,
even if the other subterm never reduces at all.

Let us focus now on reducing \ga's disjunction operator
to \ppf instead of \pcf,
because such a reduction to \ppf is simple and intuitive.

First, assume we have defined a function in \ppf
that implements the same \emph{guarded term} idiom
that we defined for \gd in \cref{sec:quant:guard},
which we can define as follows:

\[
	p \oq a \ldef \tif{p}{a}{\bot}
\]

That is, `$p \oq a$' first tests $p$,
returns the result of $a$ if $p$ evaluates to \ctrue,
and otherwise deliberately enters an infinite loop,
never yielding any result.

Given this shorthand for guarded terms,
we can readily implement \ga's non-strict logical disjunction
by reducing it to the following equivalent PPF term:

\[
	a \lor b \ldef	(a \oq \ctrue) \parallel (b \oq \ctrue) \parallel
			(\neg a \oq \neg b \oq \cfalse)
\]

In essence, our \ppf implementation of logical disjunction
essentially launches three independent threads:
the first tests $a$ and returns \ctrue if $a$ evaluates to \ctrue;
the second tests $b$ concurrently 
and returns \ctrue if $b$ evaluates to \ctrue.
Finally, the third thread tests both $a$ and $b$ in sequence
and returns \cfalse only if both evaluate to \cfalse.

This reduction essentially just amounts to a parallelized implementation
of logical conjunction corresponding to
Kleene's 3-valued ``strong logic of indeterminacy''.\later{cite}

Notice that provided that subterms $a$ and $b$
each represent at least \emph{effectively} deterministic computations --
\ie computations that can reduce to at most one concrete value
even when implemented in a nondeterministic language like \ppf{} --
then this implementation of `$a \lor b$'
preserves this effective determinism.
In order for `$a \lor b$'
to reduce to both \ctrue and \cfalse nondeterministically,
either the first and third, or the second and third,
parallel ``child threads'' launched above would have to yield a result.
But the first and third child threads cannot both complete
because of the assumption that $a$ and $\neg a$ cannot both reduce to \ctrue.
Similarly, the second and third child threads cannot both complete
provided that $b$ and $\neg b$ cannot both reduce to \ctrue.

\later{Note somewhere that Plotkin's formulation of PCF
	basically contained a ``parallel or'' extension...}

Having first reduced \ga's non-strict logical disjunction operator
to parallel evaluation in \ppf,
we can further reduce this parallel evaluation to \pcf, if desired,
using the simulation techniques discussed above in \cref{sec:pcf-ppf:eval}.
Provided the subterms $a$ and $b$ are effectively deterministic
as described above,
the PCF simulation of `$a \lor b$' will reduce to at most one concrete value
that does not depend on messy subtleties like the relative ``virtual timing''
of the various \ppf threads being simulated.

\subsubsection{Reducing the existential quantifier in \ga
		to \pcf or \ppf}
\label{sec:reduct:quant}

We now turn to the slightly tricker but conceptually similar challenge
of reducing an existential quantifier in \ga,
of the form `$\texists{x}{\tto{p}{x}}$',
to \ppf and ultimately to \pcf.

We can reduce the behavior of \ga's existential quantifier
(and by duality, its non-primitive universal quantifier)
into two parts,
one specialized to evaluating the quantifier to \ctrue,
the other specialized to evaluating it to \cfalse.
In effect,
we can view \ga's full ``two-sided'' existential quantifier,
which can evaluate to either to either \ctrue or \cfalse,
as the parallel combination of two ``one-sided'' quantifiers --
a one-sided existential quantifier `$\exists^{+}$'
and a one-sided universal quantifier `$\forall^{+}$' --
each of which only ever evaluates either to \ctrue or not at all ($\bot$).
In this way, 
we can decompose \ga's two-sided quantifiers
into one-sided quantifiers in \ppf as follows:

\begin{align*}
\texists{x}{\tto{p}{x}}	&\ldef	\texistsp{x}{\tto{p}{x}} \parallel
				\neg \tforallp{x}{\neg \tto{p}{x}} \\
\tforall{x}{\tto{p}{x}}	&\ldef	\tforallp{x}{\tto{p}{x}} \parallel
				\neg \texistsp{x}{\neg \tto{p}{x}}
\end{align*}

Focusing on the former case,
the two-sided existential quantifier `$\exists$' 
effectively launches two parallel child threads in \ppf.
The first, ``true case'' child thread
invokes the one-sided existential quantifier `$\exists^{+}$' 
to search for some concrete natural number $x$
that makes the predicate $\tto{p}{x}$ evaluate to true.
The second, ``false case'' child thread, in parallel,
invokes the one-sided universal quantifier `$\forall^{+}$'
to search for a \ga \emph{proof}
that no such natural number satisfying the predicate exists:
\ie that for \emph{all} natural numbers $x$,
the predicate $\tto{p}{x}$ evaluates to false.
If this latter search succeeds
and the one-sided universal quantifier `$\forall^{+}$'returns \ctrue,
then the second child thread negates that result
and returns \cfalse from the original two-sided existential quantifier.

The above implementation of \ga's universal quantifier
does exactly the same,
only swapping the uses of the two one-sided quantifiers.
It thus returns \ctrue
if the first child thread
successfully uses the one-sided universal quantifier
to prove that $\tto{p}{x}$ is \ctrue for all natural numbers $x$,
and returns \cfalse
if the second child thread
successfully uses the one-sided existential quantifier
to find some natural number $x$ for which $\tto{p}{x}$ yields \cfalse,
thus serving a counterexample refuting the universal quantifier's claim.

Having split the positive and negative cases of each quantifier in this way,
how do we then implement the two \emph{one-sided} quantifiers we still need?

\subsubsection{Implementing the one-sided existential quantifier `$\exists^{+}$'}
\label{sec:reduct:explus}

We can implement the one-sided existential quantifier `$\exists^{+}$'
via a recursive function that makes further use of parallel composition
to express a parallel unbounded search through all possible values of $x$
for any natural number $x$ that satisfies the predicate.
Assume we have a general-recursive metacircular evaluator $E\qb{t}$
that evaluates \ppf term $t$,
terminating with some result or not according to $t$'s behavior in \ppf.

\begin{align*}
\texistsp{x}{\tto{p}{x}}	&\ldef	\exists_0^{+}{x}\ \tto{p}{x} \\
\exists_n^{+}{x}\ \tto{p}{x}	&\ldef	\tto{p}{n} \oq \ctrue \parallel
					\exists_{n+1}^{+}{x}\ \tto{p}{x}
\end{align*}

In effect,
the one-sided existential quantifier
first launches two child threads:
one to test the zero case,
the other to test all higher cases starting from 1.
The first child thread tests $\tto{p}{0}$,
returning true if that test succeeds,
but never yielding any result at all if this test fails.
The second child thread recursively does the same,
launching a child to test the case $x=1$,
the other child thread devoted to testing all cases greater than $1$,
and so on.
Notice that by construction,
`$\exists^{+}$' can only ever return true or nothing at all:
there is no execution path by which it can ever return false.

If we knew that the predicate $\tto{p}{x}$ would always terminate
and yield \emph{some} boolean result of true or false,
then in principle the one-sided existential quantifier
could be simplified to dispense with
this forking of an unbounded number of child threads
to perform a parallel search for a natural number
satisfying the predicate $p$.
In that case,
we could instead just perform a sequential search upwards
through the natural numbers starting with 0.
We do \emph{not} wish to assume
that the predicate $p$ always terminates, however.
If $p$ doesn't always terminate,
then a sequential search would get ``stuck''
at the first natural number for which $p$ fails to terminate,
never getting around to testing larger natural numbers
for which $p$ might again terminate.
With this unbounded parallel search construction, in contrast,
we can guarantee that the one-sided existential quantifier
will terminate and yield true
if $\tto{p}{x}$ evaluates to \ctrue for any value of $x$,
even if for all \emph{other} values of $x$
$\tto{p}{x}$ never terminates at all.
In essence, 
the one-sided existential quantifier `$\texistsp{x}{\tto{p}{x}}$'
is constructed to perform non-strict evaluation
across all possible values of $x$.

\later{	Footnote: discuss the fact that Kleene's
	development of recursion theory introduced
	a similar parallel-search operator complementary to $\mu$
	that will never get stuck if the predicate is not total,
	but is not guaranteed to find the *least* satisfying value.
}

\subsubsection{Implementing the one-sided universal quantifier `$\forall^{+}$'}
\label{sec:reduct:allplus}

As discussed earlier in \cref{sec:sos:pred},
the false case of the two-sided existential quantifier,
and hence the one-sided universal quantifier,
essentially relies on a theorem prover
embedded in the evaluation logic.

Assume that we have a \ga proof-checker function $C(n_P,n_j)$
that returns \ctrue (\ie 1)
exactly when $n_P$ is the \gdl code for a valid \ga proof $P$,
$n_j$ is the \gdl code for a judgment `$H \vdash p$',
and the final judgment in proof $P$ is `$H \vdash p$'.
Since this function checks only a single proof
and performs no unbounded search,
it is primitive recursive and always terminates.

We can then implement the one-sided universal quantifier `$\forall^{+}$'
as follows:

\begin{align*}
\tforallp{x}{\tto{p}{x}}	&\ldef	\forall_0^{+}{x}\ \tto{p}{x} \\
\forall_n^{+}{x}\ \tto{p}{x}	&\ldef	\tif{C(n,\quo{x \jnat \vdash \tto{p}{x}})}{\ctrue}{\forall_{n+1}^{+}{x}\ \tto{p}{x}}
\end{align*}

This implementation just uses a standard \kif conditional and recursion
to search sequentially through all possible encoded proofs
for one that correctly proves the desired result.
We do not need \ppf's parallel composition in this case
because the proof checker $C$ is primitive recursive
and guaranteed to terminate:
it can thus never ``get stuck'' forever
checking a possible proof.

In summary,
we have seen how \ga's logical disjunction operator
and existential quantifier --
the only two computational features of \ga
that are not already in ordinary programming languages such as \pcf\ --
may be reduced first to parallel computations in \ppf
and then, by standard simulation techniques,
to sequential functional computation in \pcf.
We have not yet proven that these constructs in fact
``do'' or express what we want, or anything reasonable for that matter,
which will be the continuing task of the next section.

%later

\section{Denotational semantics of \ga}
\label{sec:domain}

Having informally sketched a broad relationship
between \ga terms and computations in the prior section,
we now focus on deeper analysis of
only a fragment of the full \ga language:
namely the quantifier-free \emph{base \ga} or \bga subset.

Since \bga omits the existential and universal quantifiers,
its operational semantics can similarly omit
reduction rules for the quantifiers.
This change mitigates the complications
highlighted in \cref{sec:sos:dep-impl},
and likely renders it more feasible
to prove interesting properties of \bga
in terms of its operational semantics.

For now, however,
we will instead analyze \bga using \emph{denotational semantics},
where we assign semantic meaning according to term structure,
rather than by induction over reduction steps.
One attractive feature of denotational semantics for our purposes
is that it readily yields determinism and consistency proofs.
A correct denotational semantics in essence automatically attaches
one and only one mathematical object as the ``meaning'' of any language term.
If \bga has a denoational semantics and a predicate $p$ denotes \ctrue,
then $p$ cannot simultaneously denote \cfalse (making \bga inconsistent)
because each term denotes one and only one mathematical object.

\later{	cite "Design Concepts in Programming Languages"
	as an example of using basic functional programming theory
	as a foundation for more sophisticated languages}

\subsection{Types and domains: borrowing from PCF}

For our purposes we will use the well-established tools of domain theory
to model and assign semantic meaning to \bga terms.
The domains we will need are in fact a subset of those
that Scott and Plotkin used to model
\emph{Programming Computable Functions} or PCF,
a simple statically-typed functional programming language.\footnote{
	See \cite{plotkin77lcf} for Plotkin's formulation of PCF,
	which was in turn based on Scott's influential
	but long-unpublished logic for computable functions or LCF
	(see \cite{scott93type}).
	For a more recent exploration
	of domain theory and functional programming foundations
	focusing on PCF,
	see \cite{streicher06domain}.
}

PCF employs a type system with types defined inductively as follows:

\begin{itemize}
\item	\tbool is the type of boolean values \ctrue and \cfalse.
\item	\tnat is the type of natural numbers 0, 1, etc.
\item	For any types $\sigma$ and $\tau$,
	`$\sigma \to \tau$' is a type representing computable functions
	that map values of type $\sigma$ to values of type $\tau$.
\end{itemize}

Associated with each of these types is a Scott domain
that we can define inductively alongside the structure
of the type it represents.
In particular:

\begin{itemize}
\item	PCF type \tbool uses the domain $\sbool_{\bot}$,
	the flat domain of boolean truth values.
\item	PCF type \tnat uses the domain $\snat_{\bot}$,
	the flat domain of natural numbers.
\item	Given domains $D_{\sigma}$ and $D_{\tau}$
	associated with PCF types $\sigma$ and $\tau$, respectively,
	$D_{\sigma \to \tau}$
	is the domain of Scott-continuous functions
	from $D_{\sigma}$ to $D_{\tau}$.
\end{itemize}

From the above domains
we finally construct a domain $D_V$
representing \bga's entire value space.
$D_V$ is simply the squashed disjoint union
of the above boolean, natural number, and function domains,
with the respective bottom elements `$\bot$' identified
(hence ``squashed'' together)
while keeping all other elements disjoint.

\later{illustrate}

We will use the notation `$\sigma^k \to \tau$' to refer to
the PCF type of a $k$-argument ``curried'' function
that returns functions in order to handle multiple arguments.
That is:

\[
\begin{array}{rcl}
	\sigma^0 \to \tau	&\leqv&	\tau \\
	\sigma^1 \to \tau	&\leqv&	\sigma \to \tau	\\
	\sigma^2 \to \tau	&\leqv&	\sigma \to (\sigma \to \tau)	\\
	\sigma^3 \to \tau	&\leqv&	\sigma \to (\sigma \to (\sigma \to \tau)) \\
			&\vdots&
\end{array}
\]

In fact the main function types we will need for \bga
are those specifically of the form `$\tnat^k \to \tbool$'.
To any $k$-parameter predicate symbol $s$ defined in \bga
via a definition of the form `$s(\vec{x}) \ldef \ttc{d}{\vec{x}}$' --
where $k = |\vec{x}|$ --
we assign the PCF type `$\tnat^k \to \tbool$'.
Constant symbol definitions for the case $k=0$
thus have a PCF type of \tbool,
single-parameter defined symbols
have PCF type `$\tnat \to \tbool$', 
and so on.

\subsubsection{Fixed points and recursive definitions}
\label{sec:domain:fixpoint}

Since a \bga definition may be recursive, however,
to assign meaning to the body of a recursive definition
we will need PCF's fixed-point combinator \ky.
The \ky combinator has PCF type $(\sigma \to \sigma) \to \sigma$
for any type $\sigma$.
The \ky combinator in essence ``invokes'' its input function $f$,
of type $\sigma \to \sigma$,
feeding $f$'s output (of type $\sigma$)
back as the argument to $f$,
so that $f$ can refer recursively to the output value that it returns.

In the context of PCF,
fixed points are typically useful
only when $\sigma$ is in turn a function type:
the \ky combinator then enables that function to invoke itself recursively.
The \ky combinator nevertheless works on all types $\sigma$
including \tbool,
as we will need if we express a constant symbol recursively,
such as the Liar paradox `$L \ldef \neg L$'.

In \bga, the denotational semantic meaning
of a $k$-ary \bga definition `$s(\vec{x}) \ldef \ttc{d}{\vec{x}}$'
is defined by the fixed point
`$\ky\tlambda{s\vec{x}}{M_A\qb{\ttc{d}{\vec{x}}}}$',
where $M_A\qb{\ttc{d}{\vec{x}}}$
is the meaning of the definition's expansion $\ttc{d}{\vec{x}}$
as defined below,
under an assignment $A$ that binds the recursively-defined symbol $s$
and its $k$ formal arguments $\vec{x}$
to the $\lambda$ expression's actual parameters.

%later

\subsection{Semantics of \bga terms and formulas}
\label{sec:domain:semantic}

Using the above type system,
we assign the type \tnat to any term $t$ in \bga's syntax,
and we assign the type \tbool to any formula $f$ in this syntax.
In a $k$-argument predicate-function application `$f(t_1,\dots,t_k)$',
each of the arguments $t_1,\dots,t_k$ has type \tnat,
and the function $f$ must have type `$\tnat^k \to \tbool$'
in order for the application to be well-typed.
(We can assume that the types of arguments to $\lambda$ terms are inferred,
or alternatively we could add explicit types
to the arguments of $\lambda$ terms as in PCF.)

We define an \emph{assignment} $A$
as a function from variable names ($v_1,v_2,\dots$)
and symbol names ($s_1,s_2,\dots$)
to elements of the domain $D_V$ representing \bga values.
We will assume that variables and defined symbols
have disjoint namespaces and are thus distinct from each other,
although an assignment
can map both variables and symbols to domain elements.
A \emph{variable assignment} $V$
maps only variable names to elements of $D_V$,
while a \emph{symbol assignment} $S$
maps only symbol names to elements of $D_V$.
Any assignment $A$ may be viewed as a composition of (or decomposed into)
a variable assignment $A_V$ and a symbol assignment $A_S$.

We can now express a semanting meaning function $M_A\qb{a}$
that assigns domain-theoretic semantics to term $a$ under assignment $A$,
as follows:

\begin{itemize}
\item	$M_A\qb{v}$ is $\lift{n}$ if $A(v)$ is $\lift{n}$
	for some natural number $n$,
	and is $\bot$ otherwise.
\item	$M_A\qb{0}$ is $\lift{0}$.
\item	$M_A\qb{\suc(a)}$ is $\lift{n+1}$
	if $M_A\qb{a}$ is $\lift{n}$,
	and $\bot$ otherwise (\ie if $M_A\qb{a}$ is $\bot$).
\end{itemize}

When applied to a term $t$,
the meaning function $M_A(t)$ clearly yields values
only in the subdomain $\snat_{\bot}$
of domain $D_V$ representing \bga's value space.
%\ie either a lifted natural number $\lift{n}$ or $\bot$.

Next, we assign a semantic meaning $M_A\qb{f}$
to any \bga formula $f$ as follows:

\begin{itemize}
\item	$M_A\qb{a \jnat}$ is $\ltrue$
	if $M_A\qb{a}$ is $\lift{n}$ for some natural number $n$,
	and $\bot$ otherwise.
\item	$M_A\qb{a=b}$ is $\ltrue$
	if $M_A\qb{a}=\lift{n}$ and $M_A\qb{b}=\lift{n}$
	for some natural number $n$.
	Otherwise,
	$M_A\qb{a=b}$ is $\lfalse$
	if $M_A\qb{a}=\lift{n_a}$ and $M_A\qb{b}=\lift{n_b}$
	for natural numbers $n_a \ne n_b$.
	Otherwise,
	$M_A\qb{a=b}$ is $\bot$.
\item	$M_A\qb{p \jbool}$ is $\ltrue$
	if $M_A\qb{p}$ is either $\ltrue$ or $\lfalse$,
	and $\bot$ otherwise.
\item	$M_A\qb{\neg p}$ is
	$\ltrue$ if $M_A\qb{p}$ is $\lfalse$,
	$\lfalse$ if $M_A\qb{p}$ is $\ltrue$,
	and $\bot$ otherwise.
\item	$M_A\qb{p \lor q}$ is $\ltrue$
	if either $M_A\qb{p}$ or $M_A\qb{q}$ is $\ltrue$.
	Otherwise,
	$M_A\qb{p \lor q}$ is $\lfalse$
	if both $M_A\qb{p}$ and $M_A\qb{q}$ are $\lfalse$.
	Otherwise,
	$M_A\qb{p \lor q}$ is $\bot$.
\com{
\item	$M_A\qb{\tlambda{\vec{x}}{\ttc{p}{\vec{x}}}}
	is the Scott-continuous curried function taking arguments $\vec{x}$
	and yielding $M_A\qb{\ttc{p}{\vec{x}}}$.
\item	$M_A\qb{\ky p}$
	is the fixed point of $p$,
	as described above.
}%com
\item	$M_A\qb{s(a_1,\dots,a_k)}$
	is $f(M_A(a_1),\dots,M_A(a_k))$
	if assignment $A$ maps $s$ to $f$
	where $f$ is a value of PCF type `$\tnat^k \to \tbool$',
	and $\bot$ otherwise.
	In the constant-definition case $k=0$,
	$M_A\qb{s}$ is $\lift{b}$
	if $A(s)=\lift{b}$ for boolean $b$,
	and $\bot$ otherwise.
\com{
	is the Scott-continuous function $M_A\qb{f}$
	applied to the arguments $M_A(a_1),\dots,M_A(a_k)$,
	each of type $\tnat$.
	Recall that the $f$ here always has the form
	`$\ky\tlambda{(f,\vec{x})}{\ttc{d}{\vec{x},f(\vec{x})}}$'
	due to the replacement of defined symbols
	with equivalent recursive lambda predicates as discussed above.
}%com
\end{itemize}

\baf{	we could use letrec notation if that would be clearer?}

\com{
In the last case, \ky is the fixed-point combinator for type $\tau$,
which passes function $f$ to $f$ itself as $f$'s first parameter,
yielding a general-recursive function of type $\tau$
taking exactly $k$ arguments from domain $\snat_{\bot}$
and yielding a result from domain $\sbool_{\bot}$.
}%com

We have thus assigned domain-theoretic semantic meaning under assignment $A$
to each of the terms and formulas in \bga,
the quantifier-free restricted syntax of \ga.

\subsection{Proving the consistency of \bga}
\label{sec:domain:consist-bga}

We now have the machinery necessary to prove \bga consistent.
We do so in the classic model-theoretic fashion,
by proving inductively over the length of quantifier-free \bga proofs
that each \bga inference rule ``preserves truth''
under any variable assignment.

We first assume a fixed symbol assignment $S$
representing the set of background definitions in effect.
We will compose this fixed symbol assignment $S$
with non-fixed variable assignments $V$ as described below
to form the composite assignments $A = S \cup V$
required as arguments to the semantic meaning function $M_A$.

We say that a judgment `$H \vdash p$' is \emph{valid}
if for all variable assignments $V$
mapping variable names to elements of the domain $\snat_{\bot}$,
yielding a composite assignment $A = S \cup V$,
if $M_A\qb{H_i} = \ltrue$ for each $1 \le i \le |H|$,
then $M_A\qb{p} = \ltrue$.
That is, a judgment holds if any assignment
that satisfies all hypotheses in $H$
also satisfies the consequent $p$,
quantified over all possible variable assignments
that may be composed with the fixed symbol assignment $S$.

To prove \bga consistent,
we must handle the inference rules for definitions, equality,
natural numbers, recursive computation, and propositional logic,
as shown in \cref{tab:ga:rules}.

We first address the conditional bidirectional rule \irl{{\ldef}IE}
for definitions and recursive computation.
Our fixed symbol assignment $S$
maps each $k$-ary defined symbol $s$
to a value of the PCF domain for type $\tnat^k \to \tbool$,
and that domain value is the semantic meaning
of the definition's body $\ttc{d}{\vec{x}}$
as discussed above in \cref{sec:domain:fixpoint}.
When operated in either direction,
the rule \irl{{\ldef}IE}
substitutes a predicate-function call with that definition's body
or vice versa.
The meaning $M_A\qb{s(\vec{a})}$ of the predicate-function call
with argument terms $\vec{a}$ is,
by the domain-theoretic properties of function application,
is equal to the meaning $M_a\qb{\ttc{d}{\vec{a}}}$
of the corresponding definition's body
with the (meanings of) the arguments $\vec{a}$ substituted accordingly.
As a result, by structural induction over the predicate $p$
that the predicate-function invocation $s(\vec{a})$
or corresponding definition body $\ttc{d}{\vec{a}}$ are embedded in,
the denotation of the overall predicate $p$
is unaffected by substitution in either direction,
and hence the inference rule preserves truth.

%\subsubsection{Inference rules for equality}

We next prove that the inference rules for equality
in \cref{tab:ga:rules}
preserve truth as follows:

\begin{itemize}
\item	\irl{{=}S}:
	Assuming that a judgment of the form `$H \vdash a=b$' is valid,
	we must show that the judgment `$H \vdash b=a$' is valid.
	Suppose that some assignment $A$ satisfies all hypotheses in $H$.
	Then by the inference rule's premise and our induction hypothesis,
	$A$ satisfies formula `$a=b$'.
	By the domain-theoretic semantics above,
	this event ($M_A\qb{a=b} = \ltrue$) occurs only if
	$M_A\qb{a} = \lift{n}$ and $M_A\qb{b} = \lift{n}$
	for some natural number $n$.
	But then by the same semantic rule,
	$M_A\qb{b=a} = \ltrue$,
	so this inference rule preserves truth under $A$.

\item	\irl{{=}E}:
	Assuming the judgments
	`$H \vdash a=b$' and `$H \vdash \tto{p}{a}$' are valid,
	we must show that `$H \vdash \tto{p}{b}$' is valid.
	By the rule's first premise and the semantics of equality,
	there is a natural number $n$
	such that $M_A\qb{a} = \lift{n}$ and $M_A\qb{b} = \lift{n}$.
	We can then prove by induction on the structure of $p$
	that $M_A\qb{\tto{p}{a}} = M_A\qb{\tto{p}{b}}$:
	\ie replacing $a$ with $b$ in $p$
	does not affect $p$'s denotation.
	Since $M_A\qb{\tto{p}{a}} = \ltrue$
	by the inference rule's second premise,
	$M_A\qb{\tto{p}{b}} = \ltrue$ as well,
	thereby satisfying the rule's conclusion.

\item	\irl{{\ne}S}:
	Assuming `$H \vdash a \ne b$' is valid,
	we must show that `$H \vdash b \ne a$' is likewise valid.
	For any assignment $A$ satisfying all hypotheses in $H$,
	$M_A\qb{a \ne b} = \ltrue$
	only if there exist natural numbers $n_a$ and $n_b$
	such that $M_A\qb{a} = \lift{n_a}$ and $M_A\qb{b} = \lift{n_b}$
	and $n_a \ne n_b$.
	But then $M_A\qb{b \ne a} = \ltrue$ as well,
	thereby satisfying the inference rule's conclusion.

\item	\irl{{\ne}IE}:
	Since this is a bidirectional inference rule,
	we must show that `$H \vdash a \ne b$' is valid
	if and only if `$H \vdash \neg(a=b)$' is valid.
	Reasoning in each direction in turn
	given some assignment $A$ satisfying all hypotheses in $H$:
	\begin{itemize}
	\item	$M_A\qb{a \ne b} = \ltrue$
		only if there exist natural numbers $n_a$ and $n_b$
		such that $M_A\qb{a} = \lift{n_a}$ and $M_A\qb{b} = \lift{n_b}$
		and $n_a \ne n_b$.
		Then $M_A\qb{a = b} = \lfalse$,
		and by the semantics of logical negation,
		$M_A\qb{\neg(a=b)} = \ltrue$,
		thereby satisfying the rule's conclusion.
	\item	$M_A\qb{\neg(a = b)} = \ltrue$
		only if $M_A\qb{a = b} = \lfalse$,
		implying in turn that
		there exist natural numbers $n_a \ne n_b$
		such that $M_A\qb{a} = \lift{n_a}$ and $M_A\qb{b} = \lift{n_b}$.
		But then $M_A\qb{a \ne b} = \ltrue$,
		thereby satisfying the rule's conclusion.
	\end{itemize}
\end{itemize}

%\subsubsection{Inference rules for natural numbers}

We next address the \bga inference rules for the natural numbers:

\begin{itemize}
\item	\irl{natIE}:
	When we apply this bidirectional inference rule
	in the forwards direction,
	$M_A\qb{H \vdash a=a} = \ltrue$
	only when there is a natural number $n$
	such that $M_A\qb{a} = \lift{n}$,
	but then $M_A\qb{H \vdash a \jnat} = \ltrue$ as well.
	Taking the rule in the reverse direction,
	the converse to the above reasoning holds.

\item	\irl{0I}:
	By the semantics above,
	$M_A\qb{0} = \lift{0}$,
	so $M_A\qb{0 \jnat} = \ltrue$,
	satisfying the inference rule's unconditional conclusion.

\item	\irl{\suc{=}IE}:
	In the inference rule's forward direction,
	if $M_A\qb{a=b} = \ltrue$,
	then there is a natural number $n$
	such that $M_A\qb{a} = \lift{n}$ and $M_A\qb{b} = \lift{n}$.
	But then by the semantics of successor $\suc$,
	$M_A\qb{\suc(a)} = M_A\qb{\suc(b)} = \lift{n+1}$,
	so $M_A\qb{\suc(a) = \suc(b)} = \ltrue$,
	thereby satisfying the inference rule's conclusion.
	Using the inference rule in the reverse direction,
	given $M_A\qb{\suc(a) = \suc(b)} = \ltrue$,
	there is an $n$ such that
	$M_A\qb{\suc(a)} = M_A\qb{\suc(b)} = \lift{n+1}$.
	But then $M_A\qb{a} = M_A\qb{b} = \lift{n}$,
	so $M_A\qb{a=b} = \ltrue$.

\item	\irl{\suc{\ne}0I}:
	Given $M_A\qb{a \jnat} = \ltrue$,
	there is a natural number $n$
	such that $M_A\qb{a} = \lift{n}$.
	By the semantics of successor $\suc$ and the $0$ term,
	$M_A\qb{\suc(a)} = \lift{n+1}$ and $M_A\qb{0} = \lift{0}$.
	Then since $n+1$ is never equal to zero for any natural number $n$,
	$M_A\qb{\suc(a) \ne 0} = \ltrue$,
	satisfying the rule's conclusion.

\item	\irl{\suc{\ne}IE}:
	The reasoning for this rule exactly mirrors
	that of \irl{\suc{=}IE} above,
	only for inequality rather than equality.

\item	\irl{Ind}:
	We assume that assignment $A$ satisfies
	all background hyptoheses in $H$.
	From the rule's third premise `$a \jnat$'
	and the semantics above
	there is a natural number $n_a$
	such that $M_A\qb{a} = \lift{n_a}$.
	We will use the rule's first two premises to prove
	by (metalogical) induction on $n_a$
	that $M_A\qb{\tto{p}{x}} = \ltrue$
	when $n_a$ as a literal natural number is substituted for $x$.

	In the base case $n_a = 0$,
	the rule's first premise dirctly ensures
	$M_A\qb{\tto{p}{0}} = \ltrue$.

	In the inductive step,
	we can assume $M_A\qb{\tto{p}{n}} = \ltrue$
	and must show $M_A\qb{\tto{p}{\suc(n)}} = \ltrue$.
	Using inference rules \irl{0I} and \irl{\suc{=}IE}
	inductively over $n$,
	we find that $M_A\qb{n \jnat} = \ltrue$.
	Using our assumption above
	that assignment $A$ satisfies the background hypotheses $H$
	and these last two results,
	we satisfy the hypotheses of the inference rule's second premise,
	and therefore can infer
	(by our overall induction hypothesis over proof length)
	that $M_A\qb{\ttc{q}{\ttmore}} = \ltrue$.
	This result then satisfies the inference rule's conclusion.
\end{itemize}

%\subsubsection{Inference rules for recursive computation}

\com{
We next address the conditional bidirectional rule \irl{{\ldef}IE}
for recursive computation.

Recall from \cref{sec:domain:defs} above
that in function calls of the form `$s(\vec{a})$' appearing in \bga terms,
we have replaced the defined symbol $s$
with an anonymous representation of the defined function itself,
in terms of a lambda predicate and the fixed-point operator \ky.
We can therefore rewrite the inference rule \irl{{\ldef}IE}
into the following form,
in which the definition is no longer needed as a prerequisite:

\[
	\infeqv[{\ldef}IE]{
%		s(\vec{x}) \ldef \ttc{d}{\vec{x}}
%	}{
		\ttc{p}{\ttc{d}{\vec{a}}}
	}{
		\ttc{p}{(\ky\tlambda{(f,\vec{x})}{\ttc{d}{\vec{x},f(\vec{x})}})\,(\vec{a})}
	}
\]

By of the semantics of the fixed-point combinator \ky
and function application:

\begin{align*}
	M_A\qb{s(\vec{a})}
&=	M_A\qb{(\ky\tlambda{(f,\vec{x})}{\ttc{d}{\vec{x},f(\vec{x})}})\,
		(\vec{a})} \\
&=	M_A\qb{(\tlambda{\vec{x}}{\ttc{d}{\vec{x},
		(\ky\tlambda{(f,\vec{x})}{\ttc{d}{\vec{x},f(\vec{x})})}\,
		(\vec{x})}})\,(\vec{a})} \\
&=	M_A\qb{\ttc{d}{\vec{a},
		(\ky\tlambda{(f,\vec{x})}{\ttc{d}{\vec{x},f(\vec{x})}})\,
		(\vec{a})}} \\
&=	M_A\qb{\ttc{d}{\vec{a},s(\vec{a})}} \\
\end{align*}

\baf{XXX not quite right but close.
	maybe easier to define the semantics in terms of PCF itself?}
}%com

We have thus proved that all of \bga's inference rules preserve truth,
yielding only valid conclusions from valid premises.
For any unconditional (hypothesis-free) judgment
of the form `$\vdash p$' provable in \bga,
in particular,
the denotational meaning $M_A\qb{p}$ of the conclusion $p$
must be $\ltrue$.
Since $M_A$ is a mathematical function
that can by definition assign only one meaning to a given formula $p$,
$M_A\qb{p}$ cannot also be \lfalse.
But if some \bga formula $p$ existed
for which both $p$ and $\neg p$ were provable,
then $M_A\qb{p}$ and $M_A\qb{\neg p}$ would both have to be $\ltrue$,
implying that $M_A\qb{p}$ would have to be both $\ltrue$ and $\lfalse$,
which is impossible because $M_A$ is a function.
We have thus proved that \bga at least is consistent.

\later{
Related work on relevant domains in constructive type theory,
including discussion and background references about LPO etc:
"The Scott model of PCF in univalent type theory"
\url{https://arxiv.org/pdf/1904.09810}
}

\subsection{The consistency of constructive \cga}
\label{sec:domain:consist-cga}

Having established the consistency
of the quantifier-free \bga system,
how can we extend our denotational semantics --
and consistency proof --
to the quantifiers?

The most important insights here
are that \bga can already express and reason about
arbitrary recursive (\ie Turing-complete) computations,
and that \ga's quantifiers are ``just computations''
as we explored already in \cref{sec:comp}.

As a result,
we need not actually \emph{extend} \bga
with quantifiers as additional primitives,
because the computations corresponding to our intended quantifiers
already exist as ordinary computations within \bga.
Instead,
we need only \emph{find} and suitably define
these quantifier computations within the existing framework of \bga.

Having defined the quantifiers suitably as computations,
we expect to find that
the constructive inference rules for \ga's quantifiers
(\cref{tab:ga:rules})
are \emph{admissible}:
they do not make any new theorems provable with respect to \bga,
but just offer more clear and convenient shortcuts
by which to prove theorems that could be proven anyway
without the quantifier rules.
\later{
For now we restrict our attention solely to constructive \ga or \cga,
excluding the non-constructive typing rules 
$\exists TI$ and $\forall TI$.
}

We will not attempt a rigorous consistency proof here,
but merely sketch an outline for such a proof,
yet to be completed and formally verified.

\subsubsection{Denotational semantics of parallel composition}

Our first step in reasoning about the quantifiers,
and perhaps the most tedious to perform with full technical rigor,
is to formulate a denotational semantics for
PPF's parallel composition operator `$\parallel$'.
This formulation can potentially be accomplished in various ways,
as long as we obtain the following three key properties of interest:

\begin{itemize}
\item	If $M_A\qb{b} = \bot$, then $M_A\qb{a \parallel b} = M_A\qb{a}$.
\item	If $M_A\qb{a} = \bot$, then $M_A\qb{a \parallel b} = M_A\qb{b}$.
\item	If $M_A\qb{a} = e$ and $M_A\qb{b} = e$,
	then $M_A\qb{a \parallel b} = e$.
\end{itemize}

One way to formulate parallel composition
is to implement it directly in PCF
in terms of a stepwise metacircular evaluator for PCF terms,
as outlined in \cref{sec:pcf-ppf:eval}.
We first implement and prove the key correctness properties
of a metacircular evaluator $E_s\qb{a}$
taking a step count $s$ and a \pcf term $a$,
which evaluates term $a$ for $s$ steps
and yields a result if $a$ terminates within $s$ steps.
We then implement parallel composition `$a \parallel b$'
by alternatively stepping
evaluations of subterms $a$ and $b$
yielding the first value that either produces.
Provided that at most one of $a$ or $b$ ever terminate
and yield a domain element other than $\bot$ --
or provided that $a$ and $b$ both terminate
but \emph{agree} on the value they produce --
we obtain the key properties we need,
despite the underlying complexity of the parallel simulation.
\later{mention this key agreement property earlier too}

An alternative way to formulate
the denotational semantics of parallel composition
would be to do so directly in terms of domain theory,
thus skipping the metacircular evaluation and parallel simulation.
We can model `$a \parallel b$'
as a Scott-continuous function $f(e_a,e_b)$
whose output has the property that
if value-domain elements $e_a$ and $e_b$ are consistent
(\ie have some upper bound in common),
then $e_a \sqcup e_b \sqsubseteq e_f$.
That is, the parallel composition function $f$
yields at least as much information as inputs $a$ and $b$ ``agree upon.''
We do not care what $f$ produces
when $a$ and $b$ produce inconsistent results,
as long as $f$ is Scott-continuous, \ie computable.
The result of such inconsistent inputs could be
an arbitrary one of the inputs,
or we could model such a result as a special ``top'' element $\top$
in a complete-lattice extension to domain theory.\later{XXX notes/cite}

\subsubsection{Denotational semantics of the quantifiers}

As outlined earlier in \cref{sec:reduct},
we consider each of \ga's full ``two-sided'' quantifiers
(which can evaluate either to \ctrue or \cfalse)
as a parallel composition of two ``one-sided'' quantifiers
$\exists^{+}$ and $\forall^{+}$,
each of which can only ever evaluate to \ctrue or nothing at all ($\bot$).
We must therefore define the denotational semantics
of these one-sided quantifiers,
and prove that their results are consistent with each other
so that their parallel composition yields useful results.

As discussed in \cref{sec:reduct:explus},
we view a one-sided existential quantifer `$\texistsp{x}{x}{\tto{p}{x}}$'
as a non-strict unbounded search for any natural number $x$
for which the predicate $\tto{p}{x}$ yields $\ltrue$.
The one property that is ultimately important
is that $\exists^{+}$ yields $\ltrue$ if such a natural number $x$ exists,
and yields $\bot$ if no such $x$ exists.

On the other hand,
as discussed in \cref{sec:reduct:explus},
we view a one-sided universal quantifer `$\tforallp{x}{x}{\tto{p}{x}}$'
as the unbounded search for a \emph{proof} in \bga'a deduction system
that every natural number $x$ satisfies predicate $\tto{p}{x}$.
We could formulate these semantics either
as an actual program that performs such a search,
or perhaps more directly in domain theory.
The key property we need, ultimately,
is that $\forall^{+}$ yields $\ltrue$ if such a universal proof exists,
and otherwise yields $\bot$.
Note that even though this represents an essential dependency
on \bga's deduction system,
it requires only the quantifier-free fragment of \bga,
relying only on the fact that free variables
are universally quantified implicitly.

Provided we have these two one-sided quantifiers as building blocks
and have proven their relevant properties,
it is a small further step to show that
\ga's two-sided existential and universal quantifiers
have the intended semantics (see \cref{sec:reduct:quant}).
Building on the semantics of \bga as defined above
and its resulting consistency proof,
a natural number $x$ satisfying a two-sided quantifier's predicate
cannot both exist and not exist.
This fact is then sufficient to ensure that
at most one of the complementary one-sided quantifiers
comprising a two-sided quantifier
can ever yield a value other than $\bot$,
thereby guaranteeing that the two-sided quantifiers are in turn consistent
and work as expected.

We can then use these established semantics of the two-sided quantifiers
to prove that their inference rules in \cref{tab:ga:rules}
preserve truth with respect to the underlying semantics of \bga.

%later

\section{Reflecting on \ga within \ga}
\label{sec:refl}

Given the similarities between \ga
and classical Peano Arithmetic (\pa),
it is natural to wonder whether and in what way Kurt \gdl's famous 
incompleteness theorems about arithmetic apply to \ga.
We now explore this question.

\com{
Although Kurt \gdl is most famous
for his two fundamental incompleteness theorems~\cite{XXX},
a contribution of his just as important as the theorems themselves
is a set of now-standard techniques
for \emph{coding} formal languages and proof systems
so that they can \emph{reflect} or ``reason about themselves.''
In this section we now adapt
some of \gdl's coding and reflection techniques to \ga,
in order to see where they lead us in this alternate logical universe.
}

\baf{maybe late this section is a good place to discuss the Berry paradox?}

\baf{consider using $S_n$ consistently instead of $L_n$
	for formal system instances.
	L is also used for languages, and the Liar paradox. }

\baf{consider metalogically treating the finite structures we use
	as identical to the natural numbers we encode them into.
	This is an alternative to the ``type disjointness'' option
	we discussed earlier, but it is a valid alternative.
	And it might be useful in simplifying the presentation and reasoning.
	If we do this, we probably need a section just on (finite) types
	laying out these two alternative perspectives.
	The main issue is how to handle quoting,
	distinguishing the use of a term or template as a literal
	versus an active part of the logic or computation.

	Might work if we clearly lay out the ``coded types'' CT discipline,
	and use quote-braces for literal quoting of syntax.
}

\subsection{Visualizing reflective reasoning via logic system instances}
\label{sec:refl:vis}

Before getting into further detail,
a visualization might help clarify in our minds what is happening
in the (reflective) study of logic in \gdl's fashion.

\begin{figure}[p]
\begin{center}
	\begin{subfigure}{0.48\textwidth}
	\begin{center}
	\includegraphics[width=\textwidth]{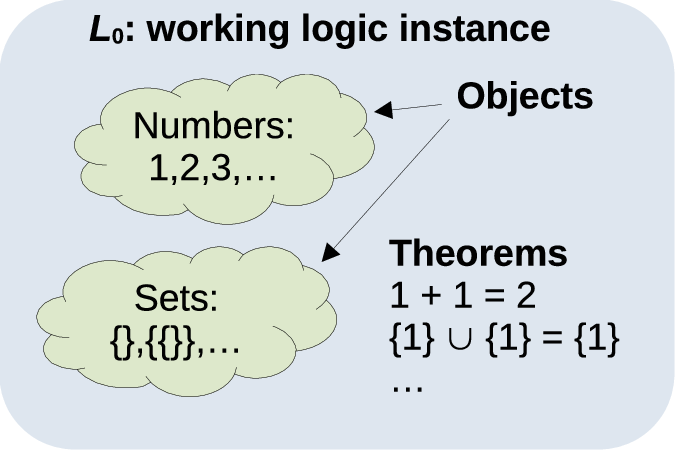}
	\end{center}
	\caption{Working use of a logic $L_0$
		to reason about ordinary objects such as numbers and sets.}
	\label{fig:refl:base}
	\end{subfigure}
\hfill
	\begin{subfigure}{0.48\textwidth}
	\begin{center}
	\includegraphics[width=\textwidth]{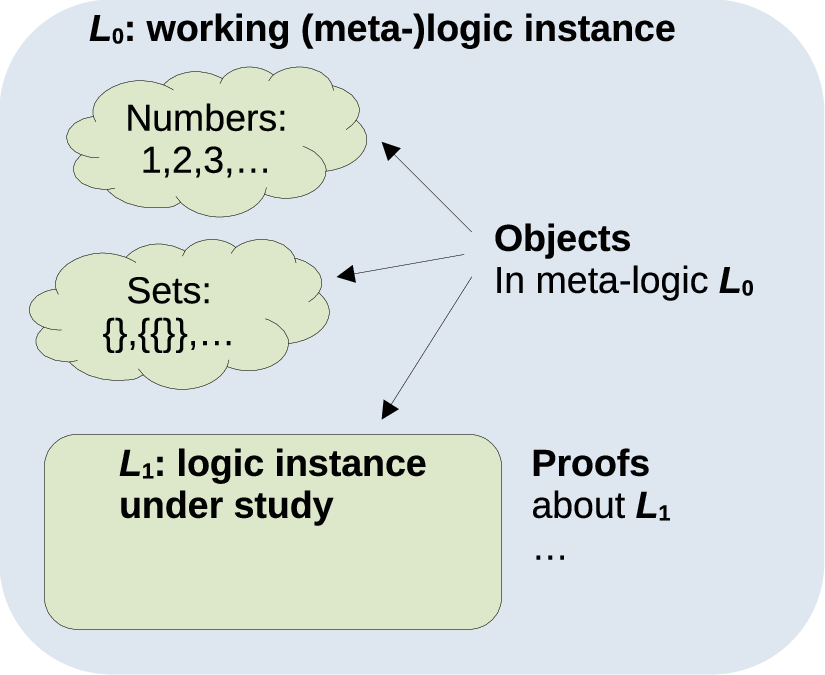}
	\end{center}
	\caption{Working within a meta-logic $L_0$
		to study a target logical system $L_1$.}
	\label{fig:refl:logic}
	\end{subfigure}
\end{center}

\begin{center}
	\begin{subfigure}{0.48\textwidth}
	\begin{center}
	\includegraphics[width=\textwidth]{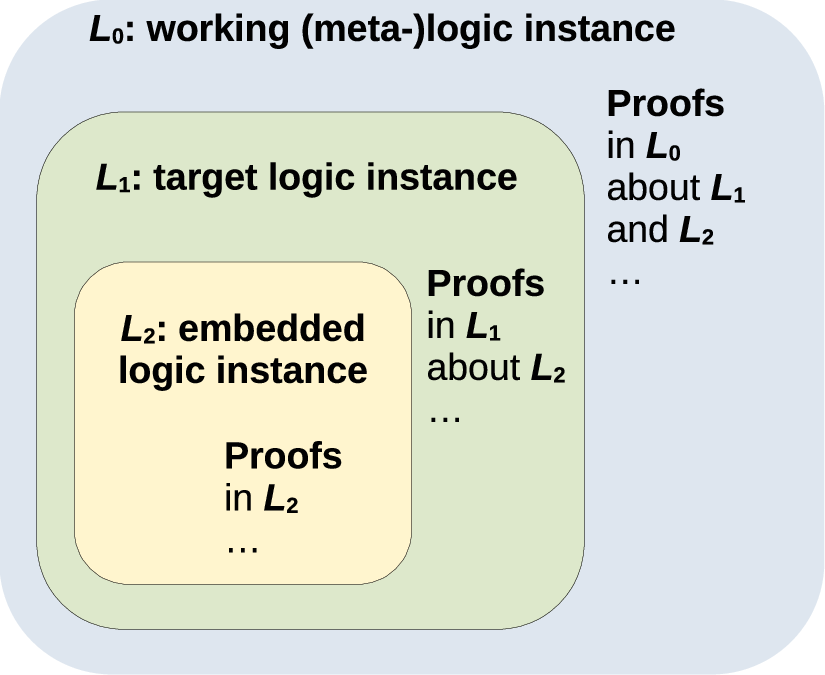}
	\end{center}
	\caption{Working within a meta-logic $L_0$
		to study a target logic system $L_1$,
		with an embedded logic system $L_2$ embedded in $L_1$.}
	\label{fig:refl:goedel}
	\end{subfigure}
\hfill
	\begin{subfigure}{0.48\textwidth}
	\begin{center}
	\includegraphics[width=\textwidth]{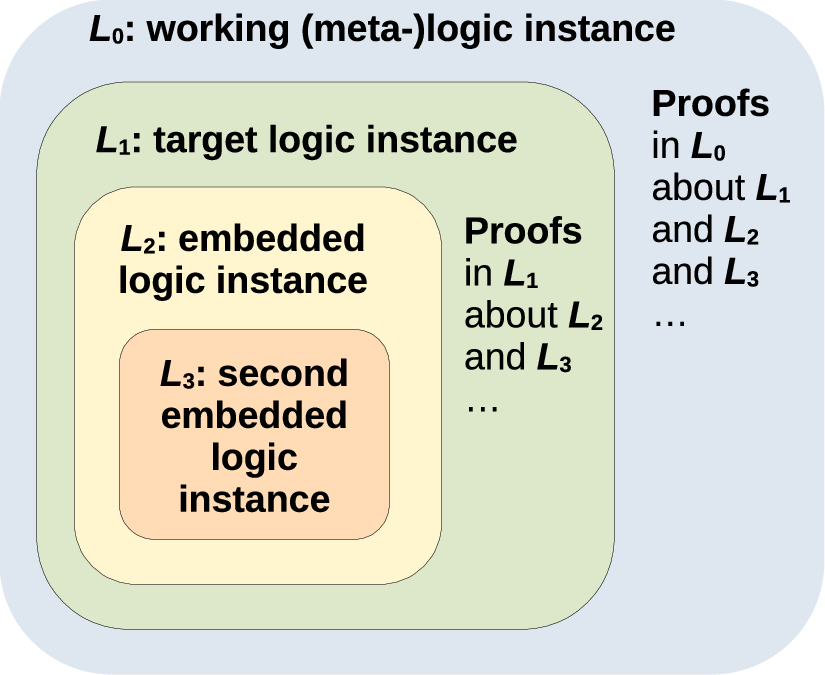}
	\end{center}
	\caption{Working within a meta-logic $L_0$
		to study a target logic system $L_1$,
		containing two embedded logic instances $L_2$ and $L_3$.}
	\label{fig:refl:goedel2}
	\end{subfigure}
\end{center}
\caption{Working use of a system of logic $L_0$
	to reason about ordinary objects such as numbers and sets.}
\label{fig:refl}
\end{figure}

\com{
\begin{figure}[t]
\begin{center}
\includegraphics[width=0.5\textwidth]{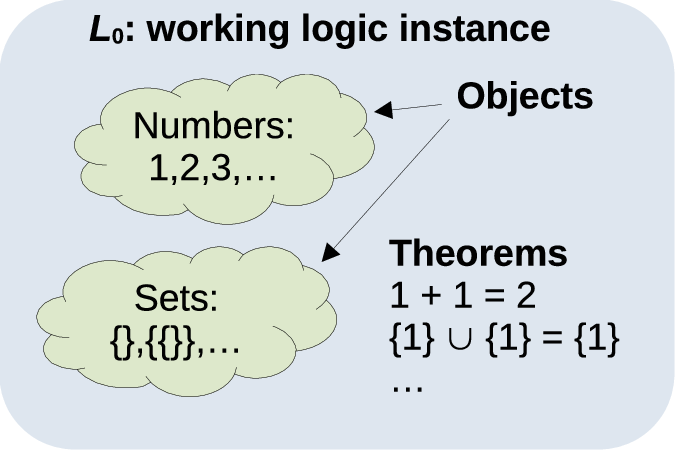}
\end{center}
\caption{Levels of reflection used in mathematical and logical reasoning.}
\label{fig:refl}
\end{figure}
}%com

First, consider as a base case our normal mathematical practice
of \emph{using} logical reasoning to study and make deductions about
other kinds of ``ordinary'' mathematical objects,
such as numbers or sets.
That is, in this base case
we are \emph{using} logic as a tool,
to study anything we like \emph{except for} systems of logic.
In this situation,
illustrated in \cref{fig:refl:base},
there is only one logical system in use:
our \emph{working logic}, which we will call $L_0$.
We might use $L_0$ to prove many things about numbers, sets, and so on,
but we are not trying study or prove anything \emph{about} logic.
This working system $L_0$ may constitute simply
the informal use of a human language (\eg English)
mixed casually with mathematical notation as needed
in writing or blackboard discussion;
$L_0$ may constitute the more rigorous but still only ``human-verified'' 
practice of traditional pencil-and-paper mathematical proof
as typically practiced in the peer-reviwed theory literature;
or $L_0$ may constitute the working use of
a mechanical proof verifier such as Coq or Isabelle,
which still more rigorously
insists that every detail be explicit and verifiably correct.
Regardless of the level of rigor in our mathematical practice,
the point for now is that when logic is \emph{not} an object of study,
there is only one logical system in play: $L_0$, the working system.

\com{
\begin{figure}[t]
\begin{center}
\includegraphics[width=0.5\textwidth]{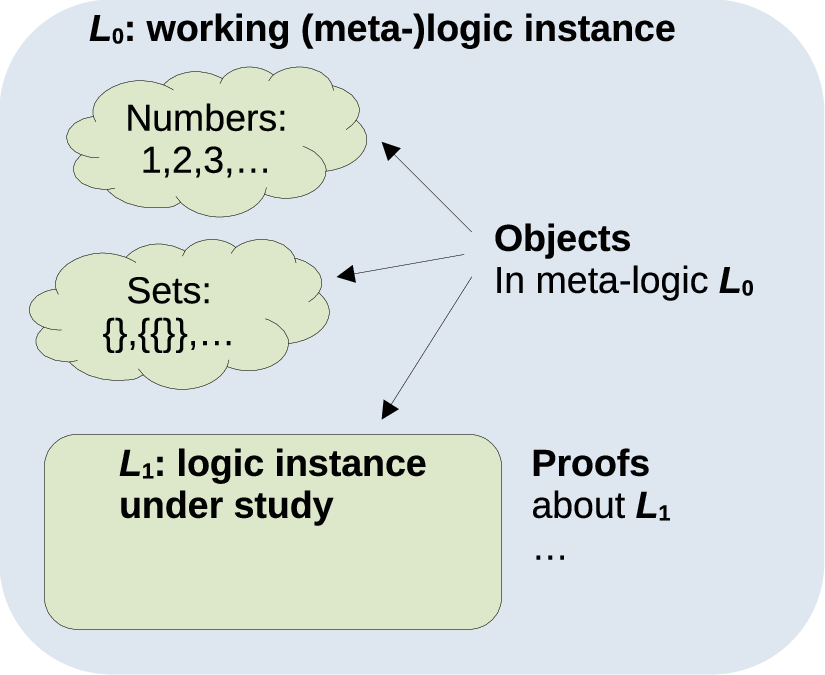}
\end{center}
\caption{Working within a meta-logic $L_0$
	to study a logical system $L_1$.}
\label{fig:refl:logic}
\end{figure}
}

\Cref{fig:refl:logic}, in contrast,
illustrates what a \emph{logician} is doing
when performing (basic) studies of a system of logic
\emph{as an object of study}. 
In this case, two \emph{instances} of logical systems are in play.
We will call the logical system actually being studied
the \emph{target} logic or $L_1$.
To reason \emph{about} logic $L_1$, however,
the logician still needs to \emph{use} a working logic
defining the rules and techniques the logician considers acceptable
as tools with which to reason \emph{about} the objects of study,
just as in all other mathematical practice.
As in the base case considered above,
we will label this working logic $L_0$.
We may also say that $L_0$ plays a role of \emph{meta-logic}
in defining the tools and rules for reasoning about target logic $L_1$.
The logician often still needs and uses ``ordinary'' mathematical tools
such as numbers and sets,
in particular as building blocks
with which to construct and reason about target logic $L_1$ --
but in this case the target logic $L_1$ itself
is the primary object of mathematical interest.

We say that $L_0$ and $L_1$ are two different logic \emph{instances}
because they play two distinct \emph{roles} --
that of meta-logic and target logic, respectively --
regardless of precisely which
language syntax, deduction rules, or axioms
we might choose for $L_0$ and $L_1$.
In particular, we might pick a single formal system $F$,
defining a particular language syntax and set of deduction rules --
\eg Peano arithmetic (PA) or Zermelo-Fr\"ankel set theory (ZF) --
and decide to use this \emph{same} system $F$
in both the meta-logic role of $L_0$
and in the target logic role of $L_1$.
In this case,
it is accurate and common, but a bit sloppy,
to say that we are using system $F$ to ``study itself.''
What we are actually doing in this case is
to use one \emph{instance} of system $F$,
in a role of meta-logic $L_0$,
as a tool to study a \emph{another logic instance} --
namely target logic instance $L_1$.
We construct and analyze target logic $L_1$ \emph{within} meta-logic $L_0$
using the ordinary mathematical tools that $L_0$ provides --
and it just so happens that we know
that both instance $L_0$ and instance $L_1$
use identical language and rules as defined by system $F$.
We usually ``know'' that $L_0$ and $L_1$ follow the same rules of system $F$
only to the extent that
(a) we are careful in constructing target logic $L_1$ within $L_0$
to ensure that this model $L_1$ of $F$ is indeed faithful to 
the definitive specification of $F$, wherever that is; and
(b) we are careful to follow the \emph{discipline} defined by system $F$
when using meta-logic instance $L_0$ to construct and reason about $L_1$.
We may informally \emph{know} that
$L_1$ is just another instance of the same logic as $L_0$,
but can we \emph{prove} this fact? If so, how?
That is not so easy.

\com{
\begin{figure}[t]
\begin{center}
\includegraphics[width=0.5\textwidth]{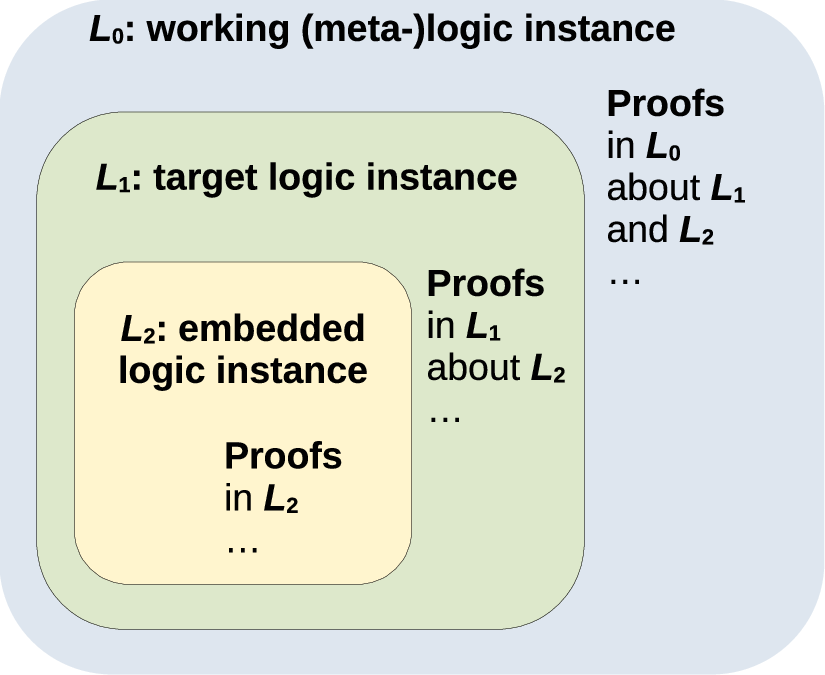}
\end{center}
\caption{Working within a meta-logic $L_0$
	to study a target logic system $L_1$,
	with an embedded logic system $L_2$ embedded in $L_1$ via coding,
	as done by \gdl in his incompleteness theorems.}
\label{fig:refl:goedel}
\end{figure}
}

\subsubsection{\gdl's approach to reflective reasoning about logic systems}

One of \gdl's central insights was to recognize that
if we can model and embed ``a logic within a logic'' once
(\eg target logic instance $L_1$ into meta-logic instance $L_0$),
then we can probably do the same thing again,
or indeed any number of times.
In particular,
\gdl used as his meta-logic $L_0$
the traditional practice of working mathematics,
as expressed in a natural human language (German in this case)
augmented with traditional mathematical notation as needed.
Using this meta-language of conventional mathematics and German,
he carefully specified and constructed
a particular formal system to study in the target logic role of $L_1$.
Within $L_1$, however,
he further showed
(using the traditional mathematical tools available in $L_1$)
that both the syntactical and deductive rules of $L_1$
could be \emph{coded} and further embedded within $L_1$
to form an \emph{embedded logic instance},
which we will call $L_2$.
\Cref{fig:refl:goedel} illustrates
\gdl's three-level nesting of
(embedded) logic within (target) logic within (working) logic.

The syntax and rules defining $L_1$ might in principle
be either the same as or different than
the syntax and rules defining $L_2$.
\gdl's famous theorems, however,
focus on the particularly interesting case in which
logic instances $L_1$ and $L_2$ represent the \emph{same} formal system $F$
using identical syntax and rules
(but with those of $L_2$ necessarily coded within $L_1$
as part of getting the details right).
Because the working meta-logic $L_0$ is available
for reasoning about both $L_1$ and the $L_2$ embedded within it,
\gdl was able to \emph{prove}, within $L_0$,
this exact correspondence between $L_1$ and the embedded $L_2$.

\gdl's theorems do \emph{not}, however,
need $L_0$ and $L_1$ to be defined by the same formal system,
and in practice usually they are not.
Usually $L_0$ is only semi-formal ``pencil-and-paper'' mathematical practice
while only $L_1$ and $L_2$ are fully, precisely-specified formal systems.
Only in the exceptionally rare cases in which
\gdl's theorems have been re-proven in mechanically-verifiable fashion,
as in Paulson's formulation of these proofs in Isabelle~\cite{XXX},
can the working $L_0$ be said to be a truly \emph{formal} system at all.
Even in these rare cases,
$L_0$ has so far always been a \emph{different} --
typically richer and hence easier to use --
formal system than the target and embedded systems $L_1$ and $L_2$.

\baf{	XXX make more precise that the above figure
	represents \gdl's first incompleteness theorem;
	his second requires still another level.}

As mentioned above,
we can in principle carry on this coding and nesting of logic instances
as deeply as we might need to.

\gdl's \emph{second} incompleteness theorem,
in particular,
actually uses a second level of embedding within the target logic instance $L_1$,
as illustrated in \cref{fig:refl:goedel2}.
This is because a key initial step of \gdl's second incompleteness theorem,
as we will discuss further below,
is to carry out a proof of his own \emph{first} incompleteness theorem
within the target logic instance $L_1$.
That is, \gdl's second incompleteness theorem
``pushes down'' his first incompleteness theorem
from metalogic $L_0$ into target logic $L_1$,
then continues with further reasoning using $L_0$
about the implications of his first theorem having been proven in $L_1$.
Thus, while \gdl's informal presentation of his first theorem
in the language of working mathermatics and German
used only instances $L_0$ through $L_2$,
the ``pushed down'' version of his first theorem
carried out as part of his second theorem
instead uses instances $L_1$ through $L_3$.

\subsubsection{Using and labeling reflective instances of logic systems}

It should be clear that the numbers 0--3
we have assigned to these nested logic instances are arbitrary:
we could just as well have instead labeled them $L_1$ through $L_3$,
$L_{-1}$ through $L_2$, or $L_{10}$ through $L_{13}$.
It should also be obvious that this nesting of logics
could in principle be repeated any number of times.
Once we have done the hard work of formally coding $L_2$ within $L_1$,
in fact,
it becomes fairly trivial in principle
to construct a further $L_3$ within $L_2$ using the same rules:
\ie just invoke the same already-constructed coding method twice.
Thus, our ability to embed logics within logics to construct an instance $L_i$
certainly extends to any finite nesting level $i \ge 0$.

The same principle is true extending in the opposite direction as well,
as \gdl's second theorem illustrates.
If we have used logic instance $L_0$
draw some interesting conclusions about instances $L_1$ and $L_2$
(or about $L_1$ through $L_k$ for any nesting depth),
but we find ourselves ``running into a wall''
in terms of what we can prove within $L_0$,
one option available to us is to ``step back'' to a new vantage point
that we might now call $L_{-1}$:
\ie formalize our reasoning in $L_0$ sufficiently
so that we can now treat $L_0$ through $L_k$ as the targets of study
and use a \emph{new} meta- (or ``meta-meta-'')logic $L_{-1}$
as our new vantage point from which to perform this study.
This is, in fact,
just another perspective on what \gdl does in his second incompleteness theorem:
he takes his first incompleteness theorem (in logics $L_0$ through $L_2$),
``steps back'' to a new working logic $L_{-1}$,
formalizes the first theorem's reasoning within $L_0$
using the tools of $L_{-1}$,
and completes the second incompleteness theorem within $L_{-1}$.
It seems equally valid to say either that
\gdl's second theorem ``pushes down'' his first theorm
from $L_0$--$L_2$ to $L_1$--$L_3$,
or that he ``steps back'' by
leaving his first theorem at $L_0$--$L_2$
but starting to work in a new meta-logic at $L_{-1}$.

A bit more philosophically,
we might argue that the moment we even \emph{mention} explicitly
and start \emph{talking about} any particular logic instance $L_i$
as an actual \emph{subject} of study or conversation --
rather than merely as a body of background understanding
that we implicitly assume and hope
is sufficiently in common with the understanding of our audience --
we are implicitly ``stepping back''
into a new (meta-)logic instance $L_{i-1}$ in order to talk about $L_i$.
Taking this perspective,
from the moment we first mentioned $L_0$ at the start of this section,
we were \emph{already} implicitly ``stepping back'' 
into an unmentioned $L_{-1}$ in order to talk about $L_0$
(and subsequently about $L_1$ and $L_2$).
Thus, from this perspective
the entire earlier part of this section
was actually ``done'' in this implicit $L_{-1}$.
The moment we first mentioned $L_{-1}$ explicitly,
we implicitly assumed the existence of and started using
a logic instance $L_{-2}$.
By virtue of mentioning \emph{that} in the above sentence,
we implicitly started using a meta-logic instance $L_{-3}$, and so on.
By this perspective,
we always have ``need'' for, and implicitly assume the existence of,
at least one meta-level logic instance ``below''
the lowest one we have ever explicitly talked about.
Since there is obviously no limit in principle to the number of times
we might force ourselves to ``step back'' in this way,
it seems we may have potential need for the entire integer line,
positive and negative, in order to label our logic instances.

\subsection{Quoted GA terms via \gdl coding}
\label{sec:refl:coding}

\baf{	Include a precise formulation?}

Having illustrated the general idea of reflective reasoning,
we now summarize the coding techniques \gdl used
in just slightly more detail,
and how they apply in the context of \ga.

To set the scene,
we will be using the informal language of English
and mostly-traditional mathematical reasoning (except when noted otherwise)
as the meta-logic $L_0$ illustrated in \cref{fig:refl:goedel}.
The main target logic we will focus on in the $L_1$ role
will of course be \ga as defined in \cref{sec:ga},
again except when otherwise noted.
Our task is now to encode an embedded instance of \ga as $L_2$
within the $L_1$ instance,
using the conventional reasoning tools at our disposal in $L_0$ of course.

\subsubsection{Coding \ga term syntax}

Our first task is to encode \ga's syntax.
Since the $L_1$ instance of \ga ``knows'' only about natural numbers,
this means translating any \ga term $t$ into a natural number.
We accomplish this by defining a function within $L_0$
that transforms any syntactically-valid term $t$ into a natural number.
For convenience and consistent with tradition,
we will express this transformation as a special form of ``quoting'':
namely, given any \ga term $t$,
we will use the notation $\quo{t}$ to represent
term $t$ encoded into a natural number.
Thus, $t$ is a term but $\quo{t}$ is a numeric code for that term.

To be clear,
we feel free to call $L_1$ terms ``terms''
becuase we are reasoning about $L_1$ in our informal meta-logic $L_0$,
which we can safely assume has a rich type system 
in which we can inductively define our own types,
such as an $L_0$ type that correctly represents
(\emph{only}) syntactically-valid \ga terms.
$L_1$, however, does \emph{not} have a rich type system --
at least not natively --
but knows only about natural numbers,
so we must define a function in $L_0$
that transforms a \ga term in $L_1$'s syntax
into a natural number within $L_1$'s domain of discourse
(\ie a quantifiable object in $L_1$).

While in his proofs \gdl used a clever prime-number encoding of terms
that was well-matched to the arithmetical reasoning used in his proofs,
we will feel free to use (informally) more generic modern practices.
In particular, 
\cref{sec:ga:struct} already pointed out some standard ways
with which we may encode general finite structures,
such as lists or strings,
into (typically rather large) natural numbers.
Assuming we fix a concrete, string-based syntax for \ga terms,
using a character set whose characters likewise map to natural numbers,
we can encode a \ga term simply using the list encoding presented earlier.
Not all natural numbers will map \emph{back} to syntactically valid terms:
that is, this transformation will be injective but not bijective.
This is not a problem, however,
as in the relevant contexts
we will generally be interested only in natural numbers
that \emph{are} the codes of valid \ga terms.

All of the computations we need to perform in order to achieve this coding,
as well as related computations
such as validating a coded term
(checking whether an arbitrary natural number is actually the code for a term)
is all readily expressible via primitive-recursive reasoning,
requiring nothing either logically or computationally special or problematic.
As a result,
we can be certain of being able to define and reason about these functions
regardless of
whether we choose an adequately-powerful classical logic for our $L_0$,
or we want to ``inhabit'' and use some form of grounded deduction as our $L_0$
(\eg, a variant of \gd with a rich type system
allowing us to define and reason about \ga terms directly within $L_0$).
Either way, we can depend on basic,
primitive-recursive computation and reasoning working as expected.

\subsubsection{Coding proofs and provability}
\label{sec:refl:coding:proofs}

Having encoded \ga terms into natural numbers usable as objects in $L_1$,
our next step is to do the same with \ga inference rules and proofs.
We consider a \ga proof to be simply a list of strings,
each string having the form of
either a definition `$s(\vec{x}) \ldef \ttc{d}{\vec{x}}$'
or an entailment `$\Gamma \vdash p$'.
In the former case of a definition,
the defined symbol $s$ must not have been already defined earlier in the proof
and expansion $d$ must have no free variables
other than the formal arguments $\vec{x}$.
An entailment `$\Gamma \vdash p$' in the proof represents a deduction,
which must be justified by one of the inference rules of \ga.
A term $t_i$ is justified by an inference rule $R$
if the term matches the conclusion of $R$,
and each of $R$'s premises match a term appearing earlier in the proof
(before position $i$).

\baf{	deal somewhere with: 
	structural rules, background assumptions,
	perhaps explicit background definitions? }

Since a proof is simply a list of strings,
we can first encode each string comprising the proof into a natural number
as discussed in \cref{sec:ga:struct},
then use our standard list encoding again
to encode the proof (now as a list of natural numbers)
into a single natural number.
Since a proof has only a finite number of terms,
and the validity of each proof step may be verified in term
by checking it against a finite set of inference rules
(each of which may have a finite number of premises),
checking an encoded proof for validity is a readily decidable
and in fact primitive-recursive computation,
again non-problematic regardless of whether we might be operating
in a classical or grounded meta-logic $L_0$.
We will call this proof-checker function $V(P,t)$,
which returns \ctrue iff $P$ is a valid proof ending in term $t$.
\baf{	clarify: while defining \ga in $L_0$
	we first form a proof checker in $L_0$
	using $L_0$'s rich type system;
	but then we need to translate that proof checker
	into a boolean predicate definable in $L_1$
	taking the natural-number code of an encoded proof --
	but this is in principle also straightforward
	since $L_1$ also supports anything in PRA.}

While checking a purported proof $P$ for validity is primitive recursive,
\emph{searching} for a proof $P$ given only a term $t$ of interest
is of course another matter.
Finding a proof of a term $t$
in principle corresponds to an unbounded search through all possible proofs:
\eg an existential proposition of the form $\texists{P}{V(P,t)}$,
where $P$ ranges over (the codes of) all possible proofs.
Proof search is effectively computable, however:
\emph{if} a valid proof of $t$ exists,
then in principle an unbounded search through
all natural numbers that might encode valid proofs of $t$
will eventually find it.
(If a proof of $t$ does not exist,
then we may expect a computational search for such a proof
simply to run forever without terminating --
but we will return to this question later.)
Since \ga's recursive definition capability is Turing complete,
it is straightforward to express an $L_1$ function that recursively
searches through all natural numbers for one encoding a proof $P$
that the proof checker (in $L_1$) accepts
as proving a given coded term $\quo{t}$.
\baf{	explore somewhere: what exactly can we prove about this program
	within $L_1$?  and what can we prove only in $L_0$?
	or just ask this as a question for now?}

\baf{	alternative: expressing proof search
	as a direct existential term in $L_1$.
	Can return \cfalse if it finds a proof that a proof cannot exist!

	maybe define
	$\funp\quo{t} \equiv \neg\fprv\quo{t} \land \neg\fprv\quo{\neg t}$
	here?
}

\baf{	need the proof checker function to take premises, \ie support entailments. 
	$V(P,p,c)$ where $P$ is the proof,
	$p$ is a term representing premises assumed true (\ctrue if none),
	and $c$ is the conclusion to be proved.}

\baf{	figure out whether fprv should take terms or codes of terms,
	and be consistent about it.}

%later

\subsection{\gdl's first incompleteness theorem}

We may briefly and informally state \gdl's first incompleteness theorem
as follows:

\begin{quote}
For any formal system $S$ that includes arithmetic and classical logic,
if $S$ is consistent, then it is incomplete:
that is,
there is some proposition $p$
such that neither `$\vdash p$' nor `$\vdash \neg p$' is provable in $S$.
\end{quote}

\subsubsection{Historical context and shifting expectations}

\gdl's work came at a point in mathematical history infused with
the optimistic hope and expectation that
all mathematical problems should ultimately be solvable sooner or later,
the main question being merely how hard a given mathematical problem is
and how long it might take to resolve.
Emerging from this spirit of optimism,
it was widely hoped and even expected that we should be able to formulate
a foundation for mathematics embodied in a formal system
powerful enough to express, and ultimately be able to resolve,
any mathematical problem.
It would appear essential for such a system to be \emph{complete}:
that is, to yield a \ctrue or \cfalse answer
to any well-formed mathematical proposition,
even if that answer might be hard or take a long time to find.

\gdl's first theorem dashed these hopes by apparently proving that
any formal system $S$ surmounting the seemingly rather low bar
of being consistent
(\ie useful \emph{at all} to distinguish truth from falsehood),
and \emph{merely} powerful enough to reason about basic arithmetic
(\ie ignoring set theory and all the rest of mathematics),
cannot possibly be complete.
As a result,
it would seem that there cannot possibly be any formal system
usable --
even given unlimited time and resources --
to resolve \emph{all} mathematical questions,
which would obviously have to include questions about arithmetic.

We are designing \ga in a very different environment,
long post-\gdl,
in which the theory of computation has matured,
real computers are ubiquitous
along with sophisticated programming languages for them,
and we are thoroughly familiar with the theory and practical reality
of software bugs such as nonterminating programs.
In particular.
we do not hope or expect \ga to be complete:
in fact we hope and expect it \emph{not} to be complete,
at least by the definition \gdl used in his incompleteness theorems.

One of \ga's central design goals is to accept (via recursive definitions)
but safely ``reason around'' paradoxes such as the Liar `$L \ldef \neg L$'
without actually falling into them.
It would be unfortunate for \ga's chances of success
if we could find a way to prove either `$L$' or `$\neg L$'
given this definition,
since inconsistency would then follow immediately.
We have thus moved the target at the outset:
while the mathematicians of \gdl's time
were hoping for a complete (and consistent) system,
we are hoping for an incomplete (but consistent) system
that can gracefully reason around paradox.

Furthermore, in the above paragraph we just informally stated
a trivial proof of \gdl's first incompleteness theorem applied to \ga.
Like classical arithmetic,
if \ga is consistent then it \emph{must} be incomplete,
otherwise the readily-definable Liar sentence would make it inconsistent,
contradicting the assumption of consistency.
QED.
We do not need \gdl's sophisticated machinery
to prove this particular point about \ga.

Nevertheless,
it remains well worth exploring where \gdl's techniques and line of reasoning
\emph{does} lead in the context of \ga,
so in that spirit we wll press on.

\subsubsection{Indirect self-reference and self-replicating code}
\label{sec:refl:quine}

\gdl's proof skirts suspiciously close
to the type of self-reference we see in paradoxical statements,
in particular constructing a formula that talks about its own unprovability.
Because the classical systems of arithmetic \gdl's proof is about
do not allow direct self-reference
in the way that natural language \ga does, however,
\gdl's proof usees \emph{indirect self-reference}:
it talks about a natural number computed in a certain way,
which -- when computed in that carefully-prescribed fashion --
turns out to be the \gdl code or natural number that encodes the formula itself.

The form of indirect self-reference that \gdl uses is aptly illustrated
by this statement by W.V. Quine,
now known as Quine's Paradox:\footnote{
	See \cite{quine82mathematical}.
}

\begin{quote}
``yields a falsehood when appended to its own quotation'' \\
yields a falsehood when appended to its own quotation.
\end{quote}

This statement does not directly refer to itself, as the Liar paradox does;
it only makes a certain claim about the text in quotes.
The claim itself is expressed in the words strictly after the quotes,
and so is clearly, textually separate and disjoint
from the quoted text that the claim talks about.
We naturally interpret the phrase ``when appended to its own quotation''
as explicit instructions -- a \emph{program} encoded in English --
specifying a certain computation that we are asked perform on the text in quotes
(\ie append it to a quoted version of itself)
before we subject the resulting text
to the ``yields a falsehood'' predicate part of the claim.
But we find that when we follow these instructions --
we mentally ``run the program'' -- on the specific quoted text,
the result miraculously turns out identical to the entire statement above.
Thus, the statement indirectly refers to itself
by ``computing'' a full copy of its own text,
based on its quoted first half and the instructions in its second half.
This form of self-replicating statement or program
has come to be known as \emph{a Quine},
the creation of which -- in diverse languages and exotic flavors --
has become a popular pastime for overly-clever programmers.
\baf{some citations, e.g., obfuscated programming contests}

We can use this method of indirect self-reference
in \gdl's style of reflective reasoning
through the use of a \emph{diagonalization} or \emph{Quine function},
which implements an operation analogous to
appending a string to its own quotation.
We define a function $Q(t)$
taking as arguments the \gdl code of 
a term $\tto{t}{v}$ containing a particular free variable $v$
having special significance.
Our Quine function $Q$ computes and returns
the \gdl code for $t' = \tto{t}{\quo{t}}$,
\ie the original term $t$ except with all free occurrences of the variable $v$
replaced with occurrences of a natural number literal
representing the \gdl code of the input $t$ itself
($t$'s \gdl quotation).
In short, invoking $Q(\quo{\dots v \dots})$ returns
$\quo{\dots \quo{\dots v \dots}\dots}$.
$Q$ itself is just a primitive-recursive function
taking a natural number as its argument and returning a natural number,
so expressing and formally reasoning about it
is not a problem in either \pa or \ga.
\baf{	notation for conversion of a natural number
	to a natural number literal in a term?}

\baf{mention and cite Hofstadter as appropriate}

\subsubsection{Formulating the \gdl sentence}

Given such a Quine function $Q$,
one way we can formulate \gdl's indirectly self-referential formula $G$
is as follows:

\begin{align*}
G	&\ldef	Q(\quo{\neg \fprv(Q(v))}) \\
	&=	\quo{\neg \fprv(Q(\quo{\neg \fprv(Q(v))}))}
\end{align*}

The function $\fprv$ here is the provability predicate in the target logic
based on the proof checker function
described earlier in \cref{sec:refl:coding:proofs}:
\ie $\fprv(t) \ldef \texists{P}{V(P,t)}$.

The formula $G$ thus asserts the \emph{unprovability}
of the formula computed by the term $Q(\quo{\neg \fprv(Q(v))})$.
But since Quine function $Q$ applied to this particular code
yields the \gdl formula $G$ itself,
$G$ is a formula that indirectly asserts its own unprovability.

\subsubsection{Summary of the setup in meta-logic $L_0$ and target logic $L_1$}

Recall from \cref{sec:refl:vis} and \cref{fig:refl:logic}
that we are working in some mega-logic $L_0$
while reasoning about a target system $L_1$ that includes arithmetic.
To recap the setup we have outlined above,
we first formulated within our meta-logic $L_0$
the terms, formulas, and proofs constituting $L_1$,
along with a predicate in $L_0$ determining whether an $L_1$ proof is valid
given the axioms and inference rules of $L_1$.
This is just the task of constructing $L_1$ within $L_0$ in the first place.
We then defined \gdl coding functions in $L_0$
that convert $L_1$'s terms, formulas, and proofs into natural numbers,
the only first-class objects we wish to assume
that $L_1$ can operate on directly.
We next constructed predicates within the target logic $L_1$ --
\ie as formulas of $L_1$ --
that test whether a natural number
encodes a valid term, formula, or proof in $L_1$ itself.
We constructed a Quine function first within $L_0$,
then a corresponding one in $L_1$ behaving the same way.
Finally, we proved within $L_0$ the correspondence 
between $L_0$'s metalogical notion of valid $L_1$ terms/formulas/proofs
and $L_1$'s \gdl-coded notion of valid $L_1$ terms/formulas/proofs,
as well as the Quine function's correspondence between $L_0$ and $L_1$.
In our meta-logic $L_0$ we now definitively ``know'' (have proven in $L_0$)
that $L_1$ is powerful enough to reflect on and reason about itself,
and (using the Quine function)
that $L_1$ formulas in particular can indirectly refer to themselves.

To complete all this setup and obtain this reflection capability,
our requirements of the target logic $L_1$ were \emph{only}
that $L_1$ be able to express primitive-recursive functions on natural numbers:
that is, $L_1$ can be Skolem's primitive-recursive arithmetic (PRA)
or more powerful system that includes PRA.
Furthermore, so far \emph{in principle} we could even use PRA
as our meta-logic $L_0$,
directly representing terms, formulas, and proofs as natural numbers
``in our heads'' in the first place --
although in practice such a prospect would be unbearably tedious and error-prone
given the mental limitations of real humans like us.
The point is that so far 
we have not yet needed to take any reasoning steps in either $L_0$ or $L_1$
that push beyond the minimalistic foundation of PRA,
whose solidity seems beyond reasonable question
and certainly applies in \ga as well.

It is only in its key final steps that \gdl's first proof
needs to push beyond primitive recursive computation
and into reasoning whose applicability to \ga is more dubious.

\subsubsection{Consistency versus $\omega$-consistency}
\label{sec:refl:consist}

\gdl's original proof of his first theorem
requires assuming that the target logic $L_1$
is not just consistent but $\omega$-consistent,
a stronger property that implies consistency.
$L_1$ is $\omega$-consistent if there is no predicate $\tto{p}{x}$
such that $L_1$ proves $\vdash \neg \tto{p}{x}$
for each individual natural number $x$,
but $L_1$ also proves $\vdash \texists{x}{\tto{p}{x}}$.
That is, an $\omega$-inconsistent system
provably insists that there is some unspecified number $x$
satisfying $\tto{p}{x}$
while also provably denying $\tto{p}{x}$
whenever we replace $x$ with any particular, concrete natural number.

Notice that
in order to express the concept of $\omega$-consistency at all,
our target logic $L_1$ appears to need the unbounded existential quantifier.
Otherwise `$\texists{x}{\tto{p}{x}}$'
would not even be a well-formed formula in $L_1$,
so there is no way $L_1$
could ever possibly prove
either `$\vdash \texists{x}{\tto{p}{x}}$' or its negation.
Thus, by expressing and making an assumption of $\omega$-consistency,
\gdl has taken a first step beyond the boundaries
of primitive-recursive arithmetic.
Expressing $\omega$-consistency is not a problem for \ga, however,
since it has unbounded quantifiers just like classical predicate logic.

Rosser later strengthened 
\gdl's first theorem to assume only ordinary consistency,
at the cost of some additional complexity
in formulating the \gdl sentence and reasoning about it.
Expressing the ordinary consistency of $L_1$
does not require that $L_1$ have unbounded quantifiers.
Instead, we need unbounded quantifiers in our meta-logic $L_0$
to express the consistency of $L_1$ in the natural and standard fashion:
\emph{there exists} no formula $f$
for which `$\vdash f$' and `$\vdash \neg f$' is provable in $L_1$:

\begin{align*}
\fcon_{L_1}
&\ldef	\neg \texists{f}{
			\fprv_{L_1}(\quo{\vdash f}) \land
			\fprv_{L_1}(\quo{\vdash \neg f})
		} \\
&\ldef	\tforall{f}{
			\neg \fprv_{L_1}(\quo{\vdash f}) \lor
			\neg \fprv_{L_1}(\quo{\vdash \neg f})
		}
\end{align*}

Regardless of which version of consistency we assume,
\gdl's line of reasoning so far remains perfectly workable in \ga
despite stepping beyond primitive recursion.

%later

\subsubsection{The diagonalization lemma}

An important explicit step in a mature, modern formulation of \gdl's proof
has come to be known as
the \emph{fixed-point theorem} or \emph{diagonalization lemma}.
Given any arbitrary formula $\ttc{f}{v}$
containing exactly one free variable is $v$,
and provided $L_1$ is a classical system such as \pa,
we can prove that that there exists a closed formula $g$
for which the following biconditional is provable in $L_1$:

\[
	\vdash_{\pa} g \liff f(\quo{g})
\]

It is at this point that \gdl's line of reasoning
ceases translating so directly into \ga.
The proof of this biconditional
makes use of the Quine or diagonalization function
discussed above in \cref{sec:refl:quine},
and that part of the reasoning is not a problem in \ga
since the Quine function is primitive recursive.

However, proving the diagonalization lemma also makes use of
the classical introduction rule \irl{{\liff}I} for the biconditional.
In classical logic,
this rule's premises impose on us only two proof obligations:
first, showing a hypothetical chain of reasoning
leading from `$g$' to `$f(\quo{g})$'
(that is, proving `$g \vdash f(\quo{g})$');
and second, showing a hypothetical chain of reasoning
leading the other way from `$f(\quo{g})$' to `$g$'
(\ie proving `$f(\quo{g}) \vdash g$').
\ga's bidirectional introduction rule \irl{{\liff}I}
as presented in \cref{sec:prop:iff}, however,
imposes two additional proof obligations:
both sides of the proposed biconditional
must also first be \emph{proven boolean}.

In \ga as in classical logic,
a biconditional `$p \liff q$' is still logical equivalent
to the pair of implications `$p \limp q$' and `$q \limp p$',
so let us break the standard proofs of the diagonalization lemma
into these two parts as usual.

Starting in the reverse direction,
in classical logic we hypothetically assume `$f(\quo{g})$'
and must show a chain of reasoning leading in $L_1$
from this assumption to `$g$'.
In \ga, however,
this chain alone is insufficient to introduce the conditional
and prove `$\vdash f(\quo{g}) \limp g$':
to do so we also first prove `$f(\quo{g}) \jbool$'.
But the standard diagonalization lemma
gives us no information about the arbitrary formula $\ttc{f}{v}$
other than that it is well-formed and has exactly one free variable $v$.
In \ga, $\ttc{f}{v}$ might denote nothing;
we have no justification to presuppose that it yields a boolean value.

To prove this direction of the diagonalization lemma in \ga, therefore,
it appears we need to add some assumption about `$\ttc{f}{v} \jbool$'.
There happens to be only one value of $v$
that the diagonalization lemma actually needs $\ttc{f}{v}$ to be boolean for,
namely $\quo{g}$
(which takes a different form within the diagonalization proof
since the proof must construct this $g$).
As such, it appears we can at least prove in $L_0$
that for any formula $\ttc{f}{v}$ there exists a $g$
for which the following hypothetical is derivable, using \ga as $L_1$:

\[
	\ttc{f}{\quo{g}}\jbool \vdash_{\ga} f(\quo{g}) \limp g
\]

\com{	unnecessary detail and not really a problem for metalogical reasoning
Requiring that $\ttc{f}{v}$ be boolean
for \emph{all} natural numbers substituted for $v$
appears sufficient, for example,
leading to this entailment being derivable in \ga:

\[
	\tforall{v}{\ttc{f}{\quo{g}}\jbool} \vdash g \liff f(\quo{g})
\]

Since the diagonalization lemma does not actually need $\ttc{f}{v}$
to be boolean \emph{for all} $v$
but only for one particular $v$ of interest,
namely $\quo{g}$,
we might try to weaken the assumption
merely to `$\ttc{f}{v} \jbool$'.
But making use of the resulting ``stronger'' derivation with this assumption
will be difficult,
because $g$ is formula that the lemma asserts the existence of,
and hence brings into the picture only \emph{after}
successful application of the lemma.
}%com

Attempting the other direction,
we encounter a more fundamental problem.
We wish to prove the implication `$\vdash g \limp f(\quo{g})$'.
To do so, however,
we must first have not only a hypothetical chain of reasoning
leading from `$g$' to `$f(\quo{g})$' --
which we do --
but also a pre-existing proof in \ga that $g$ is boolean.
How might we obtain such a proof?
We know nothing unconditionally about $g$ as yet
because we are trying to construct it.
In \ga, to prove $g$ boolean,
it appears we must first prove $g$ boolean.

We observe that
the standard proof of the diagonalization lemma in the forward direction,
in fact,
has a structure nearly identical to the chain of reasoning
leading to Curry's paradox
if we admit the recursive definition `$C \ldef C \limp P$',
as we saw long ago in \cref{sec:prop:curry}.
The corresponding chain of reasoning that forms a key step in \gdl's proof
merely takes an extra detour along the way,
through some nontrivial natural number arithmetic and back.
	% XXX "nontrivial but computable"?
In classical logic,
the ``truth'' of $g$ in `$g \limp f(\quo{g})$'
springs into existence
via the diagonalization lemma's hypothetical derivation,
in exactly the same way as the truth of $C$ does in Curry's paradox.
\ga's additional prerequisite that $g$ first be proven boolean,
in contrast,
appears to block the diagonal lemma in this direction
in exactly the same way as it blocked Curry's paradox earlier.

\baf{	We could also just ass $g$ as an explicit assumption
	and see where that goes? }

If we wish to continue trying to push \gdl's reasoning through in \ga,
we therefore seem to be constrained to use
a drastically-weaker form of the diagonalization lemma
that works at all only in one direction,
and even then only under the newly-added assumption of `$f(\quo{g}) \jbool$'.

\baf{probably need to unpack this reasoning more at some point.}

\com{
introducing the biconditional required to prove the diagonalization lemma
about the given arbitrary formula $\ttc{f}{v}$ in
requires first that we \emph{know} (or have assumed)
something about $\ttc{f}{v}$:
in particular,
that $\ttc{f}{v}$ is boolean, at least for the values of $v$
that we wish to pass to it (namely $\quo{g}$ in this case).
Thus, it does not appear feasible to prove in \ga
the unconditional version of the diagonalization lemma above,
but instead only a weaker version
with a suitable booleanness assumption on $\ttc{f}{v}$,
such as the following:

\[
	f(\quo{g}) \jbool \vdash g \liff f(\quo{g})
\]

Despite this potential issue,
let us try to push forward nevertheless.
}

\subsubsection{Wrapping up \gdl's first incompleteness theorem}

Given all this setup,
the final steps of \gdl's first incompleteness theorem
become short and fairly simple.
In classical logic,
\gdl's proof first assumes that $L_1$ is $\omega$-consistent
as discussed above in \cref{sec:refl:consist},
uses uses the diagonalization lemma above
to form the \gdl sentence $G$ that asserts its own unprovability,
and finally uses a pair of classical proofs by contradiction in $L_0$
to show that neither `$\vdash G$' nor `$\vdash \neg G$'
is provable in $L_1$.

The fact that each of these critical last steps
involve proof by contradiction is already a warning sign
that we may be in trouble,
at least if we are using grounded deduction in our metalogic $L_0$.
In each case, \gd now requires us \emph{first}
to prove that there is a boolean fact-of-the-matter
about the provability of `$\vdash G$' and `$\vdash \neg G$', respectively,
in $L_1$.
This warning sign proves to be a false alarm, however,
at least if our meta-logic is the full non-constructive formulation of $\gd$
including the type-introduction rule \rlpo (see \cref{sec:quant:type}).
Since it is primitive-recursively decidable
whether a given natural number $P$ encodes a valid \ga proof,
we can prove `$V(P,f) \jbool$' for any given $P$,
by induction over $P$
using the primitive-recursive structure of the proof verifier $V$.
Inference rule \rlpo then allows us to infer `$\texists{P}{V(P,f)} \jbool$',
\ie that there is some boolean fact-of-the-matter
about whether such a proof $P$ exists or not.
This metalogical reasoning allows us to discharge
\gd's new prerequisites for proof by contradiction
and launch into each branch
of \gdl's hypothetical reasoning towards contradiction.

In any case, perhaps we prefer to use more-familiar classical logic
as our meta-logic $L_0$,
for purposes of \emph{studying} \ga as our $L_1$.
In this case,
we will obviously have no trouble launching \gdl's proofs by contradiction
since these happen in $L_0$.

We finally hit a roadblock within each of these branches, however.
It turns out that \gdl's proof in each case
relies on the diagonalization lemma in the forward direction,
`$g \limp f(\quo{g})$',
which we were unable to carry out in \ga
without satisfying an apparently-circular proof obligation.
In each case,
\gdl's proof uses the diagonalization lemma
to take us within $L_1$,
from a hypothetical assumption that $G$ (or its negation) is provable,
to an intermediate inference that `$\neg\fprv(\quo{G})$' is provable,
and from there on to contradict the $\omega$-consistency assumption.
Proving the unprovability of either `$G$' or `$\neg G$'
when our $L_1$ is \ga, therefore,
appears to be blocked by our emasculated diagonalization lemma,
at least along the standard line of reasoning.

\com{
The first point is where we try to apply the weakened diagonalization lemma above
to the \gdl sentence.
We can certainly envision many choices of the formula $\ttc{f}{v}$
for which we will be able to discharge
the weakened lemma's booleanness prerequisite `$f(\quo{g}) \jbool$'
with no problem.
The particular, interesting choice of $\ttc{f}{v}$
that makes \gdl's theorem work in classical logic, however,
is `$\neg\fprv(v)$',
for which we want the diagonalization lemma to deliver us
an unconditional proof `$\vdash G \liff \neg\fprv(\quo{G})$'.
If we attempt to achieve this goal
by discharging the diagonalization lemma's
new booleanness requirement, however --
namely `$\neg\fprv(\quo{G})$' in the case of interest --
then we encounter the problem of \emph{first} having to prove
}%com

\baf{	relate this discussion better to the visualization earlier:
	\ie how the $\fprv$ statements are venturing into $L_2$.
}

\subsection{\gdl's second incompleteness theorem}

We informally state \gdl's second incompleteness theorem as follows:

\begin{quote}
For any formal system $S$ that includes arithmetic and classical logic,
if $S$ is consistent,
then $S$ cannot prove itself consistent.
\end{quote}

\subsubsection{Historical perspective: \gdl's curse and the tree inspector}

This second theorem,
even though it was only sketched in \gdl's original paper
and rigorously proved only later,
nevertheless struck an even more devastating blow
to the apparent prospects of constructing a solid foundation for mathematics.
We would like to \emph{know} for certain,
and ideally be able to \emph{prove}
based on a ``minimal'' and ``self-evident'' set of principles,
that the foundation of our mathematical edifice is solid.

As a bare minimum,
we would at least like to know that our foundation is consistent,
\ie does not make \emph{all} formulas statements whether true or false.
Ideally a proof of consistency would be only a first step, in fact:
we would really like to know that the theorems of our mathematical foundation
are \emph{true}, not merely consistent.
(It is easy to invent any number of beautiful theories of flying pigs,
which may be perfectly consistent,
as long as we are careful never to include any axioms or inference rules
implying that pigs don't fly.)
Our standard interpretation of \gdl's second incompleteness theorem,
unfortnately,
appears to dash our hopes
even of achieving that bare-minimum first step
of proving the \emph{consistency} of any realistic foundation for mathematics.

The situation that \gdl's second incompleteness theorem leads us to
is that a classical formal system powerful enough to include arithmetic
can be proven consistent only in a \emph{strictly more powerful} system.
Thus, to prove the consistency of any classical system $S$ of interest,
we appear to need a different system $S'$
that starts with stronger foundational assumptions than $S$.
At least when it comes to arithmetic,
it appears there can be no ``bootstrapping''
from simpler to more sophisticated foundational assumptions:
it goes only the other direction.

As an illustration of the predicament that \gdl's results place us in,
let us imagine a logician to be analogous to an official
assigned to inspect the health and safety of trees in an urban park,
and to ensure that any unsafe trees are pruned or cut down.
Our safety inspector finds himself examining a tree
with a peculiar property:
the foliage growing \emph{downward} from each branch is so dense
that it is impossible to inspect the branch's safety by looking at it from below.
To check the safety of a given branch,
our inspector must actually climb the tree and view the branch from above.
Further,
the potential safety risks in this tree are such that
merely climbing up the trunk is not sufficient:
the branches are long and gnarled enough
that the inspector must actually climb out onto each branch,
potentially to its very tip,
in order to look back on and check the interior portions of that branch.
If our safety instructor makes an incorrect guess
about how far is safe to climb out on a particular branch,
our unfortunate inspector risks breaking the branch and falling to his death.
The inspector discovers that it is possible
to inspect a lower branch \emph{completely}
if he climbs up to a higher branch and looks down on the lower branch --
but this does not help him much in managing his risks,
since the higher branches are invariably smaller and more fragile
than the lower ones --
so his risk of death only ever increases by climbing higher.
There simply seems to be no safe way to inspect the tree for safety.
Was this tree created as a cruel joke by a mischevious deity?

The conundrum that \gdl's theorems place the logician in
seem arguably even worse than our tree inspector,
in that the the tree's branches (representing formal systems)
would be infinitely long:
most interesting formal systems have an infinite number of theorems.
The logician's tree also has an infinite number of branches,
since there are clearly an infinite number of possible formal systems
exhibiting no apparent limit to their sophistication, complexity,
or fragility to breakage (\ie inconsistency).

\subsubsection{Applying the second incompleteness theorem to \ga}

Most statements of \gdl's incompleteness theorems
mention that the system $S$ in question must include arithmetic,
but they often neglect to mention the \emph{other} crucial assumption
that $S$ includes classical logic.
This habitual lapse in precision about stating the assumptions
may be attributable at least in part to the fact that
intuitionistic logic --
the only real ``competitor'' to classical logic
that gained even modest traction around \gdl's time --
was found to be equiconsistent to classical logic.
For most formulations of intuitionistic logic,
if `$\vdash p$' is classically provable,
then the double-negation of the same statement,
`$\vdash \neg\neg p$', is intuitionistically provable.
If a classical formal system is inconsistent, therefore,
then trivially so is the corresponding intuitionistic system.
Switching to intuitionistic logic is thus, obviously,
no help for the logician to escape ``\gdl's curse'' as we might call it.
\gdl's theorems would thus seem equally applicable
to essentially all of the powerful formal systems
that have ever obtained any significant traction or use
in regular working mathematics.

Let us see now how \gdl's second incompleteness theorem applies to \ga.
\gdl's second incompleteness theorem builds on his first incompleteness theorem.
The second theorem uses the first theorem twice, in fact:
first, reasoning in the meta-logic $L_0$ about the target logic $L_1$,
and second,
reasoning within the target logic $L_1$
about an embedded instance $L_2$ of the same logic within $L_1$.
That is, to prove \gdl's second theorem,
we must in essence ``replay'' the first theorem twice,
once in $L_0$ and then again in $L_1$.
The considerable technical tedium of this process
may be mitigated using a shortcut known as the
Hilbert-Bernays derivability conditions~\cite{XXX} --
but since we will not be detailing the entire second proof,
we will not need this shortcut for the moment.

The fact that our attempt to apply \gdl's proof of the first theorem in \ga
came up short might tempt us to give up already.
Recall, however,
that we did not actually need \gdl's sophisticated \emph{proof}
to arrive at the conclusion that his first incompleteness \emph{theorem}
is nevertheless true of \ga.
Because we can directly express the Liar paradox `$L \ldef \neg L$' in \ga,
we get a trivial proof that if either `$\vdash L$' or `$\vdash \neg L$'
is provable in \ga,
then \ga is inconsistent because the other is provable as well.
Let us therefore try to proceed with
translating \gdl's second incompleteness theorem into \ga,
but using this trivial proof of \ga's incompleteness
in place of \gdl's involved first proof.

Let $G$ be the \ga formula `$L$' after defining `$L \ldef \neg L$'.
Reasoning in a grounded meta-logic $L_0$ about as our $L_1$ target logic,
we deduce from the non-constructive type-introduction rule \rlpo
(\cref{sec:quant:type})
that there is a boolean fact-of-the-matter
about whether entailment `$\Gamma \vdash p$' is provable in $L_1$,
including one representing our alternative \gdl sentence `$\vdash G$'.
It is also primitive-recursively decideable
whether or not both `$\vdash p$' and `$\vdash \neg p$' are provable
for any given formula $p$ --
\ie whether formula $p$ in particular causes inconsistency in $L_1$.
The same non-constructive rule \rlpo
therefore allows us to deduce in $L_0$ that
the question of whether $L_1$ is consistent --
\ie whether \emph{there exists} such a $P$ causing inconsistency --
also has a boolean fact-of-the-matter.

Notice that all of these strings
we believe to be definitely either provable or unprovable,
however,
are \emph{entailments} in \ga.
Recall further that \ga does not permit the same freedom of movement
between entailment and implication as in classical logic.
In particular,
the entailment `$p \vdash q$' might be provable in $L_1$
while the implication `$\vdash p \limp q$' might be unprovable,
due to the booleannness test that the \irl{{\limp}I} rule imposes.
More generally,
our belief that all of these entailments are either provable or not,
reasoning in \gd as $L_0$,
does not translate into an expectation that every \ga formula $p$
is either true or false \emph{within} $L_1$.

Just as \gdl's second incompleteness theorem does,
we are confident of being able to ``push'' all of the above reasoning down
from our meta-logic \gd as $L_0$ into our target logic \ga as $L_1$.
Defining the consistency of $L_1$ as above in \cref{sec:refl:consist},
we can prove in $L_1$ all of the following:

\begin{align*}
	& \vdash \tforall{e}{\fprv(e) \jbool} \\
	& \vdash \fcon_{L_1} \jbool \\
	& \vdash \fprv(\quo{\vdash G}) \jbool \\
	& \vdash \fprv(\quo{\vdash \neg G}) \jbool \\
%G	& \vdash \neg G		\\
%\neg G	& \vdash G		\\
	& \vdash \fprv(\quo{G \vdash \neg G})		\\
	& \vdash \fprv(\quo{\neg G \vdash G})		\\
	& \vdash \fprv(\quo{\vdash G}) \limp \fprv(\quo{\vdash \neg G}) \\
	& \vdash \fprv(\quo{\vdash \neg G}) \limp \fprv(\quo{\vdash G}) \\
%\fprv(\quo{\vdash G}) \lor \fprv(\quo{\vdash \neg G})
%	& \vdash \neg \fcon_{L_1} \\
	& \vdash \fcon_{L_1} \limp	\neg\fprv(\quo{\vdash G}) \\
	& \vdash \fcon_{L_1} \limp	\neg\fprv(\quo{\vdash \neg G}) \\
	& \vdash \fcon_{L_1} \limp	\neg\fprv(\quo{\vdash G}) \land
					\neg\fprv(\quo{\vdash \neg G})
\end{align*}

The last theorem in $L_1$ states \gdl's first incompleteness theorem,
proven via our alternate \gdl sentence $G$ above based on the Liar paradox.
The booleanness of the provability predicates
allow us to invoke proof by contradiction
to prove `$\vdash G$' and `$\vdash \neg G$' unprovable
while assuming $\fcon_{L_1}$,
and the booleanness of $\fcon_{L_1}$
allow us to invoke \irl{{\limp}I} to get the final statement
of the first incompleteness theorem.

Now we try to finish \gdl's second incompleteness theorem in \gd.
Although the last theorem above states \gdl's first theorem in full,
we got this along the way:

\begin{align*}
	& \vdash \fcon_{L_1} \limp	\neg\fprv(\quo{\vdash G})
\end{align*}

If we were working in a classical $L_1$
and could have followed the first incompleteness theorem ``correctly''
along \gdl's line of reasoning,
then the diagonalization lemma would have established
the biconditional `$\vdash G \liff \neg\fprv(\quo{\vdash G})$',
making the last statement above equivalent to `$\vdash \fcon_{L_1} \limp G$'.
But then if there were a proof of `$\vdash \fcon_{L_1}$' in $L_1$,
then that would make $G$ provable in $L_1$,
thereby contradicting $G$'s own statement that it is unprovable
and thus rendering $L_1$ inconsistent.
But with \ga as our $L_1$,
we were neither able to establish the diagonalization biconditional
at full strength in either direction,
nor to use the correct \gdl sentence in the first proof,
so our attempt at the second proof runs aground in \ga as well.

Examining more modern formulations of \gdl's second incompleteness theorem
that use the Hilbert-Bernays derivability conditions
and L\"ob's theorem \cite{mendelson15mathematical},
we encounter the same problem.
L\"ob's theorem also relies on the problematic forward direction 
of the biconditional that the diagonalization lemma is supposed to establish,
but appears unable to in \ga because if \ga's booleanness proof requirements.

\subsubsection{The aftermath}

None of this exploration constitutes
conclusive evidence that \gdl's second incompleteness theorem
is inapplicable to \ga, of course.
There might be a way to repair \gdl's proof and apply it to \ga.
For now let us take the inconclusive evidence we have at face value, however.
\ga's addition of typing requirements to key inference rules
\emph{appear to} prevent it from ``falling into''
either the Liar or Curry's paradox,
whether they are expressed directly via unrestricted recursive definitions
or expressed indirectly through \gdl's arithmetical wizardry.
If true, is this property of \ga desirable or undesirable?
As software developers often ask, ``is that a bug or a feature''?

\baf{	Maudlin...
	Do we actually \emph{want} \gdl's theorems --
	the second in particular,
	to be provable in our ``ideal'' formal system?

	Let us clearly separate the two main contributions
	of \gdl's incompleteness results:
	first, the techniques of coding and reflection
	that he introduced along the way to his theorems;
	and second, the two incompleteness theorems themselves.
	His coding and reflection techniques unquestionably
	have tremendous  value and utility to
	logic and computer science especially, and mathematics generally.
	The two incompleteness theorems themselves, however,
	seem to have rather little utility \emph{as tools}
	other than to prove other impossibility results.
	That is, the main outcome of the incompleteness theorems
	is to establish two major roadblocks in mathematical reasoning,
	whose main practical utility seems to be to erect other roadblocks.
	It is hard to find a single, \emph{positive} use in mathematics
	for either of the two theorems themselves:
	\ie an important theorem about possibility rather than impossibility
	that would not be readily provable
	without the incompleteness theorems.

	If this observation is even modestly on-target,
	then should we be particularly sad if we were to find a logic
	that can prove many interesting theorems
	but can't prove \gdl's incompleteness theorems?
	What if this alternative logic cannot prove these theorems
	because they \emph{aren't true} in that logic
	and because the roadblocks they represent
	\emph{don't exist} in the alternate universe of that logic?
	What if we could live in a logic without so many roadblocks to reasoning?
	What if the \emph{unprovability} of the incompleteness theorems
	in that logic
	is because interesting new avenues of reasoning
	are \emph{possible} in that logic
	that would not be possible in classical logic?
	(what if...dragon analogy?)

	Were \gdl's theorems actually serving more as
	roadblocks in our reasoning than tools of empowerment?
	The 
	Should we perhaps be interpreting many of the ``impossibility proofs''
	derived from \gdl's theorems and similar diagonalization arguments
	not as set-in-stone monoliths of truth,
	but rather as warning signs
	that perhaps we have just been working in the \emph{wrong} logic?

	Was \gdl's dashing of our hopes
	for a ``complete'' formal system
	actually just a canary in a coal mine, 
	signaling not just that \emph{that} particular hope was na\"ive,
	but in fact that our basic rules of reasoning might be flawed?
}
\baf{	Can we establish a line of reasoning 
	we might sloganize as ``\gdl for goopd'' --
	using \gdl's reflection techniques,
	in a sufficientaly non-classical logic,
	to enable interesting new lines of reasoning
	rather than to throw up mathematical roadblocks
	as his theorems appear to do in classical logic?
}
\baf{	swords into plowshares,
	weapons of math destruction into tools of empowerment?
}

\baf{	Observation:
	our incorrect interpretation of \gdl's theorems,
	at least as applied to classical logic,
	has been like a dragon jealously guarding
	a hoard of gold it has no use for --
	namely the possible cornucopia of \emph{other}
	interesting and perhaps more constructive uses
	of \gdl's coding and reflection techniques.
}

%later

\later{Related to Berry paradox...
"A New Proof of the Godel Incompleteness Theorem"
%\url{http://topologicalmedialab.net/xinwei/classes/readings/Gödel/Boolos%20proof%20Godel%20Incompleteness%201989.pdf}

Some objections to Boolos's proof (not necessarily competent/valid):
%\url{https://www.jamesrmeyer.com/pdfs/ff_boolos.pdf}
%\url{https://math.stackexchange.com/questions/3868997/boolos-s-proof-of-gödel-s-first-incompleteness-theorem-what-am-i-getting-wrong}
%\url{https://math.stackexchange.com/questions/3885919/boolos-s-proof-of-the-first-incompleteness-theorem-predicate-cx-y-and-assum}

"On proofs of the incompleteness theorems based
on Berry’s paradox by Vopˇenka, Chaitin, and Boolos"
%\url{https://www.researchgate.net/profile/Taishi-Kurahashi/publication/264380943_On_proofs_of_the_incompleteness_theorems_based_on_Berry%27s_paradox_by_Vopenka_Chaitin_and_Boolos/links/5a390ff5aca27266cfd34247/On-proofs-of-the-incompleteness-theorems-based-on-Berrys-paradox-by-Vopenka-Chaitin-and-Boolos.pdf}
}%later

\later{Philosophical (mis)interpretation of \gdl's theorems:
e.g., Nagel/Newman, "\gdl's theorem", Scientific American 1956,
%\url{https://www.jstor.org/stable/pdf/24943884.pdf?refreqid=fastly-default%3Ad92934ac2ffe3f24eba9852e0aaf8cef&ab_segments=&origin=&initiator=&acceptTC=1}:

"\gdl's conclusions have a bearing on
the question whether a calculating machine
can be constructed that would
equal the human brain in mathematical
reasoning. Present calculating machines
have a fixed set of directives built into
them, and they operate in a step-by-step
manner. But in the light of Codel's incompleteness
theorem, there is an endless
set of problems in elementary number
theory for which such machines are
inherently incapable of supplying answers,
however complex their built-in
mechanisms may be and however rapid
their operations. The human brain may,
to be sure, have built-in limitations of its
own, and there may be mathematical
problems which it is incapable of solving.
But even so, the human brain appears
to embody a structure of rules of
operation which is far more powerful
than the structure of currently conceived
artificial machines. There is no immediate
prospect of replacing the human
mind by robots."
}%later

\section{\ga and primitive recursive arithmetic (PRA)}
\label{sec:pra}

\later{Potential title: Cantor's Paradise Lost}

In the formulation of \ga so far,
everything expressible is computable and in a sense constructive,
even if \ga does not precisely follow
the intuitionistic path of constructivism.
An interesting further question to ask is whether and how
we might consistently ``strengthen'' \ga further
to allow something more like classical reasoning in \ga{} --
to express and reason about
non-computable functions and real numbers, for example.

For example, while we obviously
cannot consistently adopt the full law of excluded middle (LEM) in \ga,
a relevant question is whether we could adopt
some other weaker but still-useful non-constructive principle,
such as the \emph{limited principle of omniscience} (LPO):
namely that given an infinite series of natural numbers $n_i$,
either every $n_i$ is zero or there exists an $i$ such that $n_i$ is nonzero.
Earlier drafts of this document included the LPO
as an optional non-constructive extension to constructive \ga (\cga),
leaving open the important question
as to whether \ga with the LPO is still consistent.
The line of reasoning below, however,
makes me pessimistic as to whether at least this particular extension --
along with numerous other attractive extensions we might envision --
could be added to \ga while preserving consistency.

The rest of this section is applicable to but is not specific to \ga.

\subsection{Primitive-recursive arithmetic ($L$) and an extension ($L'$)}

Let $L$ be the language and logic of
primitive-recursive arithmetic (PRA).
Assume that we have proven $L$ to be consistent
using a more-powerful classical metalogic sufficient for this purpose,
such as set theory or Peano arithmetic.
Assume also that we have proven \gdl's incompleteness theorems in $L$,
which means that we have \gdl codes for
the primitive-recursive functions expressible in $L$,
for the logical formulae of $L$,
and for proofs in $L$.
We will not actually need \gdl's theorems below,
only the coding and reflecction ``tooling''
that he developed along his way to them.

We now create a new logic $L'$
that slightly extends $L$
by adding a single additional 2-argument function $\mu(f,x)$.
This function takes as its first argument $f$ the \gdl code in $L$
of a 2-argument primitive-recursive function $f(x,y)$ expressible in $L$.
For any natural number $x$,
if there exists a natural number $y$
such that $f(x,y)$ returns a \emph{nonzero} result,
then let $y_\mu$ be the least such $y$,
and in this case $\mu(f,x)$ returns $1+y_\mu$.
If no such natural number $y$ exists for the given value of $x$,
then $\mu(f,x)$ returns 0.\footnote{
	This function $\mu$ is just a slight variation
	on Kleene's minimization operator~\cite{kleene52introduction},
	which searches for a value $y$ for which $f$ returns zero
	instead of nonzero.
}

We can express any computation
as a primitive-recursive function $f(x,y)$
for which $x$ serves as the computation's input
and $y$ serves as a \emph{step count}
indicating a specific number of steps to run the computation.\footnote{
	Again see \cite{kleene52introduction}, for example.
}
If the computation terminates within the first $y$ execution steps,
then $f(x,y)$ returns $1+r$,
where $r$ is the computation's output.
If the computation has not yet terminated within $y$ steps,
then $f(x,y)$ returns 0 to indicate as such.
In this latter case,
the computation may or may not eventually terminate:
\ie $f(x,y)$ might or might not returns a nonzero result given some larger $y$.

While the function $f$ is primitive-recursive
and merely verifies whether the computation it represents
has terminated within $y$ steps (and if so with what result),
asking whether a value $y$ \emph{exists} for which the computation halts
expresses an unbounded search problem.
The $\mu$ function we incorporated into $L'$
presumes to answer this unbounded search question,
by finding the least step count $y$ for which the computation terminates,
if it ever does so,
and by detecting and returning 0 if the computation never terminates.

We can see that $\mu$ is clearly a non-computable function:
it solves the halting problem,
and as such constitutes a \emph{halting oracle}.
Nevertheless,
from the viewpoint of any sufficiently-powerful classical meta-logic
in which we might model and analyze our extended logic $L'$ --
\eg from set theory or Peano arithmetic --
$\mu$ will be a readily-definable function in our metalogic.
By nearly the same reasoning with which we proved $L$ consistent,
we expect to be able to prove $L'$ consistent in our classical metalogic,
by virtue of the fact that
every statement $L'$ makes is true,
provided we interpret the truth of these statements
appropriately in our metalogic.

Since \gdl's incompleteness theorems about $L$ were already provable in $L$,
they are certainly provable about $L'$ in $L'$, as well,
using nearly-identical reasoning,
since $L'$ merely extends $L$.
Our \gdl coding of the available functions
will need slight adjustment to get from $L$ to $L'$, of course,
to incorporate the added $\mu$ function.
Our proofs of \gdl's theorems in $L'$
will not need to \emph{invoke} this function, however,
so the development of the theorems are otherwise largely unaffected
by the addition of $\mu$.
\com{
Because \gdl's proofs are purely syntactic and constructive, in fact,
we can prove \gdl's theorems \emph{about} the stronger system $L'$
\emph{in} the weaker system $L$.
}

Since the syntax and inference rules of $L'$
are primitive-recursively definable,
we can define primitive-recursive functions expressible in $L$ (PRA)
that effectively reason about $L'$,
at least to the extent of confirming specific theorems of $L'$.
We do not expect $L$ to prove $L'$ consistent, of course,
which \gdl's second incompleteness theorem assures us will be impossible.
However, within the primitive-recursive $L$
we can still model the \gdl coding and other functions
defining the deduction system of $L'$.
That is, we can prove within $L$ that such-and-such is a theorem of $L'$.
Thus, $L$ can characterize and effectively enumerate the theorems of $L'$,
even though we might say that $L$ probably does not ``believe'' these theorems.
That is,
$L$ cannot be expected to ``know'' or
prove that theorems of $L'$ are ``true''
since, by Tarski's undefinability theorem,
$L$ is too weak even to express, let alone prove, a notion of truth for $L'$.
Beyond expressing in $L$ the existence of particular theorems of $L'$,
we can take our proof of \gdl's second incompleteness theorem for $L'$
and push it into our embedded instance of $L$.
So $L$ proves that \gdl's second incompleteness theorem is \emph{provable}
in $L$'s model of $L'$.

\subsection{Computable real numbers in $L'$}

Now we will construct some semblance of the real numbers
in our logic $L'$ extended with the halting oracle $\mu$.\footnote{
	In early feedback on this draft,
	Samuel Gruetter pointed out that this line of reasoning
	can probably be simplified by focusing
	only on functions from natural numbers to booleans --
	in effect functions defining \emph{sets} of natural numbers --
	thus avoiding the need for Dedekind cuts or digit calculations.
	This simplification appears to be
	a valid and attractive alternative approach.
	The main text here maintains for now
	the development using real numbers, however,
	both because that is the way the author first conceived this argument,
	and also for alignment with Cantor's famous theorem,
	which centrally inspired this line of reasoning.
}
For this purpose we will use the method of Dedekind cuts\footnote{
	See~\cite{dedekind63essays}.
}
to express a real number as a division or \emph{cut} of the rational numbers
into a lower part and an upper part, such that:
\begin{itemize}
\item	The lower and upper parts are nonempty:
	each contains some rational number.
\item	The lower part is downward-closed:
	if $x$ and $y$ are rationals, $x \le y$, and $y$ is in the lower part,
	then so is $x$.\footnote{
	It follows that the upper part is correspondingly upward-closed.
}
\item	The lower part contains no greatest rational number:
	for any rational $x$ in the lower part,
	there is some rational $y > x$ that is also in the lower part.\footnote{
	The upper part contains a least rational number
	if the Dedekind cut represents exactly that rational number.
	Otherwise the upper part contains no least rational number,
	which means that the Dedekind cut represents an irrational number.
}
\end{itemize}

We will define a real number in $L'$
to be a pair $(f,p)$ consisting of two components:
a \gdl-coded primitive-recursive function $f$,
and a \gdl-coded proof $p$ of a specific form, to be defined shortly,
within our extended logic $L'$.\footnote{
	The representation of $f$ could in principle use
	the \gdl coding from either $L$ or $L'$,
	since only primitive-recursive functions are allowed for $f$.
	We will assume the use of $L$'s \gdl coding.
	The representation of $p$ must use the \gdl coding of $L'$,
	however,
	since the proof $p$ needs access to the $\mu$ function extension.
}

The primitive-recursive function $f(r,s)$
takes as its arguments a \gdl-coded rational number $r$ and a step count $s$,
runs some computation defined by $f$ for $s$ steps,
and returns:
$0$ if the computation has not terminated yet within $s$ steps;
$1$ if the computation has terminated and found $r$ to be in the lower part;
and
$2$ if the computation has terminated and found $r$ to be in the upper part.
Since $f$ represents a computable function,
we do not expect that $f$ can recursively enumerate in this fashion ``all''
of the real numbers that set theory assures us exist --
only some countable subset of them --
but this subset of computable reals will suit our purposes for now.

Along with the \gdl-coded function $f$, however,
we require that any real number $(f,p)$ also contain
an accompanying \gdl-coded \emph{proof} $p$ in $L'$
that the standard properties of Dedekind cuts summarized above
hold with respect to the corresponding \gdl-coded function $f$.
That is, $p$ proves in $L'$
that $f$ expresses a computation that always terminates
when given a \gdl-coded rational number as input,
and that \emph{when} $f$ terminates,
its output assigns the rationals to the upper and lower parts
consistently with the three key properties of Dedekind cuts
summarized above.\footnote{
	Our real numbers encoded in $L'$
	in essence constitute \emph{proof-carrying code}~\cite{necula97proof}.
}

While $f$ itself is primitive recursive
and cannot invoke our halting oracle $\mu$,
the associated proof $p$ in $L'$
is free to use the $\mu$ oracle in \emph{reasoning about} $f$.
This capability makes up for a key limitation of $L'$:
namely that it is quantifier-free, like $L$ (PRA),
and can hence express only the top-level quantification
implicit in formulas with free variables.
Thus, we cannot directly express in $L'$
the existentially-quantified predicate-logic statement
``there exists a step count $s$ for which $f(r,s)$ terminates'' --
but within $L'$ we can invoke $\mu(f,x)$
and reason classically about its result.
In particular, if we find a particular step count $s$
within which $f$ demonstrably terminates --
because we can actually execute $f$ to termination within $s$ steps --
then we can prove in $L'$ that $\mu(f,x)$ is nonzero and at most $s+1$.
Further, since $L'$ is based on (quantifier-free) classical logic,
we can use all the tools that come with it,
such as proof by contradiction.

Similarly, the third key property of Dedekind cuts --
that the lower part contains no greatest rational number --
is most conventionally and easily expressed as an existential property:
``for every rational $x$ in the lower part
there exists a greater rational $y$ also in the lower part''.
To satisfy this apparent existence-proof requirement
in the quantifier-free $L'$, however,
we can require the proof $p$ to specify a \emph{computation}
that takes any rational $x$ in the lower part
and \emph{computes} a greater rational $y$ still in the lower part.
This proof in $L'$
can use the $\mu$ oracle as needed to express and satisfy the requirement
that this computation terminates with a suitable result.

\subsubsection{Computing digits of real numbers in $L'$}

Given any real number $r = (f,p)$ in $L'$
as defined above,
we can express within $L'$
the computation of any particular digit
in the binary representation of $r$,
and we can similarly prove in $L'$ that this computation terminates.
For simplicity, we focus attention on real numbers strictly between $0$ and $1$,
since the rest are just negations and/or inverses of these reals,
or else trivial special cases like 1.
Assuming $0 < r < 1$, therefore,
let $r_i$ be binary digit $i$ in $r$'s fractional binary representation.

Suppose we have already computed the first $i$ fractional digits of $r$,
meaning we know a natural number $n$
such that $n/2^{i} \le r < (n+1)/2^{i}$.
To compute digit $i+1$,
we invoke 
the computation represented by primitive-recursive function $f$
with the rational number $(n+1)/2^{i+1}$ as input.
This computation effectively tests
whether this rational is in $r$'s upper part, hence $r_i = 0$,
or in $r$'s lower part, hence $r_i = 1$.

If we believe the prerequisite proof $p$ in $L'$
that $f$ always terminates with an answer suitable for a Dedekind cut,
then we must similarly believe that this digit computation always succeeds.
Further, given any particular correct proof $p$ in $L'$
of the correct behavior of cut-defining function $f$,
we can construct a proof in $L'$
that this digit computation indeed terminates for any given digit $i$.
We can perform this $L'$ proof construction process (building on $p$)
not only in $L'$,
but even in $L$,
since the essentially cut-and-paste combination of \gdl-coded proof $p$
with the extension of that proof to digit computation in $L'$
is certainly a primitive-recursive proof-composition task.
We do not expect either $L$ or $L'$ to be able to prove
that these proofs in $L'$ are ``true'' --
again due to Tarski's undefinability theorem --
but from a classical metalogic more powerful than $L'$
we can convince ourselves that these proofs in $L'$ express true theorems.

\subsubsection{Comparing real numbers in $L'$}

Given any two real numbers $r_1$ and $r_2$ as defined above,
we can readily formulate a computation
to compare $r_1$ with $r_2$,
essentially by computing successive digits of each real number
and performing an unbounded search for any difference.
That is,
we can form a primitive-recursive function $c(r_1,r_2,s)$
that takes as arguments
two real numbers $r_1 = (f_1,p_1)$ and $r_2 = (f_2,p_2)$
and a step count $s$,
and calculates successive digits of $r_1$ and $r_2$,
terminating at some step count $s$ if there is any difference
to be found in $s$ steps.
If $r_1$ and $r_2$ represent the same real number, however,
$c(r_1,r_2,s)$ will never terminate for any $s$.
Thus, while $c$ represents a computable function,
we do not claim it to be a particularly \emph{useful} computation
in practice -- only in theory.

Nevertheless, within $L'$,
we can use the $\mu$ oracle to express and reason about
$c$'s termination or lack thereof.
Assume we formulate $c(r_1,r_2,s)$ as $c((r_1,r_2),s)$:
that is, as a 2-argument primitive-recursive function
that takes $(r_1,r_2)$ together as a \gdl-encoded pair,
using Cantor's pairing function for example.
Then within $L'$ 
we can use $\mu(c,(r_1,r_2))$ as a real-number equality test,
returning the step count at which $c$ finds the first digit difference
if $r_1$ and $r_2$ are unequal,
and returning 0 if $r_1 = r_2$ and hence $c$ never terminates.

\later{do we need a predicate to find the first differing digit?}

\subsubsection{Enumerating the real numbers}

Since the syntax and proof-checking predicates
for the deduction system $L'$ (like $L$)
are primitive-recursive,
we can express and reason in $L$ about primitive-recursive functions
to verify purported proofs of theorems in $L'$.
Further, via primitive-recursive step functions,
we can express in $L$ computations
that effectively perform unbounded searches --
not only for proofs in $L$,
but for objects \emph{containing} such proofs,
such as the real numbers as we defined them above.

In particular,
since we have a primitive-recursive test
for whether a pair $(f,p)$ represents a real number as defined above,
we can form a primitive-recursive step function $N(n,s)$
taking a natural number $n$ and step count $s$,
which searches starting from $n$
for the ``next'' real number whose \gdl code $\quo{(f,p)}$
is greater than $n$.
Given that there are an infinite number of real numbers,
we expect this search should always terminate for some step count $s$ --
a fact we can likely prove within $L'$
on the grounds that there are an infinite number of rationals,
each of which corresponds to a real number.

Building on the above computational search
for the next real number with \gdl code greater than $n$,
we can build a primitive-recursive function $R(i,s)$
that expresses the unbounded search for the $i$'th real number --
essentially by stepping through $i$ sequential executions
of $N$ to termination and returning
only the result of the $i$'th successive termination of $N$.
This primitive-recursive function thus represents a computation
that recursively enumerates all of the real numbers
meeting our realness criteria above,
if we allow $R$'s step count argument $s$ to grow without bound.
Further, we can prove not only in $L'$ but even in $L$
that given any particular real number $r$,
the enumeration function $R$ will find it:
since $R$ tests one potential real number per step,
the \gdl code for $r$ trivially upper bounds the step count required
for $R$ to find $r$.

\subsubsection{Cantor's diagonalization argument in $L'$}

We now use this recursive enumeration of the reals
to construct a new real number as in Cantor's famous diagonalization proof:
a real number $C$ that essentially takes a digit number $i$,
finds the $i$'th real number in the enumeration of real numbers,
and sets bindary digit $i$ of $C$
to the opposite of digit $i$ of real number $R_i$ of the enumeration.

We first form a primitive-recursive function $C(i,s)$
that uses the real-enumeration function $R$ above to step
until $R$ finds and produces the $i$'th real number $r = (f,p)$,
then uses metacircular evaluation
to execute the primitive-recursive function $f$ embedded in this real number
for as many steps as needed to determine the $i$'th digit of $r$.
Only when $C$'s step count argument $s$ is large enough
for this computation to terminate,
$C$ returns $C_i = 1 - r_i$,
the $i$'th digit of Cantor's real,
which is the complement of the $i$'th binary digit of $r$.

We must next form the corresponding \gdl-coded proof $p$ in $L'$
that $f$ always terminates
and satisfies the properties required of a Dedekind cut.
The fact that the rationals and hence reals are infinite,
and hence the enumeration above will always find a next one,
enables us to prove in $L'$ that this search always succeeds,
and we can use the $\mu$ operator in $L'$
to find the least step count $s$ at which it does so.\footnote{
	As an alternative to proving in $L'$
	that there is always a next real number,
	we could instead simply analyze in $L'$
	the hypothetical case of there at some point
	being \emph{no} next real number.
	At such a point,
	the $\mu$ operator returns 0
	indicating that the search for a next real number fails to terminate,
	implying that we have already found all of the real numbers.
	Since by this hypothetical chain of reasoning
	we have now exhausted the supply of real numbers
	that Cantor's real $C$ needs to differentiate itself from,
	we can simply let all the remaining digits $C_i$ of $C$ be zero --
	in which case Cantor's real becomes curiously rational,
	but the rest of our reasoning still works.
}
Given that $f$'s eventual termination is provable in $L'$,
proofs of the other properties of Dedekind cuts
represent straightforward arithmetic proofs in $L'$.

We thus have the required proof $p$ in $L'$
that $f$ behaves as required to define a Dedekind cut,
so $C = (f,p)$ is a real number by our definition.
We can verify this fact not only in $L'$ but even in $L$,
since checking that an $(f,p)$ represents a real number is primitive recursive.
As a result, reasoning in $L'$,
Cantor's real number $C$ exists and has a \gdl code
that must be in the enumeration $R$ of all real numbers expressed earlier.
Yet by $C$'s construction,
each digit $i$ of $C_i$ is distinct from digit $i$
of real $i$ in the enumeration,
so $C$ is also distinct from and unequal to each real number in the enumeration,
including from $C$ itself.
Since $L'$ proves that $C$ is both equal to and unequal to itself,
it follows that $L'$ is inconsistent.

\later{	Cantor's real is in the same boat as Berry's and Yablo's paradoxes:
	each express a computation that cannot terminate
	because it is trying to pull truth values down from infinity.
}

\subsection{Cantor's paradise lost?}

At the time of writing, the above is of course merely a sketch of a proof
that demands rigorous formal verification or refutation.\footnote{
	Samuel Gruetter has already translated
	key steps of this argument into a formal proof skeleton
	using the Coq proof assistant~\cite{huet16coq,chlipala13certified}
	(recently renamed \href{https://rocq-prover.org}{Rocq}),
	although this skeleton axiomizes large reasoning steps 
	yet to be formalized or verified.
	The author of this draft is in the midst
	of a parallel formalization effort
	using Isabelle/HOL~\cite{nipkow02isabelle}.
}
If this or similar line reasoning holds up, however,
then the implications would be significant.
Not only do we find that Zermelo-Fraenkel set theory is inconsistent,
along with similarly-powerful set-theoretic foundations
such as von Neumann–Bernays–Gödel set theory (NBG)
and Morse-Kelley (MK),
but so is Peano Arithmetic (PA),
which can use its first-order quantifiers
to express the $\mu$ function above
and assign a classical truth value to its result.
We may find many intuitionistic systems to be inconsistent as well --
particularly those closely-related to a corresponding classical system
and having the property that for any formula $f$ provable classically,
its double negation (`$\neg\neg f$') is provable intuitionistically.
This certainly appears to include Heyting Arithmetic (HA), for example.

In essence,
the signs seem to suggest that even just
\emph{reasoning about arbitrary recursive computations in general}
may represent a hard barrier that classical (and perhaps intuitionistic)
formal systems cannot breach without losing consistency.
Moreover, reconsidering \gdl's incompleteness theorems in this light
suggests that classical (and perhaps intuitionistic) systems
may never even be able to \emph{reach} the point
of expressing and reasoning about
arbitrary recursive (Turing-complete) computations
without losing consistency,
but instead can at best achieve arbitrarily-close approximations.
(We might start by considering extensions to PRA
that allow higher powers of recursion,
\eg allowing expression of Ackermann's functions,
while still retaining the same underlying principles
and in particular staying restricted to terminating computations.)

Applying \gdl's incompleteness results to systems
that apparently remain consistent
but are too week to express all Turing-complete computations
(Church's simply-typed lambda calculus being another example),
it appears that the traditional interpretation of \gdl's theorems
probably still applies as before:
\ie that these systems are consistent but incomplete
and unable to prove their own consistency,
only the consistency of strictly-weaker systems.
Applying \gdl's theorems to stronger systems beyond the ``computation barrier'',
however,
it may be that the traditional interpretation of these theorems was mistaken:
in these cases, \gdl's proofs may be not so much about incompleteness at all,
but rather instead are incomplete proofs of inconsistency.

Examining ``the Halting problem'' as expressed by Turing in this light,
it appears that this ``problem'' itself may be not so much
``an unsolvable problem'' but rather an \emph{ill-defined and meaningless} one.
To reason coherently about a ``problem'' being solvable or unsolvable,
we first need a coherent definition of the problem and a test for success:
\ie a way to check whether a purported solution indeed solves the problem,
or alternatively a way to prove categorically that no solution can exist.
But what if the ``test for success'' we wish to measure solutions against
is itself just a divergent, nonterminating computation,
to which we can never attach any meaningful boolean answer
without falling into logical inconsistency?
To attach a truth value to the \emph{question}
of whether a given program $P$ ``solves'' the Halting ``problem'',
it appears we must \emph{first} answer that same ``question'' --
in the case of Turing's halting-detection program, for example,
which first requires (again) having already solved the Halting problem,
\emph{ad infinitum}.
In this sense, the expression of the Halting problem itself --
even before and independently of any purported ``solution'' to it --
appears to be in the same boat as Berry's and Yablo's paradoxes,
not to mention Cantor's divergent real number.

Again assuming this or a similar line of reasoning holds up, therefore,
it appears fairly urgent to explore alternatives --
ideally to find a consistent foundation for reasoning
that can at least \emph{reach} the ability to express and reason about
arbitrary Turing-complete computations,
even if we perhaps must give up the hope of \emph{surpassing} that barrier
without losing consistency.
Grounded deduction appears to be an interesting candidate for this purpose,
for further exploration and development.
Another alternative of course is simply to ``give up'' on consistency,
and adopt the program of paraconsistency or dialetheism.\footnote{
	See for example \cite{priest06in}.
}
This approach has many theoretical and practical issues of its own,
however,
starting with the unanswered question of how robust to inconsistency
a paraconsistent system can be while still being \emph{usable}
for practical reasoning.

More immediately and pragmatically,
the loss of not only the law of excluded middle (LEM)
from powerful (\eg set-theoretic) systems,
but also apparently even the limited principle of omniscience (LPO),
would appear to deny us basic tools with which
we are accustomed to constructing and reasoning about
useful mathematical objects such as infinite sets and real numbers.
Faced with a choice between logical consistency
and having infinite sets and real numbers,
mathematicians may be understandably tempted --
even at risk of inconsistency --
to side with David Hilbert's famous sentiment:

\begin{quote}
``No one shall drive us out of the paradise which Cantor
has created for us.''\footnote{
	From a lecture given in Münster
	to Mathematical Society of Westphalia on June, 1925.
	See \cite{hilbert26uber} for the original German version,
	or \cite{hilbert26on} for an English translation.
}
\end{quote}

If we wish to consider grounded deduction 
as a potential alternative approach to reasoning in practice,
is there a viable way to regain in any fashion
at least some of the beautiful and useful mathematical abstractions
we found in Cantor's paradise,
such as infinite sets and real numbers to start with?

\section{Reflective idealization: Cantor's paradise regained?}
\label{sec:rga}

Since not only the law of excluded middle (LEM)
but even the limited principle of omniscience (LPO)
appear incomplatible with consistent reasoning in \gd,
to what other principle might we turn in seeking to express and build
anything recognizably like familiar abstractions
such as infinite sets and real numbers?

Informally,
all the known routes to inconsistency appear to stem from
ungrounded assumptions that we can ``safely'' ascribe truth values
to expressions that can represent self-referential or divergent
logical dependency structures.
We are either pulling truth values out of dependency cycles,
as in the Liar paradox,
or we are pulling truth values down from infinity,
as in Yablo's paradox.
But taking the liberty of doing so
always appears to get us into trouble eventually.
Is there \emph{any} context in which we might idealistically presume
to assign classical truth values
to the results of complex, open-ended processes
such as Turing-complete computations?
There may well not be.

If there remains any ``safe space''
for classical mathematical idealism,
however,
then one observation is that this safe space may be
precisely at the \emph{metalogical boundary}:
that is, not in a target logic itself at all
(at least serving in the role of target logic),
but rather ``just beyond'' the target logic's ``reach''
in whatever metalogic we are using to reason about it.
From the perspective of a target logic,
the metalogic being used to analyze it
is akin to some kind of deity:
the unknowable, beyond the borders of the target logic's universe,
and hence by construction impervious to whatever attempts might be made
\emph{within the target logic}
to cause trouble through circular or divergent dependency structures.
Stated simply, nothing expressible in a target logic
can ``know'' or depend on anything about a metalogic in use to analyize it.
From another perspective,
we can view a target logic as merely a simulation
``running'' in some metalogic;
the metalogic can potentially know and analyze anything about the target logic,
while the target logic cannot know or depend on anything in the metalogic.
The potential dependencies run only one way.

This principle appears true regardless of our choices of target and metalogic.
The two might be different logics,
or they might happen to be (instances of) the same logic
playing the distinct \emph{roles} of target logic and metalogic.
If we use \ga to reason metalogically about \ga,
then we can see and know that we are doing so
only from the ``outer'' perspective of \ga in the role of metalogic.
The inner target logic, in contrast,
cannot ``know'' that it is being modeled and analyzed
by another identical instance of \ga:
the target logic's universe might just as well be
modeled in some other entirely-different metalogic,
a simulation in the mind of some other deity.

\subsection{Reflective excluded middle (REM)}

These observations lead us to wonder whether
some \emph{reflective} form of mathematical idealization
might be safely added to \gd or \ga without causing inconsistency.
Concretely, suppose we have built within ordinary, constructive \ga (\cga)
all the \gdl-style tooling for reflective reasoning,
which we already needed to do anyway
to get the computable quantifiers as discussed in \cref{sec:comp}.
Let `$\quo{\Gamma \vdash c}$' be the \gdl code
for the logical entailment or judgment `$\Gamma \vdash c$'
modeled reflectively within \ga to reason about \ga itself.
Let `$\fpr\quo{\Gamma \vdash c}$' be a predicate
in the outer instance of \ga, serving in the role of metalogic,
stating that judgment `$\Gamma \vdash c$' is provable in \ga.
Then let us consider extending \cga
with something like the following axiom, or a corresponding inference rule,
expressing a principle of \emph{reflective excluded middle} or REM:

\[
	\vdash	\fpr\quo{\ctrue \vdash p} \lor 
		\fpr\quo{p \vdash \cfalse}
	\tag{REM}
\]

Informally,
this principle makes an assertion not about \emph{truth}
but rather about \emph{provability}:
namely, that for any expressible proposition $p$,
either $p$ is provably true,
or \emph{hypothetically supposing} $p$ to be true 
provably leads to contradiction.
In a classical logic,
this principle would be equivalent to the LEM --
but in grounded reasoning this equivalence does not appear to hold
(and clearly \emph{must not} hold if \ga is consistent).

Notice that this principle applies not just to \emph{truth}
but also to \emph{booleanness} in \ga.
We can for example substitute `$p \jbool$' for `$p$' in the REM above
to obtain trivially:
\[
	\vdash	\fpr\quo{\ctrue \vdash p \jbool} \lor 
		\fpr\quo{p \jbool \vdash \cfalse}
\]

Instantiated in this way,
the REM effectively states that for any expressible proposition $p$,
either $p$ provably has a truth value (is boolean),
or else \emph{hypothetically supposing} $p$ to have either truth value
would lead us provably to contradiction.
The former case represents the ``normal'' case
where $p$ has a classical truth value,
while the latter represents the ``gappy'' case
where $p$ fails to have any truth value.

Note that if `$\fpr\quo{\ctrue \vdash p \jbool}$'
is metalogically true,
this does \emph{not} necessarily mean
that `$\fpr\quo{p \jbool \vdash \cfalse}$'
is metalogically false.
The former might be metalogically true
while the latter \emph{has no truth value},
even metalogically,
as we will see in examples below.
This possibility is crucial in avoiding the morass of
``revenge problems'' that attempting to force truth values on expressions
invariably gets us into, sooner or later.\footnote{
	See for example \cite{beall08revenge}.
}

\subsubsection{A few paradoxes in the light of REM}

Consider the Liar paradox `$L \ldef \neg L$' as a first illustrative example.
Assuming \ga is consistent,
then we must not be able to assign either truth value to $L$;
therefore the metalogical statement `$\fpr\quo{\ctrue \vdash p \jbool}$'
must be false.
Even without the REM,
we can easily use \gdl-style reflective reasoning in \ga about \ga
to prove that the contrary metalogical statement,
`$\fpr\quo{p \jbool \vdash \cfalse}$', is metalogically true.
(Ascribing either truth value to $p$ leads us immediately to the other,
a contradiction, and hence to \cfalse.)
So the Liar paradox appears entirely consistent
with the claim that the REM makes,
although we did not actually need the REM in this case.

Consider the Truthteller, in contrast, `$T \ldef T$'.
This definition is what Kripke calls
\emph{ungrounded} but not \emph{paradoxical}.
Grounded reasoning gives us no basis to assign a truth value to $T$,
but if we were to ``force'' it to have a value by an axiom for example,
then doing so would not immediately appear to cause paradox.
The Truthteller merely agrees with any truth value it already has.
Forcing a truth value on $T$ would change the system, however.
In the \emph{unmodified} \ga system,
we find no grounds to assign $T$ any value.
We can reason reflectively in \ga to prove that this is the case,
which means that `$\fpr\quo{\ctrue \vdash T \jbool}$' is metalogically false.
But if we can reflectively prove `$\neg \fpr\quo{\ctrue \vdash T \jbool}$',
then hypothetically ascribing either truth value to $T$
would contradict exactly that metalogical claim.
Thus, even without using the REM,
we may be able to use \gdl-style reflective reasoning
to prove that `$\fpr\quo{T \jbool \vdash \cfalse}$' is metalogically true.

Similar reasoning appears to hold for trickier examples
such as the Strengthened Liar,
which we may express informally as
``This sentence is not true'',
which in \ga we might formalize as
`$L' \ldef \neg \fpr\quo{\ctrue \vdash L'}$'.
That is, the Strengthened Liar asserts that
the sentence itself is not (provably) true,
while apparently leaving open the possibility of it being ``true''
\emph{by virtue of} the sentence being (apparently) not provably true.
But this is necessarily a problem only in classical reasoning,
where we presume that there must be some ``fact of the matter''
even regarding whether $L'$ has a truth value:
in grounded reasoning the situation is different.
In constructive \ga without the REM at least,
the Strengthened Liar is just another ungrounded statement:
to have a truth value, 
we would first have to prove
that the metalogical claim `$\fpr\quo{\ctrue \vdash L'}$'
has some truth value,
but our attempt to do \emph{that} in turn
leads us back to the need first to prove something about $L'$.
Nevertheless, ``stepping back'' into a further,
outer level of metalogical analysis,
we find that hypothetically ascribing any truth value to $L'$
would lead us to contradiction,
so the claim `$\fpr\quo{L' \vdash \cfalse}$'
is in fact metalogically true.
Once again, we don't necessarily need the REM
to make this metalogical judgment.

Can we find a paradoxical (or ungrounded but non-paradoxical) statement in \ga
that, even given a sufficient number of ``steps back''
into outer metalogical levels,
cannot -- without using the REM -- already be ``binned''
into one of the two metalogical categories of ``$p$ has a truth value''
or ``supposing $p$ has a truth value leads us to a contradiction
at \emph{some} level of metalogical analysis''?
This appears to be a non-obvious question for further exploration.
But perhaps this non-obviousness actually bodes well for the REM's chances.
Perhaps we don't need it.  Perhaps we do.
Perhaps, without first ``taking the plunge'' and adopting the REM,
there is nothing we can say -- \emph{no fact of the matter} --
about whether we do or do not need it.

But in any case, is REM or something like it
conceivably at least \emph{consistent} to adopt?
Classical metalogical reasoning about \ga suggests that
if \cga is already consistent,
then adding REM to it cannot make it inconsistent.
If \cga without REM is consistent,
then -- reasoning classically --
for any expressible proposition $p$,
either `$\fpr\quo{\ctrue \vdash p}$' might be true,
or `$\fpr\quo{p \vdash \cfalse}$' might be true,
or \emph{neither} of these might be true --
but \emph{both} cannot be true by our consistency assumption.
If either `$\fpr\quo{\ctrue \vdash p}$'
or `$\fpr\quo{p \vdash \cfalse}$' is already provably true for a given $p$,
then adding the REM does not change the situation:
the REM merely ``agrees with'' the already-established fact
that there is a metalogical truth one way or the other concerning $p$.
If neither `$\fpr\quo{\ctrue \vdash p}$'
nor `$\fpr\quo{p \vdash \cfalse}$' is already provably true in \ga,
then from the perspective of the classical metalogic we are reasoning in,
$p$ clearly has no truth value provable in \ga associated with it.
In this case,
let us \emph{pretend} that assigning $p$ a truth value
would lead to contradiction,
essentially by ``forcing'' the metalogical statement 
`$\fpr\quo{p \vdash \cfalse}$' to true.
This claim, even if forced,
appears to have no opposite metalogical truth
that it could contradict.
This reasoning is only preliminary and informal,
and we might be rightfully suspicious
of classical reasoning about \ga in any case,
but it is a starting point.

\subsection{Cantor's paradise regained $\dots$ maybe?}

While leaving detailed development for later,
it appears that the REM or something like it
may offer a means to rebuild something recognizably like
infinite sets or the real numbers in \ga,
even if their properties would necessarily be different
and perhaps weaker and more nuanced
than in the classical world we are accustomed to.

For example,
building on the reflective ``proof-carrying code'' techniques
used above in \cref{sec:pra},
we might define real numbers in \ga as (\gdl-coded)
computations to test whether a rational number
is in the upper or lower part of a Dedekind cut,
together with a (\gdl-coded) proof in \ga 
that this computation indeed reliably decides given any rational number.
We cannot expect to construct Cantor's real number in \ga, however,
because it expresses a divergent and hence non-terminating computation --
reinvoking ``itself'' as a precondition for computing ``its own'' digit,
wherever that may be --
and \ga's \emph{habeas quid} preconditions will thus
require us to to satisfy an infinitely-growing chain of dependencies
before we could conceivably even get started
on assigning a truth value to Cantor's real.
Since hypothetically presuming Cantor's computation yields a truth value
would lead to a contradiction, however,
a computational enumeration of the real numbers in \ga
will not only \emph{not} find Cantor's real in its enumeration of reals,
but will also find the metalogical fact-of-the-matter
that Cantor's real \emph{cannot be} a real number,
precisely because hypothetically assuming it was would lead to contradiction.

Similar techniques appear to make it conceivably possible
to construct grounded infinite sets
with some recognizable similarities to classical sets,
with important conditions and caveats.
Set comprehension will clearly require a key \emph{habeas quid} precondition:
essentially a proof that a predicate \emph{always decides} the set's membership
before we can deduce that the predicate defines a set.
These \emph{habeas quid} requirements appear to head off Russell's paradox
and various others,
again by placing circular proof-obligation barriers
on the path to proving any such ungrounded, indecisive constructions are sets.
Given particular existing sets,
it appears possible to satisfy these \emph{habeas quid} requirements
in order to build ``larger'' sets,
roughly in Zermelo-Fraenkel fashion,
though with important new limitations likely appearing.

In the resulting system it may be, for example, that 
there is not and cannot be a ``set of all real numbers''
corresponding to the familiar classical set $\sreal$.
Instead, we may be able to construct only \emph{particular subsets}
of all real numbers -- at least in terms of \emph{sets}
that always decide their membership.
A set of all rational numbers, probably.
A set of all algebraic numbers, maybe.
A set of all real numbers, less likely.
A set of all reals might become just as meaningless and impossible
as a set of all sets or a largest ordinal.
To express the concept of ``all real numbers'' coherently at all
we may have to think purely in terms of computations,
which may or may not terminate,
and hence which may or may not have meaningful results at all.
To form a set we must choose a predicate that always decides its membership,
but such decision predicates can at best recognize
particular kinds of terminating computations,
not all terminating computations.
We cannot meaningfully form ``sets''
with their familiar classical properties such as extensionality
if they cannot decide their membership;
we can talk about computations and perhaps enumerations of computations,
but performing such a recursive enumeration
will not necessarily (often perhaps cannot, consistently) form a set at all.

%com

\later{
\input{real}

\input{turing}

\input{model}
\input{cons}	% or even before refl?
\input{unprov}
\input{axiom}

%\part{For Later}

\input{fun}
\input{set}

%\later{
\input{meta}
%}%later

\input{other}	% proving consistency of other systems in \gd
\input{rich}
\input{semantic}	% applying \gd to the semantic paradoxes of self-reference
\input{pa}
\input{rel}
\input{concl}
}%later

\section{Non-conclusion: ``to be continued$\dots$''}

As mentioned in the introduction,
this is a preliminary and incomplete draft,
in particular a portion of a longer document
to be released progressively.
There is no conclusion yet.
Feedback is highly welcome and appreciated
on the portions released so far.

\subsection*{Acknowledgments}

The author wishes to thank colleagues
who have offered helpful feedback so far,
in particular:
Jeff Allen,
Cristina Basescu,
Elliot Bobrow,
Thomas Bourgeat,
Anton Burtsev,
Adam Chlipala,
Henry Corrigan-Gibbs,
Samuel Gruetter,
M. Frans Kaashoek,
Sascha Kehrli,
Viktor Kuncak,
Derek Leung,
Stefan Milenković,
John Regehr,
Robbert Van Renesse,
Deian Stefan,
and
Nickolai Zeldovich.

\subsection*{Mystery hash}

This digest covers the document's sources
including the parts not yet released:

\begin{tiny}
\begin{center}
42829743fec72377206e62d87649856aaa7f49e4e810fabedb6766c2fda27560f668de6eae6fcd146c23e3bc20218e3b3d818d1a9e0d3a80523371b39a5f90c6

\end{center}
\end{tiny}

\bibliographystyle{apalike}	% bibtex
\bibliography{logic}
%\printbibliography		% biblatex

\end{document}